\documentclass[a4paper,11pt]{article}
\usepackage{amsmath,amssymb,latexsym, enumerate}
\usepackage{amsthm,amsmath,bm}
\usepackage{color}
\usepackage{graphicx}
\newtheorem{thm}{Theorem}[section]  \newtheorem{cor}[thm]{Corollary}
\newtheorem{lem}[thm]{Lemma}   \newtheorem{defn}[thm]{Definition}
  \newtheorem{prop}[thm]{Proposition}

\def\remark{\refstepcounter{thm}\bigskip
         \noindent\bf Remark \thethm.\rm\ }
\newcommand{\preuve}[1][\!\!]{\bigskip
            \noindent{\bf Proof #1. \ \ }}
\def\fin{\hfill$\Box$\\}

\newcounter{numexo}\newcounter{numq}
\newcounter{numsq}

 \def\kkk{{\cal K}}\def\lll{{\cal L}}
\def\mmm{{\cal M}}\def\nnn{{\cal N}} \def\ooo{{\cal O}}\def\ppp{{\cal P}}
\def\sss{{\cal S}}

\def\R{\mathbb R}
\def\N{\mathbb N}

\def\D{\partial}\def\eps{\varepsilon}\def\phi{\varphi}
\def\norm#1{\big\Vert#1\big\Vert}

\def\abs#1{\left\vert#1\right\vert}
\def\set#1{\left\{#1\right\}}
\def\seq#1{\left<#1\right>}
\def\sep#1{\left(#1\right)}
 
\def\defegal{\stackrel{\text{\rm def}}{=}}

\renewcommand{\bar}[1]{\overline{#1}}

\newcommand{\un}{\textcircled{\tiny{1}}}

\def\b{{\rm{b}}}
\def\tb{{\tilde{b}}}

\newcommand{\di}{\,d}
\newcommand{\normm}[1]{| \! | \! | #1| \! | \! | }
\newcommand{\bigset}[1]{\big\{#1\big\}}
\newcommand{\inner}[1]{\left(#1\right)}
\newcommand{\biginner}[1]{\Bigl(#1\Bigr)}
\newcommand{\comi}[1]{\bigl(#1\bigr)}
\newcommand{\comii}[1]{\left<#1\right>}
\newcommand{\com}[1]{\bigl[#1\bigr]}
\newcommand{\reff}[1]{(\ref{#1})}

\definecolor{cyan}{cmyk}{1,0,0,0}
\definecolor{lightcyan}{cmyk}{0.5,0,0,0}
\definecolor{pastelcyan}{cmyk}{0.25,0,0,0}
\definecolor{magenta}{cmyk}{0,1,0,0}
\definecolor{yellow}{cmyk}{0,0,1,0}
\definecolor{lightyellow}{cmyk}{0,0,0.5,0}
\definecolor{pastelyellow}{cmyk}{0,0,0.25,0}
\definecolor{black}{cmyk}{0,0,0,1}
\definecolor{darkgray}{cmyk}{0,0,0,0.75}
\definecolor{gray}{cmyk}{0,0,0,0.5}
\definecolor{lightgray}{cmyk}{0,0,0,0.25}
\definecolor{white}{cmyk}{0,0,0,0}
\definecolor{red}{cmyk}{0,1,1,0}
\definecolor{orange}{cmyk}{0,0.5,1,0}
\definecolor{scarlet}{cmyk}{0,1,0.5,0}
\definecolor{brown}{cmyk}{0.5,0.75,1,0}
\definecolor{camel}{cmyk}{0.25,0.375,0.5,0}
\definecolor{cream}{cmyk}{0,0.2,0.3,0}
\definecolor{green}{cmyk}{1,0,1,0}
\definecolor{lightgreen}{cmyk}{0.5,0,0.5,0}
\definecolor{pastelgreen}{cmyk}{0.25,0,0.25,0}
\definecolor{mossgreen}{cmyk}{0.64,0.4,1,0}
\definecolor{yellowgreen}{cmyk}{0.5,0,1,0}
\definecolor{skyblue}{cmyk}{0.4,0.16,0,0}
\definecolor{royal}{cmyk}{1.0,0.5,0,0}
\definecolor{navyblue}{cmyk}{0.9,0.75,0.5,0}
\definecolor{lightnavy}{cmyk}{0.4,0.3,0.2,0}
\definecolor{blue}{cmyk}{1,1,0,0}
\definecolor{lightblue}{cmyk}{0.5,0.5,0,0}
\definecolor{lavender}{cmyk}{0.25,0.25,0,0}
\definecolor{violet}{cmyk}{0.75,1,0.25,0}
\definecolor{purple}{cmyk}{0.5,1,0.5,0}
\definecolor{pink}{cmyk}{0,0.5,0,0}
\definecolor{pastelpink}{cmyk}{0,0.25,0,0}
\def\black{\color{black}}
\def\red{\color{red}}

\def\red{\color{black}}

\renewcommand{\un}{1 \!\! 1}

\title{Global hypoelliptic and symbolic estimates for the linearized Boltzmann operator without angular cutoff}
\newcommand{\footnoteremember}[2]{ \footnote{#2}\newcounter{#1}\setcounter{#1}{\value{footnote}}}

\def\nantes{Laboratoire de Math\'{e}matiques Jean Leray,
Universit\'{e} de Nantes, 44322 Nantes, France}
\def\wuhan{School of Mathematics and Statistics, and Computational Science Hubei Key Laboratory, Wuhan University, 430072
Wuhan, China}

\def\brest{Arts et M\'{e}tiers ParisTech, Paris 75013 France}
\author{Radjesvarane Alexandre\footnote{\brest} 
\and Fr\'{e}d\'{e}ric H\'{e}rau\footnoteremember{nan}{\nantes}
\and Wei-Xi Li \footnote{\wuhan}
\footnote{e-mail: frederic.herau@univ-nantes.fr, wei-xi.li@whu.edu.cn, radjesvarane.alexandre@paristech.fr}
\footnote{WXL   was supported in part by NSFC(11422106) and Fok Ying Tung Education Foundation (151001), WXL and FH were supported by ANR NOSEVOL (2011 BS01019 01) }}
\date{}

\begin{document}

\maketitle

\begin{abstract}
In this article we provide  global subelliptic estimates for the linearized inhomogeneous Boltzmann equation without angular cutoff,
and show that  some global gain in the spatial
direction is available although the corresponding operator is not elliptic in this direction.
The proof is based on a multiplier method and the so-called Wick quantization, together with a careful analysis of the symbolic properties of the Weyl symbol of the Boltzmann collision operator. \smallskip



\end{abstract}

\tableofcontents

\section{Introduction}

In this paper we are interested in giving sharp subellipitic estimates for the non-homogeneous
linearized Boltzmann operator
$$
\ppp = v \cdot \D_x -\lll
$$
considered as an unbounded operator in $L^2(\R_x^3 \times \R^3_v)$,
where $\lll$ is the linearized Boltzmann without cutoff collision kernel whose precise
expression is given in \eqref{defdeL} in the next  subsection. Here $x$ in $ \R_x^{3}$ and $v $ in  $\R^3_v$ are respectively the space and velocity variable and  $\D_x$ denotes the gradient in the space variable. The main result of this paper is the sharp estimate given in Theorem \ref{thmain}. In this introduction we first present the model, then the main results including Theorem \ref{thmain} and   bibliographic comments and we conclude by giving some general comments about the interest of this work and  the methodology we followed for the proofs.

\subsection{Model and notations} \label{modnot}

Let us first recall some facts about the
 non-cutoff inhomogeneous Boltzmann equation. It reads
\begin{equation} \label{QL}
     \D_{t}F+v \cdot \D_{x}F =  Q(F,F),
\end{equation}
with $F$ standing for a probability density function, and a given Cauchy data at  $t=0$, while the position $x$ and velocity $v$ are in $\R^3$, see \cite{CerIllPul94,Vil01} and references therein for more details on Boltzmann equation.
In \eqref{QL}, the collision kernel $Q$ is defined for sufficiently smooth functions $F$ and $G$ by
\begin{equation*} 
Q(G,F)(t,x,v) =  \int_{\R^3}\int_{\mathbb S^{2}} { B}(v-v_*, \sigma) \sep{ F'  G'_* - F G_*} dv_* d\sigma
\end{equation*}
where $F'=F(t,x, v')$, $F= F(t,x, v)$, $G'_*=G(t,x,v'_*)$ and $G_*= G(t,x,v_*)$ for short. For given velocities after (or before) collision $v$ and $v_*$,   $v'$ and  $v_*'$ are the velocities before (or after) collision, with the following energy and momentum conservation rules, expressing the fact that we consider elastic collisions
 \begin{equation} \label{kmr}
 v'+v_*' = v + v_*, \ \ \ |v'|^2+|v_*'|^2 = |v|^2 + |v_*|^2.
 \end{equation}
 where $|v| $ denotes the canonical euclidian norm in $\R^3$.
We will choose the so-called $\sigma$ representation, for $\sigma $ on the sphere $\mathbb S^{2}$,
 $$
 \left\{
 \begin{array}{l}
   v'  = \frac{v+v_*}{2} + \frac{ |v-v_*|}{2} \sigma \\
  v'_* = \frac{v+v_*}{2} - \frac{ |v-v_*|}{2} \sigma ,
 \end{array}
 \right.
 $$
and define the deviation angle $\theta$ in a standard way by
 $$
 \cos \theta =  \frac{ v-v_*}{|v-v_*|} \cdot  \sigma,
 $$
 where $ \cdot $ denotes the usual scalar product in $\R^3$.
 In the case of inverse power laws, see for example \cite{CerIllPul94}, the collisional cross section $ B$ looks approximatively as follows
 \begin{equation} \label{cross-section}
 { B}(v-v_*, \sigma) =  |v-v_*|^\gamma \b (\cos \theta),
 \end{equation}
 for some real parameter $\gamma$ and some function $b$.

 Without loss of generality, we assume $B(v-v_*,\sigma)$ is supported on the set $(v-v_*) \cdot \sigma \geq 0$ which corresponds to $\theta \in  [0,\pi/2]$, since as usual, see \cite{AV02},
 $B$ can be eventually replaced by its symetrized version
 $$
 \overline{B} (v-v_*, \sigma)= {B} (v-v_*, \sigma) + {B} (v-v_*, -\sigma).
  $$
Moreover, we assume that we deal with inverse power interaction laws between particles, and thus according to \cite{CerIllPul94}, we assume that $b$ has the following singular behavior
 when $\theta \in ]0, \pi/2[$ : there exist a constant $c_b>0$ such that
 \begin{equation*} 
c_b^{-1} \theta^{-1-2s} \leq \sin \theta \b (\cos \theta)  \leq c_b  \theta^{-1-2s}, \textrm{ as } \theta \longrightarrow 0^+.
 \end{equation*}
  In the preceding formulas, we will impose the following range of parameters, coming from the
 physical derivation,
 $$
  s \in (0,1), \ \ \ \ \gamma \in ( -3, \infty )   .
 $$
Note that the last condition on $\gamma +2s$ is weaker than in \cite{AMUXY1,gres-strain} since we will deal only with the linearized part of Boltzmann collisional operator.
%
%


The behavior of this singular kernel is strongly related to the following non-integrability condition
$$
\int_0^{\pi/2} \sin\theta b(\cos\theta) d\theta = \infty,
$$
which implies some diffusion properties of the (linearized) Boltzmann operator that we will explain more in  detail later.

In some expressions involving the integral kernels, it may therefore happen that some non-integrability
arise, and in this case these integrals have to be understood as principal values (see the appendix or  \cite{AV02}).  Anyway we shall do most of the computations as if $B$ were integrable and use the principal value trick whenever needed.

In this work, we are interested in the linearized Boltzmann  operator, around a normalized Maxwellian distribution, which is described as follows. Let this normalized Maxwellian be
\[
   \mu(v)=\inner{2\pi}^{-3/2}e^{-\abs{v}^2/2}.
\]
Setting $F=\mu+\sqrt{\mu}f$,  the perturbation $f$ satisfies the equation
\[
    \partial_t  f+v\cdot\partial_x f-\mu^{-1/2}Q(\mu,~\sqrt{\mu}f)-\mu^{-1/2}Q(\sqrt{\mu}f,~\mu)=
    \mu^{-1/2}Q(\sqrt{\mu}f,~\sqrt{\mu}f),
\]
since $\partial_t  F+v\cdot\partial_x F-Q(F,~F)=0$ and $Q(\mu,~\mu)=0$. Using the notation
\[
  \tilde \Gamma(g,~f)=\mu^{-1/2}Q(\sqrt{\mu}g,~\sqrt{\mu}f),
\]
we may rewrite the above equation as
\begin{equation*}  
 \D_t f +  {\mathcal P} f=\tilde \Gamma(f,~f),
\end{equation*}
where the linearized Boltzmann operator $\mathcal P$ takes the form
\begin{equation} \label{mpintro}
   {\mathcal P}=v\cdot\partial_x-\mathcal L
\end{equation}
with
\begin{equation} \label{defdeL}
\lll = \lll_1 + \lll_2, \qquad   \mathcal L_1 f=\tilde \Gamma(\sqrt{\mu},~f),\quad \mathcal L_2 f=\tilde \Gamma(f,~\sqrt{\mu}).
\end{equation}
The   operator  $\mathcal P$ acts only in variables $(x,v)$,  is non selfadjoint, and consists of  a transport part which is skew-adjoint, and  a diffusion  part acting only in the $v$ variable.

The elliptic properties of this operator which is the autonomous linear part of  the Boltzmann equation are the main subject of this work and we present them below.

 \bigskip
 \subsubsection{Notations}
 Throughout the paper we shall adopt the following notations: we work in dimension $d=3$ and denote by $(x,v) \in  \R_{x}^{3}\times \R_{v}^{3}$ the space-velocity variables. For $v \in \R^3$  we denote  $\seq{v}=(1+\abs{ v}^2)^{1/2}$, where we recall that $|v|$ is the canonical Euclidian norm of $v$ in $\R^3$.

 The gradient in velocity (resp. space) will be denoted by $\D_v$ (resp. $\D_x$). We shall also denote  $D_v=\frac{1}{i}\partial_v$, $D_x=\frac{1}{i}\partial_x$, and denote $\xi$ the dual variable of $x$ and $\eta$ the dual variable of $v$.

 We shall extensively use the pseudo-differential theory, for which we refer to the appendix here and the reference therein. In particular operators $\seq{D_v}$ and   $\seq{v \wedge D_v}^{2s}$ denotes respectively the pseudo-differential operator with  classical  symbol $\seq{\eta}$ and  $\seq{v \wedge \eta}^{2s}$.

 We will work througout the paper in $L^2(\mathbb R_{v}^{3})$ or $L^2\inner{\mathbb R_{x}^{3}\times \mathbb R_{v}^{3} }$ for which we denote (without ambiguity depending on the sections) the scalar product by $( \cdot, \cdot)$ and the norm by $\norm{\cdot}$. We shall
 mainly work with functions in the Schwartz spaces  $\mathcal S(\mathbb R_v^3)$ or $\mathcal S(\mathbb R_x^3 \times \R_v^3)$.

 In all the article, the notation $a \approx b $ (resp. $a \lesssim b $) for $a$ and $b$ positive real means that there is some positive constant $C$ not depending on possible  free parameters such that $C^{-1} a \leq b \leq Ca$ (resp. $a \leq Cb$).


\subsection{Main results and bibliographic comments}

The main theorem of this paper deals with operator $\mathcal P$, viewed as an unbounded operator in $L^2\inner{\mathbb R_{x}^{3}\times \mathbb R_{v}^{3} }$. We adopt the conventions of notation given at the end of subsection \ref{modnot}

\begin{thm} \label{thmain}
For all $l\in \R$,
there exists a constant $C_l$ such  that for all $f \in \sss(\R_x^3\times \R_v^3)$, we have
 \begin{equation*}
 \begin{split}
&   \ \big\Vert \seq{v}^{\gamma} \seq{D_v}^{2s} f \big\Vert
  +  \big\Vert \seq{v}^{\gamma} \seq{v \wedge D_v}^{2s} f \big\Vert + \norm{\langle v \rangle^{\gamma +2s} f}    \\
&  \ \ \ \ \ \ \ \ \ \ + \big\Vert \seq{v}^{\gamma/(2s+1)} \seq{D_x}^{2s/(2s+1)} f \big\Vert
   +
  \big\Vert \seq{v}^{\gamma/(2s+1)} \seq{v \wedge D_x}^{2s/(2s+1)} f \big\Vert \\
&  \ \ \ \ \ \ \ \ \ \ \ \ \ \ \ \ \ \ \ \   \leq C_l \sep{ \norm{{\mathcal P}f} + \norm{\seq{v}^l f}}.
\end{split}
\end{equation*}
 Recall $\norm{\cdot}$ here stands for the norm $\norm{\cdot}_{L^2(\mathbb R^6)}$ in $L^2(\mathbb R^6).$
\end{thm}

Note carefully that we do not need to take into account the finite dimensional kernel associated with the linearized   Boltzmann operator \cite{AMUXY1,gres-strain} which is hidden again in the term  $\norm{\seq{v}^l f}$.

As an intermediate result, we are also able to give an explicit form of the so-called
triple norm introduced in \cite{AMUXY1}. Previous estimates from below were also given
  in \cite{mou} and \cite{mou1}, but the following coercivity estimate measures now  explicitly
  the global  weights and regularity gains of   the diffusion kernel $\lll$. Note that
  we again forget in the following result the fact that there is
  finite dimensional operator kernel.

  \begin{thm} \label{coermou}
  For all $l \in \R$, there exists a constant $C_l$ such  that for all $f \in \sss(\R_x^3\times \R_v^3)$, we have
\begin{multline*}
C^{-1}_l \sep{ \norm{ \seq{v}^{\gamma/2 } \seq{D_v}^{s}f }^2
 + \norm{\seq{v}^{\gamma/2} \seq{v \wedge D_v}^{s}f}^2
 + \norm{\seq{v}^{\gamma/2 + s}f}^2 } \\
  \leq -\sep{ \lll f, f} + \norm{\seq{v}^l f}^2 \\
  \leq C_l \sep{ \norm{ \seq{v}^{\gamma/2 } \seq{D_v}^{s}f }^2
 + \norm{\seq{v}^{\gamma/2} \seq{v \wedge D_v}^{s}f}^2
 + \norm{\seq{v}^{\gamma/2 + s}f}^2 }.
\end{multline*}
  \end{thm}

Theorem \ref{thmain}  can be extended to a time dependent version as follows, by considering the time dependent operator
$$
\widetilde{P} = \D_t+ v\cdot\partial_x-\mathcal L,
$$
the
functional spaces being now $L^2\inner{\mathbb R_t \times \mathbb R_{x}^{3}
  \times \mathbb R_{v}^{3}}$ 
  with norm  denoted by $\norm{.}_{L^2(\mathbb R^7)}$. With this setting, one can show that

\begin{thm}\label{th1}
For all $l\in \R$, there exists a constant $C_l$ such  that for all $f \in \sss\inner{\mathbb R_t \times  \mathbb
  R_{x}^{3} \times \mathbb R_v^3  }$, we have
 \begin{equation*}
 \begin{split}
&   \norm{\comii{v}^{\frac{\gamma-2s}{1+2s}}\seq{D_t}^{\frac{2s}{1+2s}}f}_{L^2(\mathbb R^7)} +\ \big\Vert \seq{v}^{\gamma} \seq{D_v}^{2s} f \big\Vert^2_{L^2(\mathbb R^7)}
  +  \big\Vert \seq{v}^{\gamma} \seq{v \wedge D_v}^{2s} f \big\Vert^2_{L^2(\mathbb R^7)}    \\
&  \ \ \ \ \ \ \ \ \ \ + \norm{\langle v \rangle^{\gamma +2s} f}^2_{L^2(\mathbb R^7)}   + \big\Vert \seq{v}^{\gamma/(2s+1)} \seq{D_x}^{2s/(2s+1)} f \big\Vert^2_{L^2(\mathbb R^7)}
    \\
&  \ \ \ \ \ \ \ \ \ \ \ \ \
 +
  \big\Vert \seq{v}^{\gamma/(2s+1)} \seq{v \wedge D_x}^{2s/(2s+1)} f \big\Vert^2_{L^2(\mathbb R^7)} \\
&  \ \ \ \ \ \ \ \ \ \ \ \ \ \ \ \ \ \ \ \   \leq C_l \sep{ \norm{\widetilde{P}f}^2_{L^2(\mathbb R^7)} + \norm{\seq{v}^l f}^2_{L^2(\mathbb R^7)}}
\end{split}
\end{equation*}
\end{thm}

The preceding results are consequences of fundamental pseudo-differential properties of the linearized
Boltzmann operator. Indeed, as we shall see in Section 3,  the operator
  $\lll = \lll_1 + \lll_2$  can be
 splitted as
 $$
 \lll_1 = -a^w -\kkk_1, \ \ \ \ \lll_2 =  -\kkk_2
 $$
 where $a\geq 0$ is real, its Weyl quantization $ a^w$ being a  pseudo-differential operator of order $2s$,
 and $\kkk = \kkk_1 + \kkk_2$ is controlled by  $a^w$ (see Proposition \ref{estaa} below and the review about Weyl-H\"ormander calculus in the appendix,  and we refer to \cite[Chapter 18]{Hormander85} and \cite{MR2599384} for  more detail on   Weyl-H\"ormander calculus).  Precise
  expressions of $a$ and $\kkk_i$ will be given in Section 3.
The most significant part  of $\lll$ is therefore of a pseudo-differential type and by the next result, we have fundamental symbolic estimates for $a$, implying in particular that operator $a^w$ is elliptic   in the symbolic class $S(a,  \abs{dv}^2 + \abs{d \eta}^2).$
   This very strong property allows to avoid the systematic use of G\aa rding type inequalities which are not available here.

In the following, we denote $\Gamma = \abs{dv}^2 + \abs{d \eta}^2$ is the flat metric in $\R^6_{v,\eta}$ (recall that $\eta $ denotes the dual variable of $v$).  Standard notions concerning symbolic estimates and the pseudo-differential calculus are explained at the beginning of section 4.

 \begin{prop} \label{estaa} Define 
  $$
  \tilde{a}(v,\eta) \defegal \seq{v}^{\gamma} (1+ |\eta|^{2} + |\eta\wedge v|^{2} +\abs{v}^{2} )^s, \mbox{ for all } (v,\eta) \in \R^6_{v,\eta}.
  $$
 Then  we can write $\lll = -a^w - \kkk$, where
\begin{enumerate}[i)]
\item
$\tilde{a}$  is admissible weight  for $\Gamma$  and
$  a, \tilde a \in S( \tilde{a} , \Gamma)  $,
and there exists a  positive constant $C$ such that
 $C^{-1} \tilde{a}(v,\eta) \leq a(v, \eta) \leq  C \tilde{a}(v,\eta)$;
\item for all $\eps > 0$ there exists $C_\eps$ such that
$$
  \norm{\kkk f} \leq \eps\norm{a^w f } + C_\eps \norm{\seq{v}^{\gamma+ 2s} f};
$$
\item for a sufficiently large constant $K$ depending only on
the dimension,  we define $a_K$   by
\begin{equation}
\label{deak}
a_K \defegal a+K\comii v^{\gamma+2s}.	
\end{equation}
Then $a_K$
belongs to $S( \tilde{a} , \Gamma)$,  is invertible as an operator in $L^2$ and  its inverse $\inner{a_K^w}^{-1}$  has the
  form
\[
    \inner{a_K^w}^{-1} =H_1 \inner{a_K^{-1}}^w=\inner{a_K^{-1}}^wH_2,
\]
with  $H_1, H_2$ belonging  to  $ \mathcal B (L^2)$,   the space of bounded operators
on $L^2$.
   \end{enumerate}
   \end{prop}
 Recall that in H\"{o}rmander's terminology, $a \in S( \tilde{a} , \Gamma)$ means that for  all multi-indices $\alpha$ and $\beta$, there exists a constant $C_{\alpha, \beta}$ such that
  $$
  | \D^\alpha_v \D^\beta_\eta a (v, \eta) | \leq C_{\alpha, \beta} \tilde{a}(v,\eta).
  $$
The temperance then implies a correct definition for the associated  operators. We postpone to section 3 and the appendix a review of these standard notions of pseudo-differential calculus.


The exponents of derivative terms and weight terms  in   Theorem \ref{thmain} and Theorem \ref{th1}  seem  to be  optimal, since the symbolic estimates provided by Proposition \ref{estaa} implies  that the operator $\mathcal P$  should behave locally like a generalized Kolmogorov type operator
\[
     \partial_t+v\cdot \partial_x+\abs{D_v}^{2s},
\]
for which the exponent $2s/(2s+1)$ for the regularity in the time
and space variables  is indeed sharp by using a simple scaling argument (see also \cite{LMP}). In the particular case $s=1$ we recover formally the Landau equation  and our exponents (both in regularity and weight) match perfectly with the exponents in \cite{HK2011}.

The main ideas of our proofs of the above theorems rely on some formal
computations of symbols in \cite{Ale99}, on  the  method by multiplier used in
\cite{HK2011, Li11}  and some microlocal  techniques developed by
Lerner while using Wick quantization \cite{MR2477145}. We refer to Section \ref{further} for
some considerations about the methodology we used, and which comes form these previous works. Let us note
that functional estimates from a series of work of Alexandre et
al. \cite{AMUXY3, AMUXY2, AMUXY1, MR2679369} and Gressman et
al. \cite{gres-strain} are also helpful for a clear understanding of the structure of the collision operator, but a nice feature of our method is that we will be able to completely avoid the use of these previous estimates.  Note that there are some other methods to study the regularity of the transport equation; for instance the average arguments used by Bouchut \cite{Bouchut02} and  a  version of the uncertainty principle used by   Alexandre et al. \cite{MR2462585} to  prove the regularity in the time and
space variables $t, x$.  However, these results do not provide any optimal hypoelliptic estimate for the spatially
inhomogeneous Boltzmann equation without angular cutoff.

We give now some bibliographical references  about the
hypoelliptic properties of the non cutoff Boltzmann equation and related kinetic models.  Note that the angular cross-section $b$ is not integrable on the sphere due to the singularity  $\theta^{-2-2s}$, which leads to the formal statement that the nonlinear collision operator should behave like a fractional Laplacian; that is,
\begin{eqnarray*}
    Q(g,f)\approx - C_g(-\triangle_v)^sf+\textrm{lower order terms},
\end{eqnarray*}
with $C_g>0$ a constant depending only on the physical properties of $g$.   Initiated by  Desvillettes \cite{Desv97,MR1324404},  there have been extensive works around this result and regarding the  smoothness of solutions for the homogeneous Boltzmann equation without angular cutoff,  c.f.  \cite{ADVW00, MR2149928, MR2557895, DesvVillani00-2, DesvWennberg04, MR2425608, MR2476686, MR2523694}.   For the inhomogeneous case the study becomes more complicated. We remark that there have been some related works  concerned with the linear model of spatially  inhomogeneous Boltzmann equation, which takes the following form
\[
  \partial_t  +v\cdot\partial_x   + e(t,x,v)(-\triangle_v)^s,\quad \inf_{t,x,v} e(t,x,v)>0.
\]
This model equation was firstly studied by Morimoto and Xu
\cite{MorimotoXu07b}, where a global but non optimal  hypoelliptic
estimate was established.  This study was then improved by Chen et
al. \cite{MR2763329}, and  also by Lerner et al. in \cite{LMP} for an
optimal local result. We also mention \cite{A1} where a  simple proof
of the subelliptic estimate for the above model operator is given.  For general inhomogeneous Boltzmann equation we refer to \cite{AMUXY3, AMUXY2, AMUXY1, MR2679369} for recent progress on its  qualitative properties. Finally, let us also mention a recent global result by Lerner et al. \cite{LMPX11}  in the radially symmetric case and the Maxwellian case (which corresponds to $\gamma =0$ in our notations), and   closely related works
\cite{gres-strain, gres-strain1,LMPX133} where the sharp  estimates for the Boltzmann collision operator were explored.

\red

\subsection{Further comments and methodolodgy} \label{further}

  In this subsection, we give some additional comments on this work and explain the general strategy of the proofs.

  \bigskip \it On the linear approach. \rm First mention that we focus in this article on a linearized Boltzmann operator. We note that a deep knowledge of the linear behavior is of great interest in the study of the non-linear case, at least in a perturbative context (see for example \cite{AMUXY1,AMUXY2,AMUXY3,gres-strain} and the references therein for this without cutoff case). These previous works are mainly concerned with the global existence of solutions close to equilibrium for the the fully non-linear Boltzmann equation, and important parts of the proofs are connected with functional properties of the linearized part of Boltzmann collisional operator. Our main goal here is to understand  the functional properties of the linearized part of the fully inhomogeneous equation.

  \bigskip \it On the kernel of the collision operator. \rm We emphasize the fact that we are absolutely not interested in the (finite-dimensional)  kernel $\nnn$ of the linearized Boltzmann collision operator. This is an a priori independent question to establish so called hypocoercive estimates on the orthogonal of $\nnn$ and related exponential return to the equilibrium of the solutions of the Boltzmann equation. We only deal here with regularity or hypoelliptic issues.

  \bigskip \it On the interest of regularization estimates. \rm
  In this article we essentially focus on global hypoelliptic estimates concerning the linearized Boltzmann operator
  $\ppp$ defined in \eqref{mpintro}. The main result in Theorem \ref{thmain} just concerns the independent of time problem and implies the following type of result. If one consider an equality $\ppp f = g$ with given $f, \ g \in L^2$,  then in fact $f$ has a better regularity and space/velocity decay given by the inequality in Theorem \ref{thmain} :   it has some weighted $H^{2s}$ regularity in velocity and $H^{2s/(2s+1)}$ regularity in space.
   Note that this kind of conclusion is \it not \rm available if one only use triple norm estimates (see the version given in  remark \ref{remarktriple} here) for which space regularity is not given.

   Mention that estimates like in Theorem \ref{thmain} and the careful study of the pseudo-differential and hypoelliptic structure of diffusive inhomogeneous kinetic equations have concrete applications; for example many ideas and tools developed here lead in \cite{HTT17a} and \cite{HTT17b} to the existence and uniqueness of solutions of the full non-linear inhomogeneous Boltzmann equation without cutoff with close to equilibrium initial data in large spaces (in the spirit of the theory developed recently in \cite{GMM13}).

 \bigskip \it A multiplier method. \rm
%
In this work we  make use of multiplier method  to explore the intrinsic hypoelliptic structure of operator $\ppp= v.\D_x -\lll $ defined in \eqref{mpintro}. By multiplier method  we mean finding a bounded selfadjoint operator $\mathcal M,$  such that on one side the commutator between the transport part and $\mathcal M$
\begin{eqnarray*}
	{1\over2}\inner{[\mathcal M, v\cdot \partial_x ]u,  ~u}_{L^2}=  {\rm Re} \inner{v\cdot \partial_x u, ~\mathcal M u}_{L^2}
\end{eqnarray*}
gives some ``elliptic'' properties in spatial variables, and on the other  side we can control the upper bound for the term
\begin{eqnarray*}
\abs{\inner{ \mathcal Lu,  ~\mathcal M u}_{L^2}	}.
\end{eqnarray*}
  For the treatment of the latter we need to the representation of $\mathcal L$ in term of pseudo-differential operators (see Proposition \ref{estaa}) which will be useful to estimate the commutators between $\mathcal L $ and $\mathcal M.$  The choice of the multiplier here is inspired by the Poisson bracket analysis for the transport part and the collision part already done for other diffusive models (see e.g. Fokker-Planck or Landau in \cite{HK2011} or \cite{Li11}).

   \bigskip \it The multiplier method explained on a toy model. \rm
  To clarify the choice of the multiplier $\mmm$ above we consider the case when  $\ppp$ is replaced by a Kolmogorov  type operator $\ppp_{kol}$
\begin{eqnarray*}
	\ppp_{kol} = v\cdot \partial_x -\Delta_v.
\end{eqnarray*}
(This corresponds to $\gamma=0$ and $s=1$ in a simplified case).
Then a direct computation gives
\begin{eqnarray*}
	\Big[v\cdot \partial_x, -\Delta_v \Big]= 2\partial_{x}\cdot\partial_v
\qquad
	\Big[v\cdot \partial_x,\Big[v\cdot \partial_x, -\Delta_v \Big]\,
\Big]=-2\Delta_x,
\end{eqnarray*}
 and we observe that this second-commutator analysis exhibit some Laplacian in $x$.  This suggests that the multiplier should be similar to the first-order commutator $2\partial_{x}\cdot\partial_v$.
   Since it is not a bounded operator on  $L^2$  we have to modify the  multiplier to guarantee its boundedness. It is then easier to see all the computation on the Fourier side : let    $\xi$ be the dual of $x$ and let $\eta$ be the dual of $v$. then operator $2\partial_{x}\cdot\partial_v$ is represented by a multiplication  by $-2 \xi\cdot \eta$ and we note that the Laplacian in velocity is a multiplication by $-|\eta|^2$ on the Fourier side. Then a good multiplier $\mmm$ is given by the quantization of the following \it bounded \rm function
\begin{eqnarray*}
m(\xi,\eta) = \frac{  \xi\cdot\eta }{\seq{\xi}^{4/3}} \chi\left( \frac{ \seq{\eta} }{ \seq{\xi}^{1/3}} \right),
\end{eqnarray*}
where $\chi\in C_0^\infty(\mathbb R;~[0,1])$ such that $\chi=1$ in
$[-1,1]$ and supp~$\chi \subset[-2,2].$ This function is clearly bounded thanks to the localization induced by $\chi$ on small $\eta$ frequencies, and it has to be considered as a (truncated and weighted) modification of the fundamental stone $\xi \cdot \eta$. Computation of all involved commutators on the fourier side give then
$$
-\Delta_v + [ v \cdot \D_x, \mmm ] \simeq -\Delta_v +(- \Delta_x)^{1/3} + errors
$$
leading after some work to subelliptic estimates of the form
 \begin{equation*}
 \begin{split}
  \ \big\Vert \seq{D_v}^{2} f \big\Vert
   + \big\Vert  \seq{D_x}^{2/3} f \big\Vert
   \leq C \sep{ \norm{{\mathcal P_{kol}}f} + \norm{ f}},
\end{split}
\end{equation*}
for (compactly) supported smooth functions.
We refer to
  \cite{HK2011} for  more developed arguments about this method, and complete computations in some simple cases.  Theorem \ref{thmain} is of the same form but global, with weights involving velocity and with regularity $2s$ or $2s/(2s+1)$ instead of $2$ or $2/3$ because of the structure of the Boltzmann collision operator without cut-off. The proof is also much more complicated than for the previous toy model.

\bigskip \it On the use of the Wick quantization. \rm
  In the example just before,  $\mmm$ was just a standard Fourier multiplier. In the case of the Boltzmann collision operator, the corresponding operator  has a more tricky structure and has to be selected into the general family of pseudo-differential operators. Its construction follows anyway exactly the same ideas as before (see Subsection \ref{subsec44} for its expression). Now in all these strategies the positivity of the symbols, multipliers and their commutators is an important point, and it appears that one cannot apply standard positivity result of operators having non-negative symbols (as the famous G{\aa}rding inequality) since they are in bad classes in the sense of H\"ormander t(see e.g. \cite{Hormander85} chapter 18 or \cite{MR2599384}).

  Anyway by choosing the Wick quantization of symbols, we can bypass this difficulty~: recall indeed that for any symbol $
  q\geq 0$  we directly have
$q^{\rm Wick}\geq 0 $ in the sense of operators.
We will use the Wick quantization   here  instead of  the classical or the  Weyl ones, and this  will simplify our  arguments substantially :  the computations and inequalities can be directly stated on symbols.

\black


\subsection{Organization of the article}
The paper is organized as follows. In Section 2, we provide precise estimates on the nice terms appearing in the splitting of the collision operator $\lll = \lll_1 + \lll_2$, involving compact parts and relatively bounded terms w.r.t. the operator of multiplication by $\seq{v}^{\gamma+ 2s}$. In Section $3$ we deal with the main terms, which appear to be of pseudo-differential type, and give precise symbolic estimates in the sense of the Weyl-H\"{o}rmander calculus. Section 4 is devoted to the proof of the main theorems. An appendix is devoted to a short review of some tools used in this work (Wick quantization, cancellation Lemma and Carleman representation).



\section{First estimates on the linearized collision operator}

In this section we study the linearized collision  part $\lll$ defined in \eqref{defdeL}. We cut it in many pieces and study each of them except the two principal ones, which study is postponed in section 3 (they are indeed of  pseudo-differential type). We look here at the  properties of the non pseudo-differential parts, and write many estimates in weighted $L^2$ spaces.

The splitting of the linearized Boltzmann operator $\lll$ is as follows. We write of $f \in \sss$,
\begin{eqnarray}\label{L1L2def}
\begin{aligned}
\lll f  & = \mu^{-1/2} Q( \mu, \mu^{1/2} f) + \mu^{-1/2} Q( \mu^{1/2} f, \mu)  \\
& =  \mu^{-1/2} \iint dv_* d\sigma B \sep{ \mu'_* (\mu')^{1/2}f' -  \mu_* \mu^{1/2}f +   \mu'(\mu'_*)^{1/2} f'_* -  \mu (\mu_*)^{1/2} f_* } \\
& = \iint dv_* d\sigma B (\mu_*)^{1/2}\sep{ (\mu'_*)^{1/2} f' -  (\mu_*)^{1/2} f +   (\mu')^{1/2} f'_* -  (\mu)^{1/2} f_* } \\
& = \iint dv_* d\sigma B (\mu_*)^{1/2} \sep{ (\mu'_*)^{1/2} f' -  (\mu_*)^{1/2} f}   \\
&\qquad+    \iint dv_* d\sigma B (\mu_*)^{1/2} \sep{ (\mu')^{1/2} f'_* -  (\mu)^{1/2} f_* } \\
& = \lll_1 f + \lll_2 f.
\end{aligned}
\end{eqnarray}

We shall study more precisely each part of $\lll$. Let us immediately point out that
they have completely different behaviors. The non local term $\lll_2$
behaves essentially like a convolution term, with nice estimates, and
is relatively compact w.r.t. the main part of $\lll_1$ which will
appear to be of pseudo-differential type.


\subsection{Study of $\lll_2$}

Starting from the expression of $\lll_2$ given by
\begin{equation*}
\begin{split}
\lll_2 f & =
\iint dv_* d\sigma B (\mu_*)^{1/2} \sep{ (\mu')^{1/2} f'_* -  (\mu)^{1/2} f_* },
\end{split}
\end{equation*}
we split it into four terms which make sense even for strong singularities of $B$, i.e. in particular for $s\geq 1/2$. This point will be clear from the proof of Lemma \ref{l2l2} below.
\begin{equation*}
\begin{split}
\lll_2 f & =    \iint dv_* d\sigma B (\mu_*)^{1/2} \sep{ (\mu')^{1/2} f'_* -  (\mu)^{1/2} f_* }  \\
& = \iint dv_* d\sigma B \sep{ (\mu^{1/2} f)'_* (\mu')^{1/2} - (\mu^{1/2} f)_* \mu^{1/2}} +
 \iint dv_* d\sigma B   (\mu ')^{1/2} \sep{ (\mu_*)^{1/2} -  (\mu'_*)^{1/2} } f'_* \\
 & = \iint dv_* d\sigma B  (\mu^{1/2} f)'_* \sep{ (\mu')^{1/2} - \mu^{1/2}} \\
 & \ \ \ \ \ \ \ + \mu^{1/2} \iint dv_* d\sigma B  \sep{ (\mu^{1/2} f)'_* -(\mu^{1/2} f) _*  } \\
&  \ \ \ \ \ \ \ \ \ \ \ \  +
\mu^{1/2} \iint dv_* d\sigma B    \sep{ (\mu_*)^{1/2} -  (\mu'_*)^{1/2} } f'_* \\
& \ \ \ \ \ \ \ \ \ \ \ \ \ \ \ \ \ \ + \iint dv_* d\sigma B  \sep{  (\mu ')^{1/2} - (\mu )^{1/2} }   \sep{ (\mu_*)^{1/2} -  (\mu'_*)^{1/2} } f'_* \\
 & = \lll_{2,r} f + \lll_{2,ca} f + \lll_{2,c} f +  \lll_{2,d} f.
\end{split}
\end{equation*}
$\lll_{2,ca}$ involves essentially a convolution term and can be treated using the cancellation lemma (see \cite{AV02} and the appendix herein), and the three other ones can be estimated by hands.
Let us note that the analysis of $\lll_2$ was already given by \cite{AMUXY1}, Lemma 2.15, but we provide a somewhat direct and shorter proof.

\begin{lem} \label{l2l2}  For all $\alpha$, $\beta \in \R$ there exists a constant
  $C_{\alpha,\beta}$ such that for all $f\in\mathcal S(\mathbb R_v^3)$ we have
 \begin{equation*} 
   \norm{ \seq{v}^\alpha \lll_2 \seq{v}^\beta f} \leq C_{\alpha,\beta} \norm{f}.
  \end{equation*}
\end{lem}

\preuve  We start with $ \lll_{2,ca} f $:
$$ \lll_{2,ca} f  = \mu^{1/2} \iint dv_* d\sigma B  \sep{ (\mu^{1/2} f)'_* -(\mu^{1/2} f) _*  } .$$
Applying the Cancellation Lemma (see  \cite{AV02} or the appendix), we get, for some constant $c$ depending only on $b$:
$$ \lll_{2,ca} f  = c \mu^{1/2} \int dv_* | v-v_*|^\gamma (\mu^{1/2} f) _*  .$$
This is an integral operator with the kernel $K(v,v_*) = c \mu^{1/2}
(\mu_*)^{1/2} |v-v_*|^\gamma$ for which we can apply Schur's Lemma to get
$$\|  \lll_{2,ca} f \| \lesssim \|  f\| .$$
Note that the assumption $\gamma >-3$ is needed at this point.

More generally, replacing $\lll_{2,ca} f$ by $\seq{v}^\alpha \lll_{2,ca} \seq{v}^\beta f$ leads to a kernel
$$K_{\alpha,\beta} (v,v_*) = c \mu^{1/2}\seq{v}^\alpha (\mu_*)^{1/2} \seq{v_*}^\beta |v-v_*|^\gamma$$
for which we can use the same argument to get
$$
 \norm{ \seq{v}^\alpha \lll_{2,ca} \seq{v}^\beta f} \leq C_{\alpha,\beta} \norm{f}.
 $$

Next, dealing with  $\lll_{2,c} f$
$$\lll_{2,c} f =\mu^{1/2} \iint dv_* d\sigma B    \sep{ (\mu_*)^{1/2} -  (\mu'_*)^{1/2} } f'_*  ,$$
we split this term into a singular and a non-singular parts. First consider the non singular part defined as
$$
\lll_{2,c,nonsing} f \defegal \mu^{1/2} \iint dv_* d\sigma B   \un_{|v'-v| \geq 1} \sep{ (\mu_*)^{1/2} -  (\mu'_*)^{1/2} } f'_*.
$$
As noticed in \cite{AMUXY1}, one has $\mu'_* \mu' = \mu_*\mu \leq (\mu'_* \mu)^{1/5}$ due to the kinetic and momentum relations in (\ref{kmr}). Therefore
$$
Af \defegal \abs{ \lll_{2,c,nonsing} f}  \lesssim  \mu^{1/10} \iint dv_* d\sigma |B|   \un_{|v'-v| \geq 1} \abs{ ( \mu^{1/10} f)'_* }
$$
which writes in Carleman representation (see the appendix)
$$
Af   \lesssim \mu^{1/10} \int_{\R^3_h} dh \int_{E_{0,h}} d\alpha  \un_{|h|\geq 1} \un_{|\alpha| \geq |h|}  \frac{|\alpha+h|^{1+\gamma + 2s}}{|h|^{3+2s}}    |(\mu^{1/10} f)(\alpha +v)|,
$$
where  $E_{0,h}$ denotes the hyperplane orthogonal to $h$ and containing $0$.
By duality, we get, for all $g \in \sss$,
\begin{equation*}
\begin{split}
|(Af, g)| & \lesssim \int_{\R^3_v} dv \int_{\R^3_h} dh \int_{E_{0,h}} d\alpha  \un_{|h|\geq 1} \un_{|\alpha| \geq |h|}  \frac{|\alpha+h|^{1+\gamma + 2s}}{|h|^{3+2s}}    |(\mu^{1/10} f)(\alpha +v)| \cdot|\mu^{1/10} g( v)|\\
& \lesssim \int_{\R^3_v} dv \int_{\R^3_h} dh \int_{E_{0,h}} d\alpha  \un_{|h|\geq 1} \un_{|\alpha| \geq |h|} \frac{|\alpha|^{1+\gamma + 2s}}{|h|^{3+2s}}    |(\mu^{1/10} f)(\alpha +v)| |(\mu^{1/10} g)( v)|
\end{split}
\end{equation*}
which upon using (\ref{interversion}) yields
\begin{equation*}
\begin{split}
|(Af, g)|  & \lesssim \int_{\R^3_v} dv \int_{\R^3_\alpha} d\alpha \int_{E_{0,\alpha}} dh \un_{|h|\geq 1} \un_{|\alpha| \geq |h|} \frac{|\alpha|^{\gamma + 2s}}{|h|^{2+2s}}    |\mu^{1/10} f(\alpha +v)| |\mu^{1/10} g( v)| \\
& \lesssim \int_{\R^3_v} dv \int_{\R^3_\alpha} d\alpha   |\alpha|^{(\gamma + 2s)^+}   |\mu^{1/10} f(\alpha +v)| |\mu^{1/10} g( v)|.
\end{split}
\end{equation*}

Therefore
$$
|(Af, g)| \lesssim \norm{ \mu^{1/20} f} \norm{ \mu^{1/20} g}
$$
from which follows that
\begin{equation} \label{nonsing}
 \norm{ \seq{v}^\alpha \lll_{2,c,nonsing} \seq{v}^\beta f} \leq C_{\alpha,\beta} \norm{f}
\end{equation}
for all real $\alpha$ and $\beta$.

For the singular part $\lll_{2,c,sing}$,  again using Carleman's representation \eqref{carl}
gives
\begin{multline*}
\lll_{2,c,sing} f
 = \mu^{1/2} \int_{\R^3_h}dh\int_{E_{0,h}} d\alpha \tb(\alpha,h) \un_{|\alpha | \geq | h|} \un_{| h| \leq 1}  \\ \sep{ \mu^{1/2}
 (\alpha +v -h) -\mu^{1/2} (\alpha +v) } { | \alpha +h |^{1+\gamma +2s} \over{| h|^{3+2s}}} f(\alpha +v).
 \end{multline*}

Changing $h\rightarrow -h$ and adding the resulting two formulas (so we see that formally we cancel higher singularities, using also that $\tb(\alpha,h) = \tb(\pm \alpha,\pm h)$) yields
\begin{multline*}
\lll_{2,c,sing} f  = {1\over 2} \mu^{1/2} \int_h dh\int_{E_{0,h}} d\alpha \tb  \un_{|\alpha | \geq | h|} \un_{|h | \leq 1}   \times \\
  \sep{ \mu^{1/2} (\alpha +v -h) + \mu^{1/2} (\alpha +v +h)  -2\mu^{1/2} (\alpha +v) } {| \alpha +h |^{1+\gamma +2s}\over{| h|^{3+2s}}} f(\alpha +v).
  \end{multline*}
Factorizing by $\mu^{1/2} (\alpha +v)$ we get
\begin{eqnarray*}
&&\lll_{2,c,sing} f \\
 &= &{1\over 2} \mu^{1/2} \int_h dh  \int_{E_{0,h}} d\alpha \tb \un_{|\alpha | \geq | h|} \un_{|h | \leq 1}
 \sep{ e^{-\inner{\abs h^2-2(\alpha+v)\cdot h }/4  } + e^{-\inner{\abs h^2+2(\alpha+v)\cdot h}/4  }  - 2  } \\
 &&\qquad \times {| \alpha +h |^{1+\gamma +2s}\over{| h|^{3+2s}}}
\mu^{1/2} (\alpha +v) f(\alpha+v).
  \end{eqnarray*}
  The term in parentheses is bounded by $|h|^2 \mu^{-1/4}(\alpha+v)$ thanks to the condition on the support for $h$, and since $|h| \leq |\alpha|$, one has
 \begin{multline*}
|\lll_{2,c,sing} f|  \lesssim  \mu^{1/2} \int_h dh\int_{E_{0,h}} d\alpha \un_{|\alpha | \geq | h|} \un_{|h | \leq 1}
    {| \alpha |^{1+\gamma +2s}\over{| h|^{1+2s}}}
\mu^{1/4} (\alpha +v) | f(\alpha+v) |.
  \end{multline*}
Using again (\ref{interversion}) and the duality argument as in the non-singular case (now the singularity in $h$ is integrable), we easily get
\begin{equation} \label{sing}
 \norm{ \seq{v}^\alpha \lll_{2,c,sing} \seq{v}^\beta f} \leq C_{\alpha,\beta} \norm{f}
\end{equation}
for all real $\alpha$ and $\beta$.


As for $\lll_{2,r} f $, recalling that

$$\lll_{2,r} f =  \iint dv_* d\sigma B  (\mu^{1/2} f)'_* \sep{ (\mu')^{1/2} - \mu^{1/2}} $$
we see immediately that, using the classical pre-post velocities change of variables that
$$( \lll_{2,r} f , g) = (f , \lll_{2,c} g )$$
and thus we are done for this term.

It remains to study $\lll_{2,d} f$ which is exactly

$$\lll_{2,d} f  = \iint dv_* d\sigma B  \sep{  (\mu ')^{1/2} - (\mu )^{1/2} }   \sep{ (\mu_*)^{1/2} -  (\mu'_*)^{1/2} } f'_*  .$$
Using the equality $a^2 -b^2 = (a-b)(a+b)$ for the Gaussian functions in the above factors, we see again that we can put some power of a Gaussian together with $f$, by using the argument of \cite{AMUXY1}: that means that for some $c>0, d>0$, one has

$$|\lll_{2,d} f | \lesssim \mu^{d}\iint dv_* d\sigma B  \abs{  (\mu ')^{1/4} - (\mu )^{1/4} }   \abs{ (\mu_*)^{1/4} -  (\mu'_*)^{1/4} } (\mu^c)'_* | f'_* | $$
and then the remaining analysis is exactly similar to the computations done for $\lll_{2,c,sing} f$.

\fin

\subsection{Splitting of  $\lll_1$}\label{sub22n}
The operator $\lll_1$ will also be cut into several pieces, which will require two different types of arguments. For some of the nice parts, tools similar to the ones in the previous section
will be sufficient. The remaining pseudo-differential parts will be treated in the next Section.

Recall first that
\begin{equation*}
\begin{split}
\lll_1 f & =
\iint dv_* d\sigma B (\mu_*)^{1/2} \sep{ (\mu'_*)^{1/2} f' -  (\mu_*)^{1/2} f}.
\end{split}
\end{equation*}

Let $0<\delta\leq 1$ be a fixed parameter in the following argument.    We first split the above integral according to whether or
not $|v'-v|\gtrsim \delta$. To this end, let $\phi$ be a positive radial function supported on the unit ball and say $1$ in the $1/4$ ball. Consider $\phi_\delta (v)  = \phi (|v|^2 /\delta^2 )$, which is therefore $0$ for $|v| \geq \delta$ and $1$ for $|v| \leq \delta/2$. By abuse of notations we shall also denote $\phi_\delta(r) = \phi_\delta(v)$ when $|v|= r$.  Set $\tilde\phi_\delta (v ) = 1- \phi_\delta (v)$, which is therefore $0$ for small values and $1$ for large values.

Then $\lll_1 f$ can decomposed as the sum of the following two terms
$$\bar \lll_{1,\delta} f= \iint dv_* d\sigma B \tilde\phi_\delta (v'-v)  (\mu_*)^{1/2} \sep{ (\mu'_*)^{1/2} f' -  (\mu_*)^{1/2} f}$$
and
$$\lll_{1,\delta} f= \iint dv_* d\sigma B \phi_\delta (v'-v)(\mu_*)^{1/2} \sep{ (\mu'_*)^{1/2} f' -  (\mu_*)^{1/2} f}.$$

Note that $\bar\lll_{1,\delta}$ is a cutoff type Boltzmann operator. We split it into two terms since there is no singularity any more
\begin{eqnarray}\label{L1detaa}
\begin{aligned}
\bar \lll_{1,\delta} f
&  = \iint dv_* d\sigma B \tilde\phi_\delta (v'-v)  (\mu_*)^{1/2}  (\mu'_*)^{1/2} f' \\
&   \ \ \ \ \ - \sep{ \iint dv_* d\sigma B \tilde\phi_\delta (v'-v)  (\mu_*)^{1/2}
 (\mu_*)^{1/2}} f \\
& = \bar \lll_{1,\delta, a} f + \bar \lll_{1,\delta, b} f.
\end{aligned}
\end{eqnarray}

As for $\lll_{1,\delta}$, again we split it into four terms:
\begin{equation} \label{splitL1+}
\begin{split}
\lll_{1,\delta} f & =
\iint dv_* d\sigma B \phi_\delta (v'-v)(\mu_*)^{1/2} \sep{ (\mu'_*)^{1/2} f' -  (\mu_*)^{1/2} f} \\
 & = \iint dv_* d\sigma B \phi_\delta (v'-v) (\mu'_*)^{1/2}  (f' - f) \sep{ (\mu_*)^{1/2} - (\mu'_*)^{1/2}} \\
&  \ \ \ \ \ \ \ + \sep{ \iint dv_* d\sigma B \phi_\delta (v'-v) (\mu'_*)^{1/2} \sep{ (\mu_*)^{1/2} - (\mu'_*)^{1/2}} } f \\
 & \ \ \ \ \ \ \ \ \ \ \ \ \ \ +  \iint dv_* d\sigma B \phi_\delta (v'-v)\mu'_* \sep{ f' -f  } \\
&  \ \ \ \ \ \ \ \ \ \ \ \ \ \ \ \ \ \ \  +
 \sep{ \iint dv_* d\sigma B \phi_\delta (v'-v)   \sep{ \mu'_* -  \mu_* }}  f \\
 & = \lll_{1,1,\delta} f + \lll_{1,4,\delta} f + \lll_{1,2,\delta} f + \lll_{1,3,\delta}f.
\end{split}
\end{equation}

Let us immediately notice that this splitting takes into account all values of $s$. However, for small singularities $0< s <1/2$, a simpler decomposition is available and avoids some of the issues dealt with below. We note that $\lll_{1,4,\delta} f$ and $\lll_{1,3,\delta} f$ are of multiplicative type,  and together with  $\bar{\lll}_{1,\delta,a} f$, they
will be studied in the next subsection. They will appear later as relatively bounded terms with respect to   $\lll_{1,1,\delta} + \lll_{1,2,\delta} f$. These last two parts will appear to be of pseudo-differential type, and we shall estimate them very precisely in  section 3.

\remark In the coming computations, we shall follow the dependence on parameter $\delta$. We point out that it could  be  fixed at value $\delta = 1$. Anyway, as we shall see in the coming sections, the explicit dependence on $\delta$ of the various estimates enlightens the fact that we are the non cutoff case. As already mentioned, the cutoff case corresponds to the case when $\lll_{1,\delta} = 0$. It can also be seen as the limiting case $\delta \rightarrow 0$ when looking e.g. at $\lll_{1,2,\delta}$, for which we give in Proposition \ref{ap2} a lower bound which would be not relevant anymore for $\delta = 0$.

\subsection{Relatively bounded terms in $\lll_1$}

\subsubsection{Study of  $\lll_{1, 3,\delta}$}

Using some arguments from the proof of the cancellation lemma, see for example \cite{AV02}, we get the following

 \begin{lem} \label{L1easy} For all $f\in\mathcal S(\mathbb R_v^3)$, we have, for all $s<1$
$$
\|  \lll_{1,3,\delta} f  \|^2\lesssim  \delta^{2-2s} \| <v>^{\gamma+2s-2} f \|
$$
and $\lll_{1,3,\delta}$ commutes with the multiplication by $\seq{v}^\alpha$ for all $\alpha \in \R$.
\end{lem}

 \preuve
The last assertion is trivial since $\lll_{1, 3,\delta}$ is a multiplication operator. In order to prove the above inequality,
recall first that
$$\lll_{1,3,\delta} f (v) = \sep{ \iint dv_* d\sigma B \phi_\delta (v'-v)   \sep{ \mu'_* -  \mu_* }}  f .$$
Going back to the proof of the cancellation Lemma, it follows that
$$ \sep{ \iint dv_* d\sigma B \phi_\delta (v'-v)   \sep{ \mu'_* -  \mu_* }}  = S\ast_{v_*} \mu (v) $$
where, writing by abuse of notation $\phi_\delta(|z|) \defegal \phi_\delta(z)$ for all $z \in \R^3$, $S$ has the following expression
\begin{equation*}
\begin{split}
S (z)  = & | z|^\gamma \int^{\pi /2}_0  \sin \theta b(\cos \theta ) \sep{ \phi_\delta ( {{|z|}\over{\cos {\theta \over 2}}}\sin {\theta \over 2}) \cos^{-3-\gamma } {\theta \over 2} -\phi_\delta ( | z| \sin {\theta \over 2})}d\theta \\
  = & | z|^\gamma \int^{\pi /2}_0  \sin \theta b(\cos \theta ) \phi_\delta ( {{|z|}\over{\cos {\theta \over 2}}}\sin {\theta \over 2})\sep{ \cos^{-3-\gamma } {\theta \over 2} -1}d\theta\\
 &+ | z|^\gamma \int^{\pi /2}_0  \sin \theta b(\cos \theta ) \sep{ \phi_\delta ( {{|z|}\over{\cos {\theta \over 2}}}\sin {\theta \over 2}) -\phi_\delta ( | z| \sin {\theta \over 2})}d\theta \\
 =& S_1(z) + S_2(z).
\end{split}
\end{equation*}
For the first part $S_1(z)$, note that the integrand is now integrable in the  $\theta$ variable, and we have
 \begin{equation} \label{S1}
 |S_1(z)| \lesssim \abs{z}^\gamma.
 \end{equation}
 The second part $S_2(z)$ is zero if $|z| \leq \delta/2$, and we can suppose therefore that  $|z| \geq \delta/2$. Note also that for $z$ bounded, say for $|z| \leq C$ where $C$ is sufficiently large to be fixed later, $S_2(z)$ is also bounded. Since
$$
{{|z|}\over{\cos {\theta \over 2}}}\sin {\theta \over 2} \geq | z| \sin {\theta \over 2},
$$
we get that if $| z| \sin {\theta \over 2} \geq \delta $, the integrand is $0$, and similarly for small values of $\theta$. In conclusion when $|z| \geq C$,  the second integral can be estimated as follows :
$$
S_2(z)  = | z|^\gamma \int^{c\delta |z|^{-1}}_{c'\delta |z|^{-1}}   \sin \theta b(\cos \theta ) \sep{ \phi_\delta ( {{|z|}\over{\cos {\theta \over 2}}}\sin {\theta \over 2}) -\phi_\delta ( | z| \sin {\theta \over 2})}d\theta,
$$
where $C$ is a posteriori chosen so that $C^{-1}  c \delta  \leq \pi/2$.
Using Taylor formulae, we get
\begin{equation*}
\begin{split}
|S_2 (z)|
&\lesssim  \delta^{-1} | z|^{\gamma +1}\int^{c\delta |z|^{-1}}_{c'\delta |z|^{-1}}    \theta^2 b(\cos \theta )[ \cos^{-1} \theta /2 -1]d\theta  \lesssim \delta^{-1} | z|^{\gamma +1}\int^{c\delta |z|^{-1}}_{c'\delta |z|^{-1}}    \theta^4 b(\cos \theta )d\theta  \\
& \lesssim \delta^{-1} | z|^{\gamma +1}\int^{c\delta |z|^{-1}}_{c'\delta |z|^{-1}}  \theta^{2-2s}d\theta \sim \delta^{-1} | z|^{\gamma +1} \delta^{3-2s} |z|^{-3+2s} \\
& \lesssim \delta^{2-2s} |z|^{\gamma +2s -2} .
\end{split}
\end{equation*}
This estimate together with (\ref{S1}) yield the proof of the Lemma. \fin

\subsubsection{Study of $\bar \lll_{1,\delta ,a}$}
 We deal now with the non singular part $\bar \lll_{1,\delta ,a}$ of $\lll_1$ for which we have the following result

\begin{lem} \label{lb1da} (i) For all $f\in\mathcal S(\mathbb R_v^3)$  and for all $\alpha$, $\beta \in \R$ such that $\alpha+ \beta + \gamma+2s \leq 0$, we have
 \begin{equation*}
   \norm{ \seq{v}^\alpha \bar \lll_{1,\delta ,a} \seq{v}^\beta f} \leq \delta^{-1-2s} C_{\alpha,\beta} \norm{f}.
  \end{equation*}
(ii) For all $f\in\mathcal S(\mathbb R_v^3)$  and for all $\tilde\alpha$, $\tilde\beta \in \R$ such that $\tilde\alpha+ \tilde\beta + \gamma+s \leq 0$, we have
 \begin{equation*}
   \norm{ \seq{v}^{\tilde\alpha} \com{\bar \lll_{1,\delta
         ,a},~\seq{v}^{\tilde \beta}} f} \leq \delta^{-2s} C_{\tilde
     \alpha,\tilde\beta} \norm{f},
  \end{equation*}
where $\com{\cdot,~\cdot}$ stands for the commutator.
\end{lem}

 \preuve
Recalling that
$$
\bar \lll_{1,\delta ,a} f= \iint dv_* d\sigma B \tilde\phi_\delta (v'-v)  (\mu_*)^{1/2}  (\mu'_*)^{1/2} f',
$$
it follows that
$$
\seq{v}^\alpha \bar \lll_{1,\delta ,a} \seq{v}^\beta f= \seq{v}^\alpha \iint dv_* d\sigma B \tilde\phi_\delta (v'-v)  (\mu_*)^{1/2}  (\mu'_*)^{1/2} (\seq{v}^\beta f)',
$$
and
\begin{eqnarray*}
\seq{v}^{\tilde \alpha} \com{\bar \lll_{1,\delta ,a},~ \seq{v}^{\tilde\beta}}
f&=&\seq{v}^{\tilde\alpha }\bar \lll_{1,\delta ,a}
\seq{v}^{\tilde\beta} f-\seq{v}^{\tilde\alpha} \seq{v}^{\tilde \beta} \bar \lll_{1,\delta ,a} f\\
&=& \seq{v}^{\tilde\alpha} \iint dv_* d\sigma B \tilde\phi_\delta (v'-v)  (\mu_*)^{1/2}  (\mu'_*)^{1/2} \inner{\seq{v'}^{\tilde\beta} -\seq{v}^{\tilde\beta} }f' .
\end{eqnarray*}
(i) We first estimate  $\seq{v}^\alpha \bar \lll_{1,\delta ,a} \seq{v}^\beta f.$ An application of Carleman's representation (see the appendix for instance) shows that
\begin{equation} \label{expba}
\begin{split}
 | \seq{v}^\alpha \bar \lll_{1,\delta ,a} \seq{v}^\beta f| & \lesssim \seq{v}^\alpha
 \int_h dh\int_{E_{0,h}} d\alpha \un_{|\alpha | \geq | h|} \un_{ | h| \geq \delta/2}    \mu^{1/2} (\alpha +v) \mu^{1/2} (\alpha +v -h)  \\
 & \ \ \ \ \ \ \ \ \ \ \ \ \ \ \ \   {| h+\alpha |^{1+\gamma +2s}\over{| h|^{3+2s}}} \seq{v-h}^\beta|f(v-h)| \\
 & \lesssim \seq{v}^\alpha \int_h dh\int_{E_{0,h}} d\alpha  \un_{ | h| \geq \delta/2}    \mu^{1/2} (\alpha +v) \mu^{1/2} (\alpha +v -h)  \\
 & \ \ \ \ \ \ \ \ \ \ \ \ \ \ \ \   {|\alpha |^{1+\gamma +2s}\over{| h|^{3+2s}}} \seq{v-h}^\beta |f(v-h)|,
\end{split}
\end{equation}
where we used the fact that $|\alpha| \geq |h|$ for the second inequality,  and recalling that $E_{0,h}$ denotes the vector
plane containing $0$ and orthogonal to $h$.  Letting $S(h)$ for the orthogonal projection onto $E_{0,h}$, we can write
$$e^{- |\alpha +v|^2 } = e^{- |\alpha + S(h) v|^2} e^{ |S(h) v|^2 - |v|^2}$$
and similarly
\begin{equation*}
\begin{split}
e^{- |\alpha +v -h|^2}
  = e^{-|\alpha +S(h) (v)|^2 } e^{| S(h) (v-h) |^2- |v-h|^2},
\end{split}
\end{equation*}
and therefore
\begin{equation*}
\begin{split}
\mu^{1/2} (\alpha +v) \mu^{1/2} (\alpha +v -h)
& =  (2\pi)^{-3/2} \sep{e^{-  |\alpha + S(h) v|^2} \sep{   e^{2( |S(h) v|^2 - |v|^2)+ | v|^2- |v-h|^2} }^{1/2}}^{1/2} \\
& =(2\pi)^{-3/2} \sep{ e^{-  |\alpha + S(h) v|^2}  \sep{ e^{2( |S(h) v|^2 - |v|^2)+2 v\cdot h -|h|^2} }^{1/2}}^{1/2}.
\end{split}
\end{equation*}
Going back to (\ref{expba}), we obtain
\begin{equation*}
\begin{split}
 | \seq{v}^\alpha \bar \lll_{1,\delta ,a} \seq{v}^\beta f|
 &  \lesssim  \seq{v}^\alpha \int_{\R^3_h} dh\int_{E_{0,h}} d\alpha \un_{|\alpha | \geq | h|} \un_{ | h| \geq \delta/2}  {|\alpha |^{1+\gamma +2s} \over{| h|^{3+2s}}} \seq{v-h}^\beta |f(v-h)| \\
 &  \ \ \ \ \ \ \ \ \ \ \ \ \ \ \ \ \ \ \ \sep{ e^{-  |\alpha + S(h) v|^2} \sep{ e^{2( |S(h) v|^2 - |v|^2)+2 v\cdot h -|h|^2} }^{1/2}}^{1/2}.
  \end{split}
\end{equation*}
Performing the integration with respect to $\alpha$, it follows that
\begin{equation*}
\begin{split}
 | \seq{v}^\alpha \bar \lll_{1,\delta ,a} \seq{v}^\beta f| & \lesssim \seq{v}^\alpha \int_{\R^3_h} dh  \un_{ | h| \geq \delta/2}  \seq{S(h)v }^{1+\gamma +2s} {1 \over{| h|^{3+2s}}} \seq{v-h}^\beta |f(v-h)| \\
 &  \ \ \ \ \ \ \ \ \ \ \ \ \ \ \ \ \ \ \ \sep{ e^{2( |S(h) v|^2 - |v|^2)+2 v\cdot h -|h|^2} }^{1/4} \\
 & \lesssim \int_{\R^3_z} dz \un_{ | v-z| \geq \delta/2} \seq{v}^\alpha \seq{S(v-z) v}^{1+\gamma +2s}   {1\over{| v-z|^{3+2s}}} \seq{z}^\beta | f|(z) \\
  &  \ \ \ \ \ \ \ \ \ \ \ \ \ \ \ \ \ \ \ \sep{ e^{2( |S(v-z) v|^2 - |v|^2)+ 2 v.(v-z) -|v-z|^2}}^{1/4} \\
 & \defegal \int_{\R^3_z} K_{\alpha,\beta}(v,z) |f|(z) dz
 \end{split}
\end{equation*}
with
\begin{multline*}
K_{\alpha,\beta} (v,z) = \un_{ | v-z| \geq \delta/2} \seq{v}^\alpha \seq{z}^\beta \seq{S(v-z) v}^{1+\gamma +2s}  {1\over{| v-z|^{3+2s}}} \\ \sep{ e^{2( |S(v-z) v|^2 - |v|^2)+ 2 v.(v-z) -|v-z|^2}}^{1/4}.
\end{multline*}
We want to apply Schur's Lemma. To this end, let's first integrate
w.r.t. to $z$, to get
\begin{equation*}
\begin{split}
\int_{\R^3_z} dz K_{\alpha,\beta} (v,z)
& = \int_{\R^3_z} dz \un_{ | v-z| \geq \delta/2} \seq{v}^\alpha \seq{z}^\beta \seq{S(v-z) v}^{1+\gamma +2s}  {1\over{| v-z|^{3+2s}}} \\
& \qquad \qquad \qquad \sep{ e^{2( |S(v-z) v|^2 - |v|^2)+ 2 v.(v-z) -|v-z|^2}}^{1/4} \\
& = \int_{\R^3_h} dh \un_{ | h| \geq \delta/2} \seq{v}^\alpha \seq{v-h}^\beta \seq{S(h) v}^{1+\gamma +2s}  {1\over{| h|^{3+2s}}} \\
& \qquad \qquad \qquad \sep{ e^{2( |S(h) v|^2 - |v|^2)+ 2 v\cdot h -|h|^2}}^{1/4},
\end{split}
\end{equation*}
so that
\begin{equation*}
\begin{split}
& \int_{\R^3_z} dz K_{\alpha,\beta} (v,z) \\
& = \int_{\R^3_h} dh \un_{ | h| \geq \delta/2} \seq{v}^\alpha  \sep{ 1+  |v|^2 - | v \cdot {h\over{| h| }}|^2 }^{(1+\gamma +2s)/2}   {1\over{| h|^{3+2s}}}\seq{v-h}^{\beta} \\
& \qquad \qquad \qquad \sep{ e^{- 2 | v \cdot {h\over{| h| }}|^2  + 2 v\cdot h -|h|^2}}^{1/4} \\
& = \int_{\R^3_h} dh \un_{ | h| \geq \delta/2} \seq{v}^\alpha  \sep{ 1+  |v|^2 - {|v|^2\over |h|^2} | {v\over{|v|}} \cdot h|^2 }^{(1+\gamma +2s)/2}   {1\over{| h|^{3+2s}}} \\
& \qquad \qquad \qquad \sep{ 1+ |v|^2 -2 |v| {v\over{|v|}} \cdot h +|h|^2 }^{\beta /2}\sep{ e^{- 2{|v|^2\over |h|^2} | {v\over{|v|}} \cdot h|^2 + 2 | v| {v\over{|v|}}\cdot h -|h|^2}}^{1/4}.
\end{split}
\end{equation*}
Shifting to polar coordinates, with an axis along direction $v/| v|$, we obtain
 \begin{multline*}
\int_{\R^3_z} dz K_{\alpha,\beta} (v,z)
  \lesssim \int^\pi_0 \int^\infty_\delta  dr d\phi \seq{v}^\alpha \sin \phi \sep{ 1  +|v|^2 - |v|^2 \cos^2 \phi }^{(1+\gamma +2s)/2}  {1\over r^{1+2s}} \\
 (1+ |v|^2 -2 |v| r \cos\phi +r^2 )^{\beta /2}
 \sep{e^{-2 |v|^2 \cos^2 \phi  + 2 |v| r\cos\phi - r^2}}^{1/4}.
\end{multline*}
Note here that if $|v| \leq 1$, then we directly get that $\int_{\R^3_z} dz K_{\alpha,\beta} (v,z) \lesssim  1$. Therefore we may as well assume that $|v| \geq 1$. Setting $t=\cos\phi$,  we get
\begin{equation*}
\begin{split}
&\int_{\R^3_z} dz K_{\alpha,\beta} (v,z) \\
&\lesssim \int^1_{-1}  \int^\infty_\delta  dr dt \seq{v}^\alpha ( 1  +|v|^2 - |v|^2 t^2  )^{(1+\gamma +2s)/2} e^{(-2 |v|^2 t^2  + 2 |v| rt - r^2)/4}   {1\over r^{1+2s}} \\
 &  \ \ \ \ \ \ \ \ \ \ \ \ \ \ \ \ \ \ \ (1+ |v|^2 -2 |v| r t +r^2 )^{\beta /2} \\
&\approx \seq{v}^\alpha |v|^{-1}\int^{|v|}_{-|v|}  \int^\infty_\delta  dr dt ( 1  +|v|^2 -  t^2  )^{(1+\gamma +2s)/2}  e^{-(r-t)^2/4}   {1\over r^{1+2s}} (1+ |v|^2 -2 r t +r^2 )^{\beta /2} \\
&\approx \seq{v}^\alpha |v|^{-1}\int^{|v|}_{-|v|}  \int^\infty_\delta  dr dt ( 1  +|v|^2 -  t^2  )^{(1+\gamma +2s)/2}  e^{-(r-t)^2/4}   {1\over r^{1+2s}}  (1+ |v|^2 -  t^2 +(r-t)^2 )^{\beta /2}.
\end{split}
\end{equation*}
In the inner term, note that $|v|^2 -  t^2  \geq 0$. We now use Peetre's inequality
\begin{equation}
	\label{peetre}
	\seq{u}^\beta \seq{u+w}^{-| \beta|} \lesssim \seq{w}^\beta \lesssim \seq{u}^\beta \seq{u+w}^{|\beta |} , \quad \beta\in\mathbb R,
\end{equation}
to get here
$$
(1+ |v|^2 -  t^2 +(r-t)^2 )^{\beta /2}  \lesssim (1+|v|^2 -  t^2  )^{\beta/2} \seq{r-t}^{|\beta |}.
$$
In addition, since $0<\delta<1$, then $r\geq \delta$ implies that $r\geq 2^{-1/2}\delta \seq{r}$.  Thus
\begin{equation} \label{inttoprecise}
\begin{split}
& \int_{\R^3_z} dz K_{\alpha,\beta} (v,z) \\
& \lesssim \delta^{-1-2s} \seq{v}^\alpha|v|^{-1}\int^{|v|}_{-|v|}  \int^\infty_{-\infty} dr dt ( 1  +|v|^2 -  t^2  )^{(1+\gamma +2s + \beta)/2} \sep{ \seq{r-t}^{|\beta|} e^{-(r-t)^2/4} }  {1\over \seq{r}^{1+2s}} \\
& \lesssim \delta^{-1-2s} \seq{v}^\alpha |v|^{-1}\int^{|v|}_{-|v|}  dt ( 1  +|v|^2 -  t^2  )^{(1+\gamma +2s +\beta)/2}   {1\over \seq{t}^{1+2s}} \\
& \lesssim \delta^{-1-2s} \seq{v}^{\alpha-1} \int_{0}^{|v|}  dt ( 1  +|v|^2 -  t^2  )^{(1+\gamma +2s +\beta)/2}   {1\over \seq{t}^{1+2s}}.
\end{split}
\end{equation}
Now for evaluating this quantity, we split the integral into two parts.
First note that
\begin{equation} \label{inttoprecisesmall}
\begin{split}
 \int_{0}^{|v|/2}  dt ( 1  +|v|^2 -  t^2  )^{(1+\gamma +2s +\beta)/2}   {1\over \seq{t}^{1+2s}}
& \qquad \lesssim \seq{v}^{1+\gamma +2s +\beta} \int_0^{|v|/2} dt {1\over \seq{t}^{1+2s}} \\
& \qquad \lesssim \seq{v}^{1+\gamma +2s +\beta}.
\end{split}
\end{equation}
For the remaining part, we write
\begin{equation} \label{inttopreciselarge}
\begin{split}
& \int_{|v|/2}^{|v|}  dt ( 1  +|v|^2 -  t^2  )^{(1+\gamma +2s +\beta)/2}   {1\over \seq{t}^{1+2s}} \\
& \qquad \lesssim \seq{v}^{-1-2s} \int_{|v|/2}^{|v|}  dt ( 1  +|v|^2 -  t^2  )^{(1+\gamma +2s +\beta)/2} \\
& \qquad \lesssim \seq{v}^{-1-2s} \int_{|v|/2}^{|v|}  dt ( 1  + (|v| -  t)(|v|+t)  )^{(1+\gamma +2s +\beta)/2} \\
& \qquad \lesssim \seq{v}^{-1-2s} \int_{|v|/2}^{|v|}  dt ( 1  + |v|(|v| -  t)  )^{(1+\gamma +2s +\beta)/2} \\
\intertext{ Posing $s= |v|(|v| -  t)$, $ds = - |v|dt$, we get }
& \int_{|v|/2}^{|v|}  dt ( 1  +|v|^2 -  t^2  )^{(1+\gamma +2s +\beta)/2}   {1\over \seq{t}^{1+2s}} \\
& \qquad \lesssim \seq{v}^{-1-2s} |v|^{-1} \int_0^{|v|^2/2}  ds ( 1  + s  )^{(1+\gamma +2s +\beta)/2} \\
& \qquad \lesssim \seq{v}^{-1-2s} |v|^{-1}  \seq{v}^{(1+\gamma +2s +\beta) + 2} \\
& \qquad \lesssim \seq{v}^{1+\gamma+\beta}.
\end{split}
\end{equation}
Putting estimates (\ref{inttoprecisesmall}) and (\ref{inttopreciselarge}) in (\ref{inttoprecise}) we get
\begin{equation*}
\begin{split}
 \int_{\R^3_z} dz K_{\alpha,\beta} (v,z)
& \lesssim  \delta^{-1-2s} \seq{v}^{\alpha-1} \seq{v}^{1+\gamma+\beta+2s} \\
& \lesssim \delta^{-1-2s} \seq{v}^{\alpha+ \gamma +\beta+2s} \\
& \lesssim \delta^{-1-2s} \ \ \ \textrm{ if } \alpha+\beta + \gamma +2s\leq 0.
\end{split}
\end{equation*}
In conclusion, we have obtained that if $\alpha+\beta + \gamma+2s \leq 0$, then
\begin{equation}\label{intz}
\int_{\R^3_z} dz K_{\alpha,\beta} (v,z) \lesssim  \delta^{-1-2s}.
\end{equation}


Now we look for the integration w.r.t. variable $v$ of $K_{\alpha,\beta}$. We have
\begin{equation*}
\begin{split}
 \int_{\R^3_v} dv K_{\alpha,\beta} (v,z)
& = \int_{\R^3_v} dv \un_{ | v-z| \geq \delta/2}  \seq{v}^\alpha \seq{z}^\beta \seq{S(v-z) v}^{1+\gamma +2s} \\
& \qquad \qquad \qquad \sep{ e^{ |S(v-z) v|^2 - |v|^2+| S(v-z) (z) |^2- |z|^2}}^{1/4}   {1\over{|v-z|^{3+2s}}},
\end{split}
\end{equation*}
since by direct computation
\begin{multline*}
2( |S(v-z) v|^2 - |v|^2)+ 2 v.(v-z) -|v-z|^2 \\
  =  |S(v-z) v|^2 - |v|^2+| S(v-z) (z) |^2- |z|^2.
\end{multline*}
Taking $h=v-z$, $dh=dv$, we get
\begin{equation*}
\begin{split}
 \int_{\R^3_v} dv K_{\alpha,\beta} (v,z)
& = \int_{\R^3_h} dh \un_{ | h| \geq \delta/2}  \seq{z+h}^\alpha \seq{z}^\beta \seq{S(h) (z+h)}^{1+\gamma +2s} \\
&  \qquad \qquad \qquad \sep{ e^{ |S(h) (z+h)|^2 - |z+h|^2+ | S(h) (z) |^2- |z|^2}}^{1/4}   {1\over{|h|^{3+2s}}} \\
& = \int_{\R^3_h} dh \un_{ | h| \geq \delta/2}  \seq{z+h}^\alpha \seq{z}^\beta \seq{S(h) z}^{1+\gamma +2s} \\
& \qquad \qquad \qquad \sep{ e^{ |S(h) z|^2 - |z+h|^2+ | S(h) z |^2- |z|^2}}^{1/4}  {1\over{|h|^{3+2s}}},
\end{split}
\end{equation*}
so that expanding again the brackets, we get
\begin{equation*}
\begin{split}
& \int_{\R^3_v} dv K_{\alpha,\beta} (v,z) \\
 & = \int_{\R^3_h} dh \un_{ | h| \geq \delta/2}
  \sep{ 1+ |z|^2 + 2 z.h  + |h|^2 }^{\alpha/2} \seq{z}^\beta
  \sep{ 1+ |z|^2 - |z.{h\over{|h|}} |^2 }^{(1+\gamma +2s)/2} \\
  &  \ \ \ \ \ \ \ \ \ \ \ \ \ \ \ \ \ \ \ \ \sep{ e^{- |z.{h\over{|h|}}  |^2 -2z.h -|h|^2}  e^{- |z.{h\over{|h|}}  |^2}}^{1/4}   {1\over{|h|^{3+2s}}}.
\end{split}
\end{equation*}
%

We shift to spherical coordinates (along axis w.r.t $z$) ($h=r\omega$) to get
\begin{equation*}
\begin{split}
& \int_{\R^3_v} dv K_{\alpha,\beta} (v,z) \\
& = \int^\pi_0 \int^\infty_\delta d\phi \sin \phi dr \seq{z}^\beta (1+|z|^2 + 2|z| r \cos \phi + r^2)^{\alpha/2} (1+ |z|^2 - |z|^2 \cos^2 \phi )^{(1+\gamma +2s)/2}\\
 &  \ \ \ \ \ \ \ \ \ \ \ \ \ \ \ \ \ \ \ \ \sep{ e^{- |z|^2\cos^2\phi -2|z|r\cos\phi  -r^2}  e^{- |z|^2 \cos^2\phi} }^{1/4}  {1\over{r^{1+2s}}}.
\end{split}
\end{equation*}
Set $t=\cos\phi$ to get
\begin{equation*}
\begin{split}
& \int_{\R^3_v} dv K_{\alpha,\beta} (v,z) \\
&  = \int^1_{-1} dt  \int^\infty_\delta  dr \seq{z}^{\beta}(1+|z|^2 + 2|z| r t + r^2)^{\alpha/2} (1+ |z|^2 - |z|^2 t^2 )^{(1+\gamma +2s)/2} \\
&  \qquad \qquad \qquad  \sep{ e^{- |z|^2t^2 -2|z|rt  -r^2}  e^{- |z|^2 t^2} }^{1/4}  {1\over{r^{1+2s}}}.
\end{split}
\end{equation*}
We note again that if $|z|\leq 1$, this integral is bounded uniformly. We therefore assume in the following that $|z| \geq 1$ and change variable $t' =|z|t$ to deduce that
\begin{equation*}
\begin{split}
& \int_{\R^3_v} dv K_{\alpha,\beta} (v,z) \\
&  = |z|^{-1} \int^{|z|}_{-|z|} dt  \int^\infty_\delta  dr \seq{z}^{\beta}(1+|z|^2 + 2 r t + r^2)^{\alpha/2} (1+ |z|^2 - t^2)^{(1+\gamma +2s)/2} \\
 &  \qquad \qquad \qquad e^{- (t+r)^2/4}  e^{- t^2/4}   {1\over{r^{1+2s}}} \\
& \lesssim  |z|^{-1} \int^{|z|}_{-|z|} dt  \int^\infty_\delta  dr \seq{z}^{\beta}
 (1+ |z|^2 - t^2)^{(1+\gamma +2s+\alpha)/2}  \seq{r+t}^{|\alpha|} e^{- (t+r)^2/4}  e^{- t^2/4}   {1\over{r^{1+2s}}},
\end{split}
\end{equation*}
where the last inequality is a consequence of Peetre's inequality \eqref{peetre}.
With exactly the same argument as before for the integration w.r.t. $r$, for small $\delta$, we obtain
\begin{equation*}
\begin{split}
& \int_{\R^3_v} dv K_{\alpha,\beta} (v,z) \lesssim \delta^{-1-2s} \seq{z}^{\alpha + \beta + \gamma+2s}
\end{split}
\end{equation*}
and thus
\begin{equation} \label{intv}
\begin{split}
& \int_{\R^3_v} dv K_{\alpha,\beta} (v,z) \lesssim \delta^{-1-2s}
\end{split}
\end{equation}
when $\alpha+\beta+\gamma +2s\leq 0$.
From (\ref{intz}) and
(\ref{intv}),  we use Schur's Lemma to obtain conclusion (i) in
Lemma  \ref{lb1da}.

\bigskip
(ii) Now we prove the second estimate about the commutator in  Lemma \ref{lb1da}.  Using the $\sigma$ representation between $v,v_*$ and $v', v_*'$, (see Figure 2 in Subsection 5.1 of Appendix),  we have, for $\theta\in ]0,\pi[,$
\begin{eqnarray*}
\abs{v-v_*}=\frac{\abs{v'-v_*}}{\cos{\theta\over 2}}\leq \sqrt 2\abs{v'-v_*}\leq   \sqrt 2\abs{v_\lambda-v_*},
\end{eqnarray*}
where
\[v_\lambda=v+\lambda (v'-v), \quad \lambda\in [0,1]. \]
As a result,
\begin{eqnarray*}
\seq{v}\leq \seq{v-v_*}+\seq{v_*} \leq \sqrt 2 \seq{v_\lambda-v_*}+\seq {v_*}\leq(1+ \sqrt 2 )\seq{v_\lambda}\seq {v_*},
\end{eqnarray*}
which along with the estimate
\[  \seq{v_\lambda}\leq (1+\sqrt 2) \seq{v}\seq{v_*} \]
due to the fact that $\abs{v'-v}=\abs{v-v_*}\sin{\theta\over 2}\leq \frac{\sqrt 2}{2}\abs{v-v_*},$  implies
\[
 \forall ~\kappa\in\mathbb R,\quad  \seq{v_\lambda}^\kappa \lesssim  \seq{v}^\kappa \seq{v_*}^{\abs\kappa }.
\]
Therefore, we have
\begin{eqnarray*}
\abs{\seq{v'}^{\tilde \beta} -\seq{v}^{\tilde \beta}}&\lesssim& \int_0^1 \comii{v_\lambda}^{\tilde \beta-1}
d\lambda \abs{v-v'} \\
&\lesssim& \comii{v}^{\tilde \beta-1}\comii{v_*}^{\abs{\tilde \beta-1}}
 \abs{v-v'}.
\end{eqnarray*}
Then
\begin{eqnarray*}
&& \abs{\seq{v}^{\tilde \alpha} \com{\bar \lll_{1,\delta ,a},~ \seq{v}^{\tilde\beta}}
f}
= \abs{\seq{v}^{\tilde \alpha }\iint dv_* d\sigma B \tilde\phi_\delta (v'-v)
  (\mu_*)^{1/2}  (\mu'_*)^{1/2} \inner{\seq{v'}^{\tilde\beta }-\seq{v}^{\tilde \beta}
  }f'}\\
   &&\lesssim \seq{v}^{\tilde \alpha+\tilde\beta-1} \iint dv_* d\sigma B \tilde\phi_\delta (v'-v)  (\mu_*)^{1/2}  (\mu'_*)^{1/2}\comii{v_*}^{\abs{\tilde \beta
-1}}
 \abs{v-v'}\abs{ f'}.
\end{eqnarray*}
Using  Carleman's representation (see the appendix for instance) shows that
\begin{eqnarray*}
&& \abs{\seq{v}^{\tilde \alpha} \com{\bar \lll_{1,\delta ,a},~ \seq{v}^{\tilde\beta}}
f}\\
 & \lesssim & \seq{v}^{\tilde \alpha+\tilde\beta-1}
 \int_h dh\int_{E_{0,h}} d\alpha \un_{|\alpha | \geq | h|} \un_{ | h| \geq {\delta\over 2}}    \mu^{1\over 2} (\alpha +v) \mu^{1\over 2} (\alpha +v -h) \\
  & & \qquad \qquad \qquad \qquad \qquad \qquad \times {\seq{\alpha+v }}^{\abs{\tilde\beta-1}} { |\alpha + h|^{1+\gamma + 2s} \over{| h|^{2+2s}} }|f(v-h)| \\
 & \lesssim & \seq{v}^{\tilde \alpha+\tilde\beta-1}  \int_h dh\int_{E_{0,h}} d\alpha  \un_{ | h| \geq {\delta\over 2}}    \mu^{1 \over 4} (\alpha +v) \mu^{1\over 4} (\alpha +v -h)    { |\alpha|^{1+\gamma + 2s} \over{| h|^{2+2s}} }  |f(v-h)|.
\end{eqnarray*}
The last term is quite similar as the one on the right hand side of \reff{expba}, with $\alpha$ and  $\beta$ there replaced  respectively by $\tilde\alpha+\tilde\beta-1$ and  $0$,  and $\mu^{1/2}$, $\abs h^{-(3+2s)}$ there replaced respectively  by $\mu^{1/4}$, $\abs h^{-(2+2s)}$.  Then repeating the arguments after \reff{expba}, we conclude
\begin{eqnarray*}
\abs{\seq{v}^{\tilde \alpha} \com{\bar \lll_{1,\delta ,a},~ \seq{v}^{\tilde\beta}}
f}\lesssim \int \tilde K_{\tilde \alpha, \tilde \beta} (v,z) \abs{f}(z)dz
\end{eqnarray*}
with
\begin{multline*}
\tilde K_{\tilde \alpha, \tilde \beta} (v,z) = \un_{ | v-z| \geq \delta/2} \seq{v}^{\tilde \alpha+\tilde\beta-1}   \seq{S(v-z) v}^{1+\gamma +2s}  {1\over{| v-z|^{2+2s}}} \\ \sep{ e^{2( |S(v-z) v|^2 - |v|^2)+ 2 v.(v-z) -|v-z|^2}}^{1/4}.
\end{multline*}
Arguing as for the analysis of $K_{\alpha,\beta}$ in (i), with $\alpha=\tilde\alpha+\tilde\beta-1$ and $\beta=0$,   we obtain a similar estimate as \reff{inttoprecise}, that is,
\begin{eqnarray*}
\int_{\R^3_z} dz \tilde K_{\tilde \alpha, \tilde \beta} (v,z) \lesssim \delta^{-2s} \seq{v}^{(\tilde \alpha+\tilde\beta-1)-1} \int_{0}^{|v|}  dt ( 1  +|v|^2 -  t^2  )^{(1+\gamma +2s )/2}   {1\over \seq{t}^{2s}}.
\end{eqnarray*}
It's clear that
\[
\int_{\R^3_z} dz \tilde K_{\tilde \alpha, \tilde \beta} (v,z) \lesssim  \delta^{-2s}
\]
for all $v$ such that $\abs v\leq 1$.
%
%

We can therefore assume $\abs v\geq1$ in the following. We split the integration into three parts as follows.  First
\begin{eqnarray*}
\int_{0}^{1/2}  dt ( 1  +|v|^2 -  t^2  )^{(1+\gamma +2s )/2}   {1\over \seq{t}^{2s}}\lesssim \seq{v}^{1+\gamma +2s }.
\end{eqnarray*}
Next, for any $\eps_0>0,$
\begin{eqnarray*}
\int_{1/2}^{\abs v/2}  dt ( 1  +|v|^2 -  t^2  )^{(1+\gamma +2s )/2}   {1\over \seq{t}^{2s}} &\lesssim & \seq{v}^{1+\gamma +2s }\int_{1/2}^{\abs v/2}    {{dt  }\over t^{2s}}\\
&\lesssim &
\left\{
\begin{array}{lll}
 \seq{v}^{1+\gamma +2s }\inner{\abs v^{-2s+1}+1},\quad
 s\neq 1/2\\[5pt]
\seq{v}^{1+\gamma +2s }\inner{\ln \abs v+1}\lesssim \seq{v}^{2+\gamma+\eps_0},\quad
 s= 1/2
\end{array}
\right.\\
&\lesssim &
 \seq{v}^{2+\gamma+\eps_0}  + \seq{v}^{1+\gamma +2s }
\end{eqnarray*}
Finally,  repeating the arguments used to get the estimate
\reff{inttopreciselarge}, we have
\begin{eqnarray*}
\int_{\abs v/2}^{\abs v
}  dt ( 1  +|v|^2 -  t^2  )^{(1+\gamma +2s )/2}   {1\over \seq{t}^{2s}} &\lesssim &
 \seq{v}^{-2s }\abs v^{-1} \int_0^{\abs v^2/2} \inner{1+\lambda}^{(1+\gamma+2s)/2} d\lambda \\
  &\lesssim &
 \seq{v}^{2+\gamma}.
\end{eqnarray*}
Combining these inequalities gives, for $\abs v\geq 1,$
\begin{eqnarray*}
\int_{\R^3_z} dz \tilde K_{\tilde \alpha, \tilde \beta} (v,z) &\lesssim &\delta^{-2s} \seq{v}^{(\tilde \alpha+\tilde\beta-1)-1} \inner{\seq{v}^{2+\gamma+\eps_0}  + \seq{v}^{1+\gamma +2s  } }\\
 &\lesssim & \delta^{-2s} \seq{v}^{\tilde
   \alpha+\tilde\beta+\gamma+\eps_0} +\delta^{-2s} \seq{v}^{\tilde
   \alpha+\tilde\beta+\gamma+2s-1} .
\end{eqnarray*}
Then choosing $\eps_0=s$ and using the assumption that
$\tilde\alpha+\tilde\beta+\gamma+s\leq 0,$ we conclude
\begin{eqnarray*}
\int_{\R^3_z} dz \tilde K_{\tilde \alpha, \tilde \beta} (v,z)
 \lesssim  \delta^{-2s}.
\end{eqnarray*}
Similarly as in (i), we can show that
\begin{eqnarray*}
\int_{\R^3_z} dv \tilde K_{\tilde \alpha, \tilde \beta} (v,z)
 \lesssim  \delta^{-2s}.
\end{eqnarray*}
Then Schur's Lemma applies and this completes  the proof of conclusion (ii)  in Lemma \ref{lb1da}.
\fin

\subsubsection{Study of $\lll_{1,4,\delta}$}

 \begin{lem} \label{l14d}  For all $f\in\mathcal S(\mathbb R_v^3)$, we have
$$
\|  \lll_{1,4,\delta} f  \|^2\lesssim  \delta^{2-2s} \| <v>^{\gamma +2s} f \|,
$$
and $\lll_{1,4,\delta}$ commutes with the multiplication by $\seq{v}^\alpha$ for all $\alpha \in \R$.
\end{lem}

 \preuve
The last assertion is again trivial since $\lll_{1, 4,\delta}$ is a multiplication operator.
Using the formula $ 2a(b-a) = b^2 - a^2 - (b-a)^2$, we get
\begin{equation*}
\begin{split}
\lll_{1,4,\delta} f & = {1\over 2} f  \iint dv_* d\sigma B  \phi_\delta (v'-v)   \sep{ (\mu_*) - (\mu'_*)}  \\
& \ \ \ \ \ \ \ \ \ -{1\over 2} f  \iint dv_* d\sigma B  \phi_\delta (v'-v)   \sep{ (\mu_*)^{1/2} - (\mu'_*)^{1/2}}^2  \\
&  = -\frac{1}{2}\mathcal L_{1,3,\delta} f + D(v) f.
\end{split}
\end{equation*}
It suffices to estimate $D(v)$ in view of Lemma \ref{L1easy}.
To do so we essentially follow the same process, except that we don't need to use a symmetrizing argument to kill higher singularities. We write
\begin{equation*}
\begin{split}
|D(v)|
 & = \frac{1}{2}  \int_{\R^3_h}dh\int_{E_{0,h}} d\alpha \un_{|\alpha | \geq | h|} \phi_\delta (h)
 \sep{ \mu^{1/2}
 (\alpha +v -h) -\mu^{1/2} (\alpha +v) }^2 { | \alpha +h |^{1+\gamma +2s} \over{| h|^{3+2s}}} \\
 & \lesssim  \int_{\R^3_h}dh\int_{E_{0,h}} d\alpha \un_{|\alpha | \geq | h|} \phi_\delta (h)
 \sep{ e^{ (\alpha+v)\cdot h -h^2/2} -1}^2 \mu (\alpha +v) { | \alpha +h |^{1+\gamma +2s}  \over{| h|^{3+2s}}} \\
  & \lesssim  \int_{\R^3_h}dh\int_{E_{0,h}} d\alpha \un_{|\alpha | \geq | h|} \phi_\delta (h)
  \mu^{1/2} (\alpha +v) { | \alpha |^{1+\gamma +2s}  \over{| h|^{1+2s}}} \\
  &  \lesssim \delta^{2-2s}\seq{v}^{\gamma+ 2s},
 \end{split}
\end{equation*}
following the same arguments as before.
From the estimates on $\mathcal L_{1,3,\delta}$ and $D(v)$, the proof is complete. \fin

\section{Pseudo-differential parts}
In this section we deal with the remaining parts of $\lll_1$, namely:

- a multiplicative operator $\bar \lll_{1,\delta,b}$;

- the principal term
$\lll_{1,2,\delta}$ which will appear to be of pseudo-differential
type;

- and the term $\lll_{1,1, \delta}$ which is also of pseudo-differential type but with lower order (and we therefore call it subprincipal).

\noindent Our goal in this section is to prove Proposition \ref{estaa} about the behavior of these pseudo-differential parts of $\lll$.

\noindent In the following, we keep the notation for  $\phi_\delta$, the
positive compactly supported function equal to 1 in a
$\delta$-neighborhood of $0$ as introduced previously in the
definitions of the operators, and let $E_{0,\omega}=\omega^\perp$ for the hyperplane containing $0$ and orthogonal to $\omega$. We study each
operator separately. Proposition \ref{estaa} will be obtained as a
direct consequence of Proposition \ref{am2} and  Proposition \ref{ap2} below and Definition \ref{defaapm}.

\subsection{Study of the principal term  $\lll_{1,2,\delta}$}

Recall that
$$ \lll_{1,2,\delta} f =
 \iint dv_* d\sigma B \phi_\delta (v'-v)\mu'_* \sep{ f' -f  }
$$
where
$$
B(v,\sigma) =   |v-v_*|^{\gamma} b \sep{ \seq{ \frac{v-v_*}{|v-v_*|}, \sigma}}.
$$
This will appear to be a genuine pseudo-differential operator of order $2s$ for which we can control the weights. Namely one has

 \begin{prop} \label{ap2} We can write
 \begin{equation*}
 \lll_{1,2,\delta} f =  - a_p(v, D_v) f,
 \end{equation*}
 where $a_p$ is a real symbol in $(v, \eta)$ (see \eqref{equa} below for the definition of $a_p$) satisfying:
\begin{enumerate}[i)]
\item there exists  $C >0 $ such that for all $0<\kappa<1$,
\begin{eqnarray} \label{minmaxap}
\begin{split}
 &C^{-1} \delta^{2-2s} \sep{ -\kappa \seq{v}^{\gamma+2s} + \kappa \seq{v}^{\gamma} (|\eta|^{2s} + |\eta\wedge v|^{2s} ) } \\ &\qquad\quad  \qquad \leq a_p(v, \eta) \leq  C \seq{v}^{\gamma} (1+ |\eta|^{2s} + |\eta\wedge v|^{2s} );
 \end{split}
 \end{eqnarray}
 \item
$ a_p \in S \sep{ \seq{v}^{\gamma} (1+ |v|^{2s} +  |\eta|^{2s} + |\eta\wedge v|^{2s} ) , \Gamma} .$ Recall $\Gamma =\abs{dv}^2+\abs{d\eta}^2$ is the flat metric.
 \end{enumerate}
 \end{prop}

    \remark The first estimates in \reff{minmaxap} explain why we don't have regularity estimate for the Boltzmann equation with angular cutoff,  since it corresponds to the case $\delta\rightarrow 0$ and thus we lose the regularity operator $\comii v^\gamma \comii {D_v}^{2s}+\comii v^\gamma \comii {D_v\wedge v}^{2s}.$  Observe we exclude the case $s=1$, and this corresponds the Landau equation, which is the grazing limit of Boltzmann equation without angular cutoff and still admits the diffusion structure.

 \preuve From the expression of  $\lll_{1,2,\delta}$, using Carleman's transformation  as in previous arguments and as in \cite{Ale99} (see also the Appendix), we get
  \begin{equation*}
  \begin{split}
  \lll_{1,2,\delta} f & = \int_{\R^3_h} dh \int_{E_{0,h}} d\alpha \tb(\alpha, h) \un_{|\alpha |\geq | h |} \phi_\delta (h)  \mu (\alpha +v) | \alpha +h|^{1+\gamma +2s} \sep{ f(v-h) - f(v)}{1\over{| h|^{3+2s}}},
\end{split}
\end{equation*}
{\red where $\tb(\alpha, h)$ is a function of $\alpha$ and $h$ which is bounded from below and above by positive constants, and satisfies that
$ \tb(\alpha, h) = \tb(\pm \alpha, \pm h)$.}

 This integral is typically undefined for large values of $s$, and we have to use its symmetrized version in order to give a meaning in the principal value  sense:  for this purpose, we change $h$ to $-h$ and add the two expressions to obtain
 \begin{equation*}
  \begin{split}
  \lll_{1,2,\delta} f & = {1\over 2} \int_h dh \int_{E_{0,h}} d\alpha \tb \un_{|\alpha |\geq | h |} \phi_\delta (h)  \mu (\alpha +v) | \alpha +h|^{1+\gamma +2s} \\
  & \qquad \qquad \qquad \qquad \sep{ f(v-h) +f(v+h) - 2f(v)}{1\over{| h|^{3+2s}}} \\
 &  \defegal -a_p(v,D_v) f (v) \defegal -\int_{\R^3_\eta}  a_p(v,\eta ) \hat f (\eta )e^{i\eta .v} d\eta
 \end{split}
 \end{equation*}
 with
\begin{equation} \label{equa}
  \begin{split}
  a_p (v,\eta ) & \defegal  -{1\over 2} \int_{\R^3_h} dh \int_{E_{0,h}} d\alpha \tb \un_{|\alpha |\geq | h |} \phi_\delta (h)  \mu (\alpha +v) | \alpha +h|^{1+\gamma +2s} \\
 &   \qquad  \qquad  \qquad \qquad  \qquad \sep{ e^{-i\eta \cdot h} + e^{i\eta \cdot h}  - 2}{1\over{| h|^{3+2s}}} \\
 &  = \int_{\R^3_h} dh \int_{E_{0,h}} d\alpha \tb \un_{|\alpha |\geq | h |} \phi_\delta (h) \mu (\alpha +v) | \alpha +h|^{1+\gamma +2s} \\
 & \qquad  \qquad  \qquad \qquad  \qquad  \sep{ 1- \cos (\eta \cdot h)  }{1\over{| h|^{3+2s}}} .
 \end{split}
 \end{equation}

%
%
%
%

The non-negativity of $a_p(v,\eta)$ is clear and we shall now work on some properties of this symbol.
 First recall that on the support of the integrand, we have $|h| \leq
 \delta \leq 1$ and that $\alpha \perp h$, so that
$$ 0 \leq  a_p(v,\eta )   \lesssim \int_{\R^3h} dh \int_{E_{0,h}} \un_{|\alpha |\geq | h |} \un_{| h| \leq \delta} d\alpha \mu (\alpha +v) \seq{ \alpha}^{1+\gamma +2s} \sep{ 1- \cos (\eta \cdot h)  }{1\over{| h|^{3+2s}}}.
$$
Now we can shift to spherical coordinates $h= r\omega$, and (forgetting the truncation in $\alpha$) we get
$$ a_p(v,\eta )   \lesssim  \int^\delta_0\int_{S^2_\omega} dr d\omega \int_{E_{0,\omega}} d\alpha \mu (\alpha +v) \seq{ \alpha}^{1+\gamma +2s} \sep{ 1- \cos (r \eta . \omega)  }{1\over{r^{1+2s}}}.
$$
It is possible to integrate directly w.r.t.  $r$, and use the fact that
\begin{equation}\label{upbmarch30}
\int_0^\delta \sep{ 1- \cos (r \eta . \omega)  }{1\over{r^{1+2s}}} dr \leq C_{s} |\omega\cdot\eta|^{2s}.
\end{equation}
In fact, note that
$$
\int_0^\delta \sep{ 1- \cos (r \eta . \omega)  }{1\over{r^{1+2s}}} dr =\int_0^\delta  {{ 1- \cos (r \abs{\eta . \omega}) }\over{r^{1+2s}}} dr =| \omega\cdot\eta |^{2s} \int_0^{\delta |\omega\cdot\eta |} {{ 1- \cos r }\over{r^{1+2s}}} dr.
$$
Next, we choose a small constant $c$ such that $1-\cos r \gtrsim r^2$ if $r\leq c$.

If $ |\omega\cdot\eta | \geq c$, then we get
$$
\int_0^\delta \sep{ 1- \cos (r \eta . \omega)  }{1\over{r^{1+2s}}} dr \gtrsim |\omega\cdot\eta|^{2s} \int_0^{c\delta }\sep{ 1- \cos (r )  }{1\over{r^{1+2s}}} dr \gtrsim \delta^{2-2s} |\omega\cdot\eta|^{2s},
$$
while if $|\omega\cdot\eta | \leq c$, then we get
$$
\int_0^\delta \sep{ 1- \cos (r \eta . \omega)  }{1\over{r^{1+2s}}} dr  \gtrsim |\omega\cdot\eta |^2 \int^\delta_0 r^2 {1\over{r^{1+2s}}} dr \gtrsim \delta^{2-2s} |\omega\cdot\eta |^2.
$$
On the whole, we get
\begin{equation}\label{minimum}
\int_0^\delta \sep{ 1- \cos (r \eta . \omega)  }{1\over{r^{1+2s}}} dr \gtrsim \delta^{2-2s} \min\lbrace |\omega\cdot\eta |^2 , | \omega\cdot\eta|^{2s} \rbrace .
\end{equation}
This along with \eqref{upbmarch30} gives
\begin{equation}\label{min-max}
\delta^{2-2s} \min\lbrace |\omega\cdot\eta |^2 , | \omega\cdot\eta|^{2s} \rbrace \lesssim \int_0^\delta \sep{ 1- \cos (r \eta . \omega)  }{1\over{r^{1+2s}}} dr \leq C_s   | \omega\cdot\eta|^{2s}.
\end{equation}

\bigskip
Next, we deal with the upper bound on $a_p$. A crude estimate is enough and we get
\begin{equation} \label{typiquea}
 a_p(v,\eta )  \lesssim \int_{S^2_\omega}  d\omega \int_{E_{0,\omega}}  d\alpha \mu (\alpha +v) |\omega\cdot\eta|^{2s} \seq{ \alpha}^{1+\gamma +2s}.
\end{equation}
Splitting $ v = S(\omega) v + (\omega\cdot v) \omega$, we have
\begin{equation} \label{splittingorth}
|\alpha+ v |^2 = |\alpha + S(\omega) v + (\omega\cdot v) \omega |^2 = |\alpha + S(\omega) v|^2 + |(\omega\cdot v)|^2
\end{equation}
since $\alpha $ and $\omega$ are orthogonal. We can therefore write
$$ \mu (\alpha +v) =  (2\pi)^{-3/2} \sep{ e^{-|\alpha +S(\omega)v|^2} e^{-|(\omega\cdot v)|^2}}^{1/2}
$$
to get
\begin{equation} \label{typique0}
a_p(v,\eta )  \lesssim \int_{S^2_\omega}  d\omega \int_{E_{0,\omega}}  d\alpha  \sep{ e^{-|\alpha +S(\omega)v|^2} e^{-|(\omega\cdot v)|^2}}^{1/2}  |\omega\cdot\eta|^{2s} \seq{ \alpha}^{1+\gamma +2s}.
\end{equation}
Next, note that
\begin{equation*} 
\beta (v,\omega ) = \int_{E_{0,\omega}}  d\alpha  \sep{ e^{-|\alpha +S(\omega)v|^2} }^{1/2}  \seq{ \alpha }^{1+\gamma +2s} \sim <S(\omega) v>^{1+\gamma +2s}
\end{equation*}
and thus
\begin{equation} \label{+typiqueplus++}
a_p(v,\eta )  \lesssim \int_{S^2_\omega } d\omega  e^{-|(\omega\cdot v)|^2/2} <S(\omega) v>^{1+\gamma +2s} |\omega\cdot\eta|^{2s}.
\end{equation}

We introduce polar coordinates in a coordinate system where ${\mathbf i} = S(v) \eta/ |S(v)\cdot\eta|$, ${\mathbf k} = v/|v|$.

\begin{figure}[!t]
\centering
\begin{picture}(0,0)%
\includegraphics{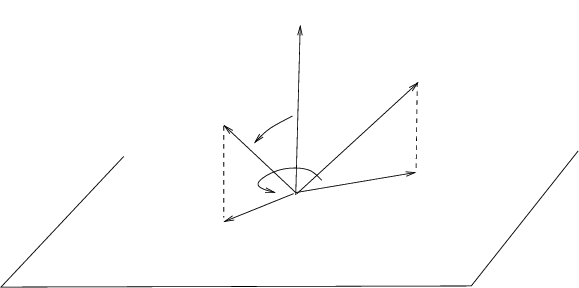}%
\end{picture}%
%
%
\setlength{\unitlength}{2368sp}%
\begingroup\makeatletter\ifx\SetFigFont\undefined%
\gdef\SetFigFont#1#2#3#4#5{%
  \reset@font\fontsize{#1}{#2pt}%
  \fontfamily{#3}\fontseries{#4}\fontshape{#5}%
  \selectfont}%
\fi\endgroup%
\begin{picture}(7719,3815)(469,-3958)
\put(3541,-1726){\makebox(0,0)[lb]{\smash{{\SetFigFont{7}{8.4}{\rmdefault}{\mddefault}{\updefault}{\color[rgb]{0,0,0}$\omega$}%
}}}}
\put(5686,-1261){\makebox(0,0)[lb]{\smash{{\SetFigFont{7}{8.4}{\rmdefault}{\mddefault}{\updefault}{\color[rgb]{0,0,0}$\eta$}%
}}}}
\put(4366,-316){\makebox(0,0)[lb]{\smash{{\SetFigFont{7}{8.4}{\rmdefault}{\mddefault}{\updefault}{\color[rgb]{0,0,0}$v$}%
}}}}
\put(5761,-2641){\makebox(0,0)[lb]{\smash{{\SetFigFont{7}{8.4}{\rmdefault}{\mddefault}{\updefault}{\color[rgb]{0,0,0}$S(v)\eta$}%
}}}}
\put(3601,-3136){\makebox(0,0)[lb]{\smash{{\SetFigFont{7}{8.4}{\rmdefault}{\mddefault}{\updefault}{\color[rgb]{0,0,0}$S(v)\omega$}%
}}}}
\put(7651,-3316){\makebox(0,0)[lb]{\smash{{\SetFigFont{7}{8.4}{\rmdefault}{\mddefault}{\updefault}{\color[rgb]{0,0,0}$i = S(v) \eta / |S(v)\eta|$}%
}}}}
\put(7651,-3601){\makebox(0,0)[lb]{\smash{{\SetFigFont{7}{8.4}{\rmdefault}{\mddefault}{\updefault}{\color[rgb]{0,0,0}$k=v/|v|$}%
}}}}
\put(4036,-1711){\makebox(0,0)[lb]{\smash{{\SetFigFont{7}{8.4}{\rmdefault}{\mddefault}{\updefault}{\color[rgb]{0,0,0}$\varphi$}%
}}}}
\put(4201,-2326){\makebox(0,0)[lb]{\smash{{\SetFigFont{7}{8.4}{\rmdefault}{\mddefault}{\updefault}{\color[rgb]{0,0,0}$\theta$}%
}}}}
\end{picture}%
\caption{{\red spherical} coordinates}
\label{polar}
\end{figure}
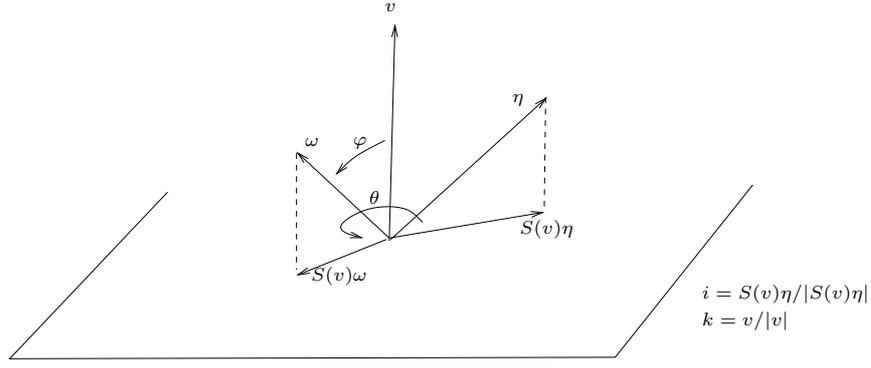

\noindent \bigskip In this system,  we note that $(\omega\cdot{\mathbf k}) = \cos(\varphi)$. Besides we have
$\eta = (\eta. {\mathbf k}) {\mathbf k} + S(v) \eta $ so that
\begin{equation*}
\begin{split}
\omega\cdot\eta & = (\eta\cdot{\mathbf k})({\mathbf k}\cdot\omega) + (S(v)\eta )\cdot  \omega \\
& = (\eta\cdot{\mathbf k})({\mathbf k}\cdot\omega) + (S(v)\eta )\cdot  (S(v) \omega) \\
& = (\eta\cdot{\mathbf k})({\mathbf k}\cdot\omega) + \sep{ {\mathbf i}\cdot (S(v) \omega)} |S(v) \eta| \\
& = \eta\cdot{\mathbf k}\cos(\varphi ) + |S(v) \eta| \sin(\varphi) \cos(\theta).
\end{split}
\end{equation*}
and in a similar way
$$
|S(\omega)v|^2  = |v|^2 - |(v.\omega)|^2 = |v|^2( 1 - \cos^2(\varphi)) = |v|^2 \sin^2(\phi).
$$
We therefore get
\begin{multline*}
a_p(v,\eta )  \lesssim \int_0^\pi d\varphi \int_0^{2\pi} d\theta   \sin(\varphi)   e^{-|v|^2 \cos^2(\varphi)} \sep{ 1 + |v|^2 \sin^2(\varphi)}^{(1+\gamma +2s)/2} \\
 |\eta\cdot{\mathbf k}\cos(\varphi ) + |S(v) \eta| \sin(\varphi) \cos(\theta)|^{2s}.
\end{multline*}
Setting $\cos \varphi = t$ in the preceding formula, we get
\begin{multline} \label{exactexpression}
a_p(v,\eta ) \lesssim \int_0^{2\pi} d\theta \int_0^1 dt     e^{-|v|^2 t^2} \sep{ 1 + |v|^2 (1-t^2)}^{(1+\gamma +2s)/2} \\
 |\eta\cdot{\mathbf k}t + |S(v) \eta| \sqrt{1-t^2} \cos(\theta)|^{2s}.
\end{multline}
If we  bound roughly $1-t^2$ and $\cos(\varphi)$ by $1$ and use the estimates  that
\begin{eqnarray*}
 e^{-|v|^2 t^2} \sep{ 1 + |v|^2 (1-t^2)}^{(1+\gamma +2s)/2}  \lesssim  e^{-|v|^2 t^2} \sep{ 1 + |v|^2}^{(1+\gamma +2s)/2}
\end{eqnarray*}
for $1+\gamma+2s\geq 0$ or $0\leq t\leq 1/2$, and that
\begin{equation*}
\begin{split}
 e^{-|v|^2 t^2} \sep{ 1 + |v|^2 (1-t^2)}^{(1+\gamma +2s)/2}
  \lesssim e^{-|v|^2 t^2}  \lesssim e^{-|v|^2 t^2/2} \sep{1+ v^2}^{(1+\gamma +2s)/2}
 \end{split}
\end{equation*}
for $1+\gamma+2s<0$ and uniformly w.r.t. $1/2\leq t\leq 1$,
\color{black}


then we get
$$
a_p(v,\eta )  \lesssim  \int_0^{2\pi} d\theta \int_0^{1} dt   \ \    e^{-|v|^2 t^2/2} \sep{ 1 + |v|^2}^{(1+\gamma +2s)/2}
 \sep{ |\eta\cdot{\mathbf k} t|^{2s}   + |S(v) \eta|^{2s}  }.
$$
If we set  $y = |v|t$,  we get
\begin{equation} \label{estt0}
\begin{split}
a_p(v,\eta ) &  \lesssim   \frac{1}{|v|}\seq{ v}^{1+\gamma +2s} \int_0^{2\pi} d\theta \int_0^{|v|} dy   \ \   e^{-y^2/2}  \sep{
 |\eta\cdot{\mathbf k} |^{2s} \frac{y^{2s}}{|v|^{2s}}  + |S(v) \eta|^{2s}  } \\
&  \lesssim   \frac{1}{|v|}\seq{ v}^{1+\gamma +2s}    \sep{
 |\eta\cdot{\mathbf k} |^{2s} \frac{1}{|v|^{2s}}  + |S(v) \eta|^{2s}  } \\
 & \lesssim    \frac{\seq{ v}^{1+\gamma +2s}}{|v|^{1+ 2s}} |\eta|^{2s}     + \frac{\seq{ v}^{1+\gamma +2s}}{|v|} |S(v) \eta|^{2s}.
\end{split}
\end{equation}
For $|v| \geq 1$, we therefore get
\begin{equation*}
\begin{split}
a_p(v,\eta ) &  \lesssim      \seq{ v}^{\gamma} |\eta|^{2s}     + \seq{ v}^{\gamma +2s} |S(v) \eta|^{2s},
\end{split}
\end{equation*}
and thus
 \begin{equation*} 
\begin{split}
a_p(v,\eta ) &  \lesssim    \seq{ v}^{\gamma} \sep{ |\eta|^{2s}     +  |v \wedge \eta|^{2s}  },
\end{split}
\end{equation*}
since $|v \wedge \eta| = |v| |S(v) \eta|$.
For $|v| \leq 1$, a rough estimate gives directly $|a(v,\eta )| \leq \seq{\eta}^{2s}$ so that the preceding estimate is also true. The proof of the upper bound is complete.

\bigskip
Now we  deal with the lower bound. To this end, we shall use the formula (\ref{equa})
\begin{equation*}  
 a_p(v,\eta )  = \int_{\R^3_h} dh \int_{E_{0,h}} d\alpha \tb \phi_\delta(h) \un_{|\alpha| \geq |h|}  \mu (\alpha +v) | \alpha +h|^{1+\gamma +2s} \sep{ 1- \cos (\eta \cdot h)  }{1\over{| h|^{3+2s}}}.
 \end{equation*}
As we want a lower bound  we can  restrict the integration
range to $\set{|\alpha| \geq 10}$ since the integrand is non negative.
We use also the facts that $\tb $
is bounded from below by a positive constant and that $|\alpha+h| \sim |\alpha|$ since $\alpha \perp
h$ and $\abs h\leq \abs\alpha$ in the preceding integral.
Therefore, we have
$$
a_p(v,\eta ) \gtrsim \int_{\R^3_h} dh \int_{E_{0,h}} \phi_\delta(h) \un_{|\alpha| \geq 10} d\alpha \mu (\alpha +v) \comii{\alpha}^{1+\gamma +2s} \sep{ 1- \cos (\eta \cdot h)  }{1\over{| h|^{3+2s}}}.
$$
We can use some of the previous computations, and from \reff{minimum}-\reff{min-max} we get as in (\ref{typique0}),
$$
a_p(v,\eta )   \gtrsim  \delta^{2-2s} \int_{S^2_\omega}  d\omega \int_{E_{0,\omega}}  d\alpha \un_{|\alpha| \geq 10}  e^{-|\alpha +S(\omega)v|^2/2} e^{-|(\omega\cdot v)|^2/2} { \min \lbrace | \omega\cdot\eta |^2 ,|\omega\cdot\eta|^{2s} \rbrace }\seq{ \alpha}^{1+\gamma +2s} .
$$
Note that
$$ 
\beta_{10} (v, \omega) \defegal \int_{E_{0,\omega}}  d\alpha  1_{|\alpha| \geq 10} e^{-|\alpha +S(\omega)v|^2/2}   \seq{ \alpha }^{1+\gamma +2s} \sim <S(\omega) v>^{1+\gamma +2s} .
$$ 


Therefore
\begin{equation} \label{typique}
a_p(v,\eta )  \gtrsim \delta^{2-2s} \int_{S^2_\omega}  d\omega  e^{-|(\omega\cdot v)|^2/2} <S(\omega) v>^{1+\gamma +2s}{\min \lbrace | \omega\cdot\eta |^2 ,|\omega\cdot\eta|^{2s} \rbrace } .
\end{equation}

We now consider an arbitrary real $0 < \kappa < 1$. Using the fact that
$$
{\min \lbrace | \omega\cdot\eta |^2 ,|\omega\cdot\eta|^{2s} \rbrace } \geq |\omega\cdot\eta|^{2s}-1,
$$
and that the right member in (\ref{typique}) is non-negative, we get that
\begin{equation} \label{typiquebis}
\begin{split}
a_p(v,\eta )  & \gtrsim   \kappa \delta^{2-2s}   \int_{S^2_\omega}  d\omega  e^{-|(\omega\cdot v)|^2/2} <S(\omega) v>^{1+\gamma +2s}{\min \lbrace | \omega\cdot\eta |^2 ,|\omega\cdot\eta|^{2s} \rbrace }  \\
 & \gtrsim \kappa \delta^{2-2s}   \int_{S^2_\omega}  d\omega  e^{-|(\omega\cdot v)|^2/2} <S(\omega) v>^{1+\gamma +2s}(|\omega\cdot\eta|^{2s} - 1) \\
& \gtrsim \kappa \delta^{2-2s}   \int_{S^2_\omega}  d\omega  e^{-|(\omega\cdot v)|^2/2} <S(\omega) v>^{1+\gamma +2s} |\omega\cdot\eta|^{2s} \\
& \qquad \qquad \qquad - \kappa \delta^{2-2s}   \int_{S^2_\omega}  d\omega  e^{-|(\omega\cdot v)|^2/2} <S(\omega) v>^{1+\gamma +2s} \\
& \defegal \kappa \delta^{2-2s}   a_{pp} -\kappa \delta^{2-2s}   a_{pr}.
\end{split}
\end{equation}

We split the study of the two terms $a_{pp}$ and $ a_{pr}$.
For $a_{pp}$, we can use previous computations yielding to
(\ref{+typiqueplus++}). More precisely, we have
\begin{multline} \label{exactexpression2}
a_{pp}(v,\eta )   =  \int_0^{2\pi} d\theta \int_0^1 dt     e^{-|v|^2 t^2} \sep{ 1 + |v|^2 (1-t^2)}^{(1+\gamma +2s)/2} \\
 \abs{\eta\cdot{\mathbf k}t + |S(v) \eta| \sqrt{1-t^2} \cos(\theta)}^{2s}.
\end{multline}
Now an easy remark is that the symbol $a_{pp}$ has the following parity properties:
$$
a_{pp}(\pm v,\pm \eta) = a_{pp}(v,\eta ).
$$
We can therefore assume that $\eta\cdot{\mathbf k} \geq 0$ in all the
computations.  Moreover we can restrict
the above integration to the following subsets
\begin{equation} \label{segm}
t \in [0, \sqrt{3}/2], \ \ \ \ \theta \in [0,\pi/3],
\end{equation}
which implies that all terms inside the absolute value
$$|\eta\cdot{\mathbf k}t + |S(v) \eta| \sqrt{1-t^2} \cos(\theta)|  $$
are non-negative. We therefore get, when (\ref{segm}) is fulfilled, that
\begin{equation*}
\begin{split}
 \sep{ 1 + |v|^2 (1-t^2)}^{(1+\gamma +2s)/2} & \geq \sep{ 1 + \frac{|v|^2}{4} }^{(1+\gamma +2s)/2} \geq c_{s,\gamma} \comii{ v }^{1+\gamma +2s} \end{split}
 \end{equation*}
 and
\begin{equation*}
\begin{split}
 |\eta\cdot{\mathbf k}t + |S(v) \eta| \sqrt{1-t^2} \cos(\theta)|^{2s} & \geq 4^{-2s} |\eta\cdot{\mathbf k}t + |S(v) \eta| |^{2s} \\
 & \geq c_s \sep{ |\eta\cdot{\mathbf k}t|^{2s} + |S(v) \eta|^{2s} } .
 \end{split}
 \end{equation*}
Therefore putting the above estimate into \eqref{exactexpression2} gives
$$
a_{pp}(v,\eta )  \gtrsim  \int_0^{\pi/3} d\theta \int_0^{\sqrt{3}/2} dt   \ \    e^{-|v|^2 t^2} \comii{ v}^{1+\gamma +2s}
 \sep{ |\eta\cdot{\mathbf k} t|^{2s}   + |S(v) \eta|^{2s}  }.
$$
As in the case of the upper bound,  we set  $y = |v|t$,  and get for $|v| \geq 1$ that
\begin{equation*}
\begin{split}
a_{pp}(v,\eta ) &  \gtrsim  \frac{1}{|v|}\seq{ v}^{1+\gamma +2s} \int_0^{\pi/3} d\theta \int_0^{\sqrt{3}|v|/2} dy   \ \   e^{-y^2}  \sep{
 |\eta\cdot{\mathbf k} |^{2s} \frac{y^{2s}}{|v|^{2s}}  + |S(v) \eta|^{2s}  } \\
 &  \gtrsim  \frac{1}{|v|}\seq{ v}^{1+\gamma +2s} \int_0^{\pi/3} d\theta \int_0^{\sqrt{3}/2} dy   \ \   e^{-y^2}  \sep{
 |\eta\cdot{\mathbf k} |^{2s} \frac{y^{2s}}{|v|^{2s}}  + |S(v) \eta|^{2s}  } \\
&  \gtrsim \frac{1}{|v|}\seq{ v}^{1+\gamma +2s}    \sep{
 |\eta\cdot{\mathbf k} |^{2s} \frac{1}{|v|^{2s}}  + |S(v) \eta|^{2s}  } \\
& \gtrsim \sep{   \seq{ v}^{\gamma} |\eta|^{2s}     + \seq{ v}^{\gamma +2s} |S(v) \eta|^{2s}  },
\end{split}
\end{equation*}
where in  the last inequality we use that $\eta\cdot{\mathbf k} \geq
0$ and the fact that if  $\eta\cdot{\mathbf k} \leq  |\eta|/2 $ then
$$
|S(v) \eta| \geq \sqrt{3} |\eta|/2.
$$
Since $|v \wedge \eta| = |v| |S(v) \eta|$ we get for $|v| \geq 1$ the desired result
 \begin{equation} \label{app}
\begin{split}
a_{pp}(v,\eta ) &  \gtrsim   \seq{ v}^{\gamma} \sep{ |\eta|^{2s}     +  |v \wedge \eta|^{2s}  } .
\end{split}
\end{equation}
For $|v| \leq 1$,  a direct check,  without the change of
variables $\abs v t\rightarrow y$,  gives
\begin{eqnarray*}
a_{pp}(v,\eta )  \gtrsim  \int_0^{\pi/3} d\theta \int_0^{\sqrt{3}/2} dt      e^{-t^2}
 \sep{ |\eta\cdot{\mathbf k} t|^{2s}   + |S(v) \eta|^{2s}  }&\gtrsim&
 |\eta\cdot{\mathbf k} |^{2s}   + |S(v) \eta|^{2s}\\
&\gtrsim& |\eta|^{2s}   + |v\wedge \eta|^{2s}.
\end{eqnarray*}
So the preceding estimate  \reff{app} is also true for $\abs v\leq 1$.

For the remainder term in (\ref{typiquebis}), we can use similar computations as the ones done for the upper bound for $a_p$, and we easily get
$$
a_{pr} \lesssim \seq{v}^{\gamma + 2s}.
$$
 Putting this estimate and (\ref{app}) together into (\ref{typiquebis}) completes the proof of the lower bound in (\ref{minmaxap}).

\bigskip Now we deal with estimates on the  derivatives in $\eta $ and $v$ of $a_p$.
 Recall that
 \begin{equation*}  
  \begin{split}
  a_p (v,\eta )
 &  = \int_{\R^3_h} dh \int_{E_{0,h}} d\alpha \tb \un_{|\alpha |\geq | h |} \phi_\delta(h) \mu (\alpha +v) | \alpha +h|^{1+\gamma +2s} \sep{ 1- \cos (\eta \cdot h)  }{1\over{| h|^{3+2s}}}
 \end{split}
 \end{equation*}
which is clearly smooth with respect to $v$ and $\eta$. Let us
consider for $\nu_1, \nu_2\in \N^3$ the derivative
 \begin{equation*}
  \begin{split}
 \D^{\nu_1}_v \D^{\nu_2}_\eta a_p (v,\eta )
 &  = \int_h dh \int_{E_{0,h}} d\alpha \tb \un_{|\alpha |\geq | h |} \phi_\delta(h) \sep{ \D^{\nu_1}_v \mu (\alpha +v)} | \alpha +h|^{1+\gamma +2s} \\
  & \ \ \ \ \ \ \ \ \ \ \ \ \ \ \ \
   \ \ \ \ \ \ \ \ \ \ \ \sep{  \D^{\nu_2}_\eta \sep{ 1- \cos (\eta \cdot h)  }} {1\over{| h|^{3+2s}}}.
 \end{split}
 \end{equation*}
Setting again $h= r\omega$, and (forgetting the truncation in $\alpha$) we get
 \begin{equation}
  \begin{split} \label{equa3}
\abs{  \D^{\nu_1}_v \D^{\nu_2}_\eta a_p(v,\eta ) }   \lesssim
\int_0^\delta \int_{S^2_\omega} dr d\omega \int_{E_{0,\omega}} d\alpha | \D^{\nu_1}_v \mu (\alpha +v) |
 \seq{ \alpha}^{1+\gamma +2s} \abs{ \D^{\nu_2}_\eta\sep{ 1- \cos (r \eta . \omega)  }}{1\over{r^{1+2s}}}.
\end{split}
\end{equation}
Since  $r \in [0, \delta]$ we claim that we have the following rough estimate

\begin{lem} \label{derivdelt}  Let $0<s, \delta<1, $ and let $\omega\in\mathbb S^2$ and  $\eta\in\mathbb R^3$ be given.  Then
$
\forall ~\nu_2 \in \N^3$, $\int_0^\delta dr \abs{  \D^{\nu_2}_\eta \sep{ 1- \cos (r \omega\cdot\eta)  }} {1\over{r^{1+2s}}} \leq
C_{\delta,s} \seq{ \omega\cdot\eta}^{2s}.$
\end{lem}
\preuve[of the Lemma] This is clear for $\nu_2 = 0$ from the previous upper bound computation.

For $|\nu_2|= 1$ we have to estimate
$$
I(\nu_2) = \int_0^\delta dr \abs{ \sep{ \D^{\nu_2}_\eta \sep{ 1- \cos (r \omega\cdot\eta)  }} } {1\over{r^{1+2s}}} \leq \int_0^\delta dr \abs{  \sin (r \omega\cdot\eta)  } {1\over{r^{2s}}} .
$$
Firstly, when $0<s <1/2$, we directly get
$$
I(\nu_2) \leq \int_0^\delta dr {1 \over r^{2s}} \leq C_{s \delta} \leq C_{s \delta} \seq{ \omega\cdot\eta}^{2s} .
$$
When $s=1/2$ then
\begin{equation*}
\begin{split}
I(\nu_2) &  \leq  \int_0^{\delta |\omega\cdot\eta|} \frac{|\sin t|}{t} dt \leq \int_0^{\seq{\delta \omega\cdot\eta}} \frac{|\sin t|}{t} dt
 \leq  \int_0^{1} \frac{|\sin t|}{t} dt
+  \int_1^{\seq{\delta\omega\cdot\eta}} 1 dt \\
& \leq 1 + C_\delta \seq{\omega\cdot\eta} = C_s \seq{\omega\cdot\eta}^{2s}
\end{split}
\end{equation*}
When $ 1/2<s<1$ we have
$$
I(\nu_2) \leq |\omega\cdot\eta|^{2s -1} \int_0^\infty \frac{|\sin t|}{t^{2s}} dt \leq C_{s }|\omega\cdot\eta|^{2s -1} \leq C_{s} \seq{ \omega\cdot\eta}^{2s} .
$$
Thus we obtain  the estimate for $|\nu_2|=1.$

It remains to  consider the case when $|\nu_2| \geq 2$.   Observe $0<s<1,$ and thus
$$
I(\nu_2)  =  \int_0^\delta dr \abs{ \sep{ \D^{\nu_2}_\eta \sep{ 1- \cos (r \omega\cdot\eta)  }}} {1\over{r^{1+2s}}} \leq {\red \int_0^\delta  {{ dr}\over{r^{2s-1}}} }\leq  C_{s \delta}  \leq  C_{s \delta} \seq{ \omega\cdot\eta}^{2s}.
$$
The proof of the lemma is complete. \fin

{\bigskip\noindent {\bf End of the proof  of Proposition \ref{ap2}}}
Now we go back to (\ref{equa3}). We have also to estimate the term $\sep{ \D^{\nu_1}_v \mu (\alpha +v)}$ in this integral. For this purpose, we directly use the fact that
for all $\nu_1$,
\begin{equation} \label{unundemi}
\abs{ \D^{\nu_1}_v \mu (\alpha +v)} \leq C_{\nu_1} \mu^{1/2} (\alpha +v) .
\end{equation}
Thanks to Lemma \ref{derivdelt} and the preceding estimate, we get from (\ref{equa3}) that
\begin{equation*}
  \begin{split} 
\abs{  \D^{\nu_1}_v \D^{\nu_2}_\eta a_p(v,\eta ) }  & \lesssim   \int_\omega  d\omega
 \int_{E_{0,\omega}}  d\alpha \mu^{1/2} (\alpha +v)
 \seq{ \alpha}^{1+\gamma +2s}   \seq{ \omega\cdot\eta}^{2s} .
\end{split}
\end{equation*}
For the final estimates, we can repeat exactly the proof of the case
$\nu_1 = \nu_2=0$,   to  get the desired result.
 The proof of Proposition \ref{ap2} is complete. \fin

For further use, we shall also need the following estimate
{\red
\begin{prop} \label{moredera} The symbol $a_p$ also satisfies the following estimate: for any $0<\eps<1,$
 \begin{equation*} 
 \begin{split}
   \partial_\eta {a_p}  \in S \sep{  \eps\comii{v}^{\gamma} (1+
  |\eta|^{2s} + |\eta\wedge v|^{2s} )+\eps^{-1}\comii{v}^{\gamma+2s} , \Gamma },
\end{split}
   \end{equation*}
   with semi-norms (see Subsection \ref{subsec41} for the definition of semi-norms) independent of $\eps.$
 \end{prop}

\preuve We can again rely on the preceding arguments. We begin with (\ref{equa3}) and we can write for $|\nu_2| \geq 1$,
\begin{multline*} 
\abs{  \D^{\nu_1}_v \D^{\nu_2}_\eta a_p(v,\eta ) }   \leq C \int^\delta_0  dr \int_{S^2_\omega}d\omega \int_{E_{0,\omega}} d\alpha \sep{ \D^{\nu_1}_v \mu (\alpha +v)} | \\
 \seq{ \alpha}^{1+\gamma +2s} \sep{ \D^{\nu_2}_\eta\sep{ 1- \cos (r \eta . \omega)  }}{1\over{r^{1+2s}}}.
\end{multline*}

  \underline{\it  Suppose  that  $|\nu_2| \geq2$.} We can verify directly that, observing $0<s<1$ and $\abs\omega=1,$
 $$
\int_0^\delta dr \abs{ \sep{ \D^{\nu_2}_\eta \sep{ 1- \cos (r \omega\cdot\eta)  }}} {1\over{r^{1+2s}}}  \leq  \int_0^\delta  {{ dr}\over{r^{2s-1}}} \leq
C_{\delta,s}.  
$$
Therefore,  using also (\ref{unundemi}),
 \begin{equation*}
  \begin{split} 
\abs{ \D^{\nu_1}_v \D^{\nu_2} a_p(v,\eta ) }  & \lesssim   \int_{S^2_\omega } d\omega
 \int_{E_{0,\omega}}  d\alpha \mu^{1/2} (\alpha +v)
 \seq{ \alpha}^{1+\gamma +2s}  \lesssim
\seq{ v}^{\gamma+2s},
\end{split}
\end{equation*}
the last inequality following the same computation as that after   (\ref{typiquea}) with $|\omega\cdot\eta|^{2s}$ there replaced here by  $1.$

\underline{\it  Consider the case when  $|\nu_2| =1$.} Then we have
$$
\int_0^\delta dr \abs{ \sep{ \D^{\nu_2}_\eta \sep{ 1- \cos (r \omega\cdot\eta)  }}} {1\over{r^{1+2s}}} \leq \int_0^\delta{{ \abs{     \sin (r \omega\cdot\eta)  }}\over{r^{2s}}}dr.
$$
Furthermore  if $0<s<1/2$ then
\begin{eqnarray*}
\int_0^\delta{{ \abs{     \sin (r \omega\cdot\eta)  }}\over{r^{2s}}}dr \leq	C_{\delta,s},
\end{eqnarray*}
and if   $ 1/2<s<1$ then
\begin{eqnarray*}
\int_0^\delta{{ \abs{     \sin (r \omega\cdot\eta)  }}\over{r^{2s}}}dr& \leq	& \abs{\omega\cdot\eta}^{2s-1} \int_0^{\delta\abs{\omega\cdot\eta}}{{ \abs{     \sin \theta  }}\over{\theta^{2s}}}d\theta\leq C_{\delta,s} (1+ \seq{ \omega\cdot\eta}^{2s-1})\\
&\lesssim&  \eps \seq{ \omega\cdot\eta}^{2s}+\eps^{-(2s-1)}\lesssim \eps \seq{ \omega\cdot\eta}^{2s}+\eps^{-1}
\end{eqnarray*}
for any $0<\eps<1,$
and finally if   $ s=1/2$ then
\begin{eqnarray*}
\int_0^\delta{{ \abs{     \sin (r \omega\cdot\eta)  }}\over{r^{2s}}}dr& \leq&	   \int_0^{\delta\abs{\omega\cdot\eta}}{{ \abs{     \sin \theta  }}\over{\theta}}d\theta\leq C_{\delta,s} (1+ \ln\seq{ \omega\cdot\eta})\leq C_{\delta,s} (1+  \seq{ \omega\cdot\eta}^{s})\\
&\lesssim& \eps \seq{ \omega\cdot\eta}^{2s}+\eps^{-1}
\end{eqnarray*}
for any $\eps>0.$    Thus combining the above estimates we conclude, for $0<s<1$ and for any $0<\eps<1,$
\begin{eqnarray*}
	\int_0^\delta dr \abs{ \sep{ \D^{\nu_2}_\eta \sep{ 1- \cos (r \omega\cdot\eta)  }}} {1\over{r^{1+2s}}} \lesssim  \eps \seq{ \omega\cdot\eta}^{2s}+\eps^{-1}.
\end{eqnarray*}
Therefore we get that, using again (\ref{unundemi}) the arguments after  \eqref{typiquea},
 \begin{equation*}
  \begin{split} 
\abs{ \D^{\nu_1}_v \D^{\nu_2} a_p(v,\eta ) }  & \lesssim   \int_{S^2_\omega } d\omega
 \int_{E_{0,\omega}}  d\alpha \mu^{1/2} (\alpha +v)
 \seq{ \alpha}^{1+\gamma +2s}  \inner{ \eps \seq{ \omega\cdot\eta}^{2s}+\eps^{-1}}\\
 &\lesssim   \eps
\seq{ v}^{\gamma} \sep{1+ |\eta|^{2s}     +  |v \wedge \eta|^{2s}  }+\eps^{-1}\seq{ v}^{\gamma+2s},
\end{split}
\end{equation*}
with $\abs{\nu_2}=1.$

Combining the estimates for $\abs{\nu_2}\geq2$ and for $\abs{\nu_2}=1$ we obtain the statement in  Proposition \ref{moredera}, completing the proof. \fin
}

\subsection{Study of the multiplicative term $\bar \lll_{1,\delta,b}$}

Recall that the multiplicative part $\bar \lll_{1,\delta,b}$ has the following form
 \begin{equation*} 
 \bar \lll_{1,\delta,b} f =  -\sep{ \iint dv_* d\sigma B \tilde\phi_\delta (v'-v) \mu_*} f.
 \end{equation*}
 A nice feature of the multiplicative function defining  $\bar \lll_{1,\delta,b}$ is its good symbolic properties.
\begin{prop} \label{am2}
We can write
\begin{equation*} 
 \bar \lll_{1,\delta,b} f=  -a_m(v) f,
 \end{equation*}
where  $a_m$ is a function in $v$  satisfying  the following symbolic estimates:
 \begin{enumerate}[i)]
 \item there exists $C >0 $ such that
$ C^{-1} \seq{v}^{\gamma + 2s} \leq a_m(v, \eta) \leq  C \seq{v}^{\gamma + 2s} $;
 \item $ a_m \in S(\seq{v}^{\gamma+ 2s}, \Gamma)$.
   \end{enumerate}
   \end{prop}

\preuve
Let us again use Carleman's representation. We get
\begin{equation} \label{amfirst}
\begin{split}
a_m(v) = \int_{\R^3_h} dh \int_{E_{0,h}} d\alpha \tb \un_{|\alpha |\geq | h |} \tilde\phi_\delta(h) \mu (v+\alpha -h) | \alpha +h|^{1+\gamma +2s} {1\over{| h|^{3+2s}}}.
\end{split}
\end{equation}
In this integral $h \perp \alpha$ and $|\alpha| \geq |h|$ so that
there exists $C_s$  such that
\begin{equation} \label{ecadrealphah}
C_s^{-1} |\alpha|^{1+\gamma+ 2s} \leq  |\alpha+h|^{1+\gamma+ 2s} \leq  C_s|\alpha|^{1+\gamma+ 2s}.
\end{equation}
Therefore, shifting to spherical coordinates, and recalling that we write $\phi_\delta(h) = \phi_\delta(r)$ for $r=|h|$ by abuse of notation,  we have
\begin{eqnarray*}
a_m(v) &\lesssim &  \iint d \omega dr \int_{E_{0,\omega}} d\alpha 1_{|\alpha |\geq r}
\tilde\phi_\delta(r) \mu (v+\alpha-r\omega) |\alpha|^{1+\gamma + 2s} {1\over{ r^{1+2s}}}\\
&\lesssim &  \iint d \omega dr \int_{E_{0,\omega}} d\alpha 1_{|\alpha |\geq r}
\tilde\phi_\delta(r) \mu (v+\alpha-r\omega) |\alpha|^{1+\gamma + 2s} {1\over{ r^{1+2s}}}.
\end{eqnarray*}
Note that
$$
|v+ \alpha -r \omega|^2 = |\alpha + S(\omega)v |^2 + |(\omega\cdot v)-r|^2
$$
exactly as in (\ref{splittingorth}) so that
$$
e^{-|v+ \alpha -r \omega|^2} = e^{-| \alpha+ S(\omega)v|^2} e^{-|(\omega\cdot v)-r|^2}.
$$
Moreover, we have
$$
\int_{E_{0,\omega}} d\alpha |\alpha|^{1+\gamma + 2s} \mu (\alpha+ S(\omega) v ) \sim \seq{S(\omega)v}^{1+\gamma + 2s},
$$
and we get (forgetting the truncation function in $\alpha$)
$$
a_m(v) \lesssim  \iint d \omega dr
\tilde\phi_\delta(r) \seq{S(\omega)v }^{1+\gamma + 2s} e^{-|(\omega\cdot v)-r|^2} {1\over{ r^{1+2s}}}.
$$
We can now integrate w.r.t. $r$ and compute by virtue of
Peetre's inequality \eqref{peetre} (forgetting now the
dependence on $\delta$ for the constants)
\begin{eqnarray*}
\int dr \tilde\phi_\delta(r) e^{-|(\omega\cdot v)-r|^2} {1\over{ r^{1+2s}}}
& \lesssim &\int dr \tilde\phi_\delta(r) e^{-|(\omega\cdot v)-r|^2}
\comii{r-\omega\cdot v}^{1+2s}\comii{\omega\cdot v}^{-(1+2s)}\\
& \lesssim & \comii{\omega\cdot v}^{-(1+2s)},
\end{eqnarray*}
and thus
$$
a_m(v) \lesssim \int_{S^2_\omega} d\omega \comii{\omega\cdot v}^{-(1+2s)}\seq{S(\omega) v}^{1+\gamma + 2s}.
$$
We therefore have a similar integral as in (\ref{+typiqueplus++}) and
using exactly the same change of polar coordinates and computations as
therein with $e^{-|(\omega\cdot v)|^2}$ replaced by $\comii{\omega\cdot
  v}^{-(1+2s)}$ (see Figure \ref{polar}), we get,  just repeating the arguments between \reff{+typiqueplus++} and \reff{exactexpression},
\begin{eqnarray*}
a_m(v)& \lesssim & \int_0^\pi d\varphi \int_0^{2\pi} d\theta
\comii{\abs v \cos\varphi}^{-(1+2s)}  \sin\varphi  \sep{ 1 + |v|^2
  \sin^2\varphi}^{(1+\gamma +2s)/2}   \\
  & \lesssim & \int_0^1 dt \int_0^{2\pi} d\theta
\comii{t \abs v }^{-(1+2s)}   \sep{ 1 + |v|^2(1-t^2)
  }^{(1+\gamma +2s)/2} \\
  & \lesssim & \int_0^{1/2} dt
\comii{t \abs v }^{-(1+2s)}   \sep{ 1 + |v|^2(1-t^2)
  }^{(1+\gamma +2s)/2}\\
  && + \int_{1/2}^1 dt
\comii{t \abs v }^{-(1+2s)}   \sep{ 1 + |v|^2(1-t^2)
  }^{(1+\gamma +2s)/2}\\
 &\defegal& a_{m,1}+a_{m,2}.
\end{eqnarray*}
One has
\begin{eqnarray*}
a_{m,1}\lesssim \int_0^{1/2} dt
\comii{t \abs v }^{-(1+2s)}   \sep{ 1 + |v|^2
  }^{(1+\gamma +2s)/2} \lesssim \comii v^{\gamma+2s},
\end{eqnarray*}
and for the term $a_{m,2}$  we have, by changes of variables and using the fact that $\gamma>-3$,
\begin{eqnarray*}
a_{m,2}
&\lesssim& \int_{1/2}^1 dt
\comii{ v }^{-(1+2s)}   \sep{ 1 + |v|^2(1-t)
  }^{(1+\gamma +2s)/2}\\
  & \lesssim & \comii v^{-(1+2s)}\abs v^{-2}\int_{0}^{\abs v^2/2} d\tilde t
   \sep{ 1 + \tilde t
  }^{(1+\gamma +2s)/2}\\
  & \lesssim & \comii v^{-(1+2s)}\abs v^{-2}\int_{0}^{\abs v^2/2} d\tilde t
   \sep{ 1 + \tilde t
  }^{-(1+s)}\sep{ 1+ \tilde t
  }^{(3+\gamma +4s)/2}\\
  & \lesssim & \comii v^{-(1+2s)}\comii v^{-2} \comii v^{3+\gamma+4s}\int_{0}^{\abs v^2/2} d\tilde t
   \sep{ 1 + \tilde t
  }^{-(1+s)}\\
  & \lesssim & \comii v^{\gamma+2s}.
  \end{eqnarray*}
Combining these inequalities we conclude
\begin{eqnarray*}
a_m\lesssim \comii v^{\gamma+2s}.
\end{eqnarray*}
For the lower bound we can do essentially the same computations : because of the non-negative sign
of $a_m$ we can restrict the computations to the following subdomains in $(\alpha, h)$
$$
\set{|\alpha| \geq 10} \ \ \textrm{ and } \ \ \ \set{ |h| \leq 10},
$$
and following (\ref{amfirst}) and using (\ref{ecadrealphah}) we get
\begin{eqnarray*}
a_m(v) &\gtrsim&   \iint d \omega dr \int_{E_{0,\omega}} d\alpha \un_{|\alpha |\geq 10}
\un_{1 \leq r \leq 10} \mu (v+\alpha-r\omega) |\alpha|^{1+\gamma + 2s}
{1\over{ r^{1+2s}}}\\
&\gtrsim&  \iint d \omega dr \int_{E_{0,\omega}} d\alpha \un_{|\alpha |\geq 10}
\un_{1 \leq r \leq 10} \mu (\alpha+S(\omega)v) e^{-\abs{\omega\cdot v}^2/2} |\alpha|^{1+\gamma + 2s} {1\over{ r^{1+2s}}}
\end{eqnarray*}
since $\bar{\phi}_\delta = 1$  in the set $\set{  1 \leq r \leq 10}$
(recall $0< \delta <1 $),  and
 \[
\abs {v+\alpha-r\omega}^2=\abs {S(\omega) v+\alpha}^2+\abs
{\omega\cdot v-r}^2 \leq \abs {S(\omega) v+\alpha}^2+\abs
{\omega\cdot v}^2+100\]
for $r\leq 10.$  Then as before we can use the fact that
$$
\int d\alpha \un_{|\alpha| \geq 10}|\alpha|^{1+\gamma + 2s} \mu (\alpha+ S(\omega) v ) \sim \seq{S(\omega)v}^{1+\gamma + 2s}
$$
and
$$
\int dr \un_{1 \leq r \leq 10} {1\over{ r^{1+2s}}} \sim C
$$
and we get for a new constant $C$ that
$$
a_m(v) \geq  C^{-1} \int d \omega  \seq{S(\omega)v }^{1+\gamma + 2s}
e^{-|(\omega\cdot v)|^2/2} ,
$$
and again we can follow the computations as in (\ref{typique}) and thereafter to get
$$
a_m(v) \geq  C^{-1}  \seq{ v}^{\gamma + 2s}.
$$
The proof of i) is thus complete.

As for the proof of ii), we use \reff{amfirst} to get
\begin{eqnarray*}
\partial_v^\alpha a_m(v) = \int_{\R^3_h} dh \int_{E_{0,h}} d\alpha \tb \un_{|\alpha |\geq | h |} \tilde\phi_\delta(h) \inner{\partial_v^\alpha\mu (v+\alpha -h) }| \alpha +h|^{1+\gamma +2s} {1\over{| h|^{3+2s}}},
\end{eqnarray*}
which gives
\begin{eqnarray*}
\abs{\partial_v^\alpha a_m(v) }\lesssim \int_{\R^3_h} dh \int_{E_{0,h}} d\alpha \tb \un_{|\alpha |\geq | h |} \tilde\phi_\delta(h)  \mu \inner{\frac{v+\alpha -h}{2}} | \alpha +h|^{1+\gamma +2s} {1\over{| h|^{3+2s}}}.
\end{eqnarray*}
Then repeating the arguments as in i), we conclude that
\begin{eqnarray*}
\abs{\partial_v^\alpha a_m(v) }\lesssim \seq{ v}^{\gamma + 2s}.
\end{eqnarray*}
This completes the proof of ii).
\fin

\subsection{Proof of Proposition \ref{estaa} i)}\label{sec33}

In this subsection we prove part i) of Proposition \ref{estaa} concerning
the so-called symbol $a$. We first give its definition, then prove the Proposition, and we shall end this section by giving additional properties of $a$ which will be needed in the sequel.

\begin{defn} \label{defaapm}
We define $a$ to be the following real symbol:
$$
a = a_p + a_m,
$$
where $a_p$ is defined in Proposition \ref{ap2} and $a_m$ is defined
in Proposition \ref{am2}.
\end{defn}

We now give the proof of Proposition \ref{estaa} i).
From Proposition \ref{ap2} and Proposition \ref{am2} we know respectively that
$$
C^{-1} \seq{v}^{\gamma + 2s} \leq a_m(v, \eta) \leq  C \seq{v}^{\gamma + 2s}
$$
 and for all $0 < \kappa \leq 1$,
 $$
  C^{-1} \sep{ -\kappa \seq{v}^{\gamma+2s} + \kappa \seq{v}^{\gamma} (1+  |\eta|^{2s} + |\eta\wedge v|^{2s} ) } \leq a_p(v, \eta) \leq  C \seq{v}^{\gamma} (1+ |\eta|^{2s} + |\eta\wedge v|^{2s} ),
  $$
where in both cases $C$ denotes a constant independent of $\kappa$ (but depending on $\delta$, $s$). Choosing  $\kappa$ sufficiently small and fixed from now on,  and adding the two inequalities gives
$$
 C^{-1} \sep{  \seq{v}^{\gamma+2s} +  \seq{v}^{\gamma} (  |\eta|^{2s} + |\eta\wedge v|^{2s} ) } \leq a(v, \eta) \leq  C \seq{v}^{\gamma+2s}+ C \seq{v}^{\gamma} (1 + |\eta|^{2s} + |\eta\wedge v|^{2s} ).
  $$
  so that
  $$
 C^{-1} \seq{v}^{\gamma}\sep{ 1+ \abs{v}^{2} +   |\eta|^{2} + |\eta\wedge v|^{2}  }^s \leq a(v, \eta) \leq  C \seq{v}^{\gamma}\sep{ 1+ \abs{v}^{2} +   |\eta|^{2} + |\eta\wedge v|^{2}  }^s
  $$
for a new constant $C$. This proves the lower and upper bounds for $a$.
Using the definition of $\tilde{a}$
\begin{equation} \label{atilde}
\tilde{a}(v,\eta) = \seq{v}^\gamma \sep{1+\abs v^2+
  |\eta|^2 + |\eta\wedge v|^2 }^{s}
  \end{equation}
we get
\begin{equation} \label{pegal}
C^{-1} \tilde{a} \leq a \leq C \tilde{a}.
\end{equation}

\bigskip
From   Proposition \ref{ap2} and Proposition  \ref{am2}, we also directly get by addition that
$$
a \in S( \tilde{a}, \Gamma).
$$
Moreover, we claim that
\begin{equation}\label{caonsym}
\tilde{a} \in S( \tilde{a}, \Gamma).
\end{equation}
{\red To see this we use induction on $\abs{\alpha+\beta}$ to prove that
  for any $\kappa\in \mathbb R$ and any $\abs{\alpha+\beta}\geq 0,$
  \begin{equation}\label{eqin+1}
  \abs{\partial_v^\alpha\partial_{\eta}^\beta   \left(1+\abs v^2+
  |\eta|^2 + |\eta\wedge v|^2 \right)^{\kappa }  }
  \lesssim  \left(1+\abs v^2+
  |\eta|^2 + |\eta\wedge v|^2 \right)^{\kappa},
  \end{equation}
 which obviously holds for $\abs{\alpha+\beta}=0.$  Now suppose  $\abs{\alpha+\beta}\geq1$ then we have either $\abs\alpha\geq 1$ or $\abs\beta\geq1$,  and suppose $\abs\beta\geq 1$ without loss of generality.  So we can write    $\partial_\eta^\beta=\partial_\eta^{\tilde\beta}\partial_{\eta_j}$ with $|\tilde\beta|=\abs\beta-1$ and thus
\begin{multline*}
\partial_v^\alpha \partial_{\eta}^\beta\Big[  \left(1+\abs v^2+
  |\eta|^2 + |\eta\wedge v|^2 \right)^{\kappa }\Big] \\
  =\partial_v^\alpha\partial_{\eta}^{\tilde\beta}\Big[\kappa  \left(1+\abs v^2+
  |\eta|^2 + |\eta\wedge v|^2 \right)^{\kappa -1} \inner{2\eta_j+2\inner{\eta\wedge v}\partial_{\eta_j} \inner{\eta\wedge v}}\Big],
\end{multline*}
which along with  Leibniz'  formula and the induction assumption yields
 \begin{eqnarray*}
&&\abs{\partial_v^\alpha \partial_{\eta}^\beta\Big[  \left(1+\abs v^2+
  |\eta|^2 + |\eta\wedge v|^2 \right)^{\kappa }\Big] }\\
  &\lesssim & \left(1+\abs v^2+
  |\eta|^2 + |\eta\wedge v|^2 \right)^{\kappa -1} \inner{1+\abs\eta+ \abs{\eta\wedge v} \abs{  v} +\abs{\eta\wedge v} +\abs\eta\abs v+\abs v^2}\\
  &\lesssim & \left(1+\abs v^2+
  |\eta|^2 + |\eta\wedge v|^2 \right)^{\kappa} .
 \end{eqnarray*}
 We have proven \eqref{eqin+1}. Now using \eqref{eqin+1} and Leibniz'  formula we conclude
 \begin{eqnarray*}
 \abs{\partial_v^\alpha \partial_{\eta}^\beta\Big[ \comii v^\gamma \left(1+\abs v^2+
  |\eta|^2 + |\eta\wedge v|^2 \right)^{s}\Big] }	\leq C_{\alpha,\beta}  \comii v^\gamma \left(1+\abs v^2+
  |\eta|^2 + |\eta\wedge v|^2 \right)^{s}.
 \end{eqnarray*}
 This gives the statement in \eqref{caonsym}.}

It only remains to check the temperance of $a$ and $\tilde{a}$. From (\ref{pegal}) it is sufficient to verify that there exist two constants $N$ and $C,$  both depending only on $s$ and $\gamma$,  such that for all $Y = (y,\eta)$, $Y'=(y',\eta')$ we have
$$
\tilde{a}(Y) \leq C \tilde{a} (Y') (1 + \Gamma(Y-Y'))^N.
$$
This is a direct consequence of Peetre's inequality \eqref{peetre} since we have powers of polynomial type quantities. Indeed,  we have
\begin{eqnarray*}
\frac{\tilde a(Y)}{\tilde a(Y')} \leq  \frac{\comii y^\gamma}{\comii{y'}^\gamma} \inner{\frac{\comii y^2+\comii \eta^2+\abs{y\wedge\eta}^2}{1+\abs{y'}^2+\abs{\eta'}^2+\abs{y'\wedge\eta'}^2}}^s.
\end{eqnarray*}
On the other hand,
\begin{eqnarray*}
\frac{\comii y^\gamma}{\comii{y'}^\gamma} \leq 2^{\abs\gamma}  \comii{y-y'}^{\abs\gamma}
\end{eqnarray*}
due to  Peetre's inequality \eqref{peetre}.  Similarly,
\begin{eqnarray*}
\frac{\comii y^2+\comii \eta^2}{1+\abs{y'}^2+\abs{\eta'}^2+\abs{y'\wedge\eta'}^2}\leq 4\comii{y-y'}^2+4\comii{\eta-\eta'}^2.
\end{eqnarray*}
Moreover  using the relation
\begin{eqnarray*}
y\wedge\eta=\inner{y-y'}\wedge\inner{\eta-\eta'}+\inner{y-y'}\wedge \eta'+ y'\wedge\inner{\eta-\eta'}+y'\wedge\eta',
\end{eqnarray*}
we compute
\begin{eqnarray*}
&&\frac{\abs{y\wedge\eta}^2}{1+\abs{y'}^2+\abs{\eta'}^2+\abs{y'\wedge\eta'}^2}\\
&\leq &\frac{4\abs{y-y'} ^2\abs{\eta-\eta'}^2 +4\abs{y-y'}^2\abs{ \eta'}^2+4\abs{y'}^2\abs{\eta-\eta'}^2+4\abs{y'\wedge\eta'}^2}{1+\abs{y'}^2+\abs{\eta'}^2+\abs{y'\wedge\eta'}^2}\\
&\leq &4\abs{y-y'} ^2\abs{\eta-\eta'}^2+4\abs{y-y'}^2+4\abs{\eta-\eta'}^2+4\\
&\leq &10 \inner{\comii{y-y'}+\comii{\eta-\eta'}}^4.
  \end{eqnarray*}
Thus,
\begin{eqnarray*}
\frac{\comii y^2+\comii \eta^2+\abs{y\wedge\eta}^2}{1+\abs{y'}^2+\abs{\eta'}^2+\abs{y'\wedge\eta'}^2}
&\leq &18 \inner{\comii{y-y'}+\comii{\eta-\eta'}}^4.
\end{eqnarray*}
Combining  the above inequalities, we get
\begin{eqnarray*}
\frac{\tilde a(Y)}{\tilde a(Y')} \leq C_{s,\gamma} \inner{\comii{y-y'}+\comii{\eta-\eta'}}^{4s+\abs\gamma}\leq \tilde C_{s,\gamma} \inner{1+\Gamma(Y-Y')}^{4s+\abs\gamma}
\end{eqnarray*}
with $C_{s,\gamma}$ and $\tilde C_{s,\gamma}$ two constants depending only on $s$ and $\gamma.$  The temperance  of $\tilde a$ follows.
The proof is complete. \fin

\bigskip
For further use we also give here two propositions concerning $a$ and $\tilde{a}
$, which will be of great interest in the next section.

\begin{prop}\label{propsym} Recall $\tilde a(v,\eta)=\comii{v}^{\gamma} \left(1+\abs v^2+
  |\eta|^2 + |\eta\wedge v|^2 \right)^{s}.$ We have \\
\begin{enumerate}[i)]
\item for any $\abs\alpha\geq 0$ and any $\abs\beta\geq 1$ ,  there exist two
constants $C_{\alpha,\beta}>0$ and $C_\beta$ such that
\begin{equation*}   \abs{\partial_v^\alpha \partial_\eta^\beta a} \leq {\red C_{\alpha,\beta} \inner{\eps \tilde a
  +\eps^{-1} \comii v^{2s+\gamma}}}
 \end{equation*}
and
\[
\abs{\partial_\eta^\beta \tilde{a} }  \leq  C_{\beta}
\comii{v}^{\gamma+1} \left(1+\abs v^2+
  |\eta|^2 + |\eta\wedge v|^2 \right)^{s-1/2};
\]
\item the following estimate is true for any $0 <\eps \leq 1$, with semi-norms (see Subsection \ref{subsec41} for the definition of semi-norms) independent of $\eps$:
\begin{equation}\label{etader}
  \partial_\eta \tilde a, ~ \partial_\eta a \in  { \red  S( \eps  a+\eps^{-1} \comii v^{2s+\gamma}, ~\Gamma)};
\end{equation}

\item we have
\begin{equation} \label{derx}
\abs{\xi\cdot\partial_\eta \tilde a} \lesssim \comii{v}^{\gamma} \left(1+\abs v^2+
  |\eta|^2 + |\eta\wedge v|^2 \right)^{s-\frac{1}{2}} \sep{ \abs{
   \xi}^2+\abs{v\wedge \xi }^2}^{1/2}.
   \end{equation}
\end{enumerate}
\end{prop}

\preuve
The  point i)  for $a$ is just an immediate consequence of   Proposition
\ref{moredera}.  Now  we check for $\tilde a$.  Recall
$$\tilde a(v,\eta)=\comii{v}^{\gamma} \left(1+\abs v^2+
  |\eta|^2 + |\eta\wedge v|^2 \right)^{s}.$$
  We claim, for any $\kappa\in \mathbb R$ and any $\abs\beta\geq 1,$
  \begin{multline*}
  \abs{\partial_{\eta}^\beta\Big[\comii{v}^{\gamma} \left(1+\abs v^2+
  |\eta|^2 + |\eta\wedge v|^2 \right)^{\kappa } \Big]}
  \lesssim \comii{v}^{\gamma+1} \left(1+\abs v^2+
  |\eta|^2 + |\eta\wedge v|^2 \right)^{\kappa-\frac{1}{2}},
  \end{multline*}
which can be deduced by induction on $\abs\beta. $   Indeed, by direct computation we see the above estimate holds for $\abs{\beta}=1$, since
\begin{eqnarray*}
&&\abs{\partial_{\eta_j} \Big[\comii{v}^{\gamma} \left(1+\abs v^2+
  |\eta|^2 + |\eta\wedge v|^2 \right)^{\kappa } \Big]}\\
  &= &\abs{\kappa\comii{v}^{\gamma} \left(1+\abs v^2+
  |\eta|^2 + |\eta\wedge v|^2 \right)^{\kappa -1} \inner{2\eta_j+2\inner{\eta\wedge v}\partial_{\eta_j} \inner{\eta\wedge v}}}\\
  &\lesssim& \comii{v}^{\gamma} \left(1+\abs v^2+
  |\eta|^2 + |\eta\wedge v|^2 \right)^{\kappa -1} \inner{\abs{\eta}+ \abs{\eta\wedge v} \abs{v}}\\
  &\lesssim & \comii{v}^{\gamma+1} \left(1+\abs v^2+
  |\eta|^2 + |\eta\wedge v|^2 \right)^{\kappa -{1\over 2}}.
  \end{eqnarray*}
  Moreover for any $\abs\beta\geq 2$, we may write $\partial_\eta^\beta=\partial_\eta^{\tilde\beta}\partial_{\eta_j}$ with $|\tilde\beta|=\abs\beta-1$ and thus
\begin{multline*}
\partial_{\eta}^\beta\Big[\comii{v}^{\gamma} \left(1+\abs v^2+
  |\eta|^2 + |\eta\wedge v|^2 \right)^{\kappa }\Big] \\
  =\partial_{\eta}^{\tilde\beta}\Big[\kappa\comii{v}^{\gamma} \left(1+\abs v^2+
  |\eta|^2 + |\eta\wedge v|^2 \right)^{\kappa -1} \inner{2\eta_j+2\inner{\eta\wedge v}\partial_{\eta_j} \inner{\eta\wedge v}}\Big].
\end{multline*}
As a result,  by Leibniz's formula and the induction assumption  on $\abs\beta$,  we obtain
\begin{eqnarray*}
&&\abs{\partial_{\eta}^\beta\Big[\comii{v}^{\gamma} \left(1+\abs v^2+
  |\eta|^2 + |\eta\wedge v|^2 \right)^{\kappa }\Big] }\\
  &\lesssim&   \bigg[\comii{v}^{\gamma+1} \left(1+\abs v^2+
  |\eta|^2 + |\eta\wedge v|^2 \right)^{\kappa-1-{1\over 2}} \bigg]\inner{1+\abs{\eta}+ |\eta\wedge v|\cdot \abs v+\abs v^2}\\
  &\lesssim&    \comii{v}^{\gamma+1} \left(1+\abs v^2+
  |\eta|^2 + |\eta\wedge v|^2 \right)^{\kappa-1-{1\over 2}}  \inner{1+\abs{\eta}^2+ |\eta\wedge v|^2+\abs v^2}\\
  &\lesssim  &  \comii{v}^{\gamma+1} \left(1+\abs v^2+
  |\eta|^2 + |\eta\wedge v|^2 \right)^{\kappa-1/2}.
\end{eqnarray*}
Applying the above inequalities  for $\kappa=s$, we obtain the desired  estimate for $\tilde a.$

Next we prove Point ii).  The conclusion for $\partial_\eta a$ follows from the    estimates in  i).  And we have to check $\partial_\eta\tilde a, $ and we have shown in i) that
\begin{eqnarray*}
\abs{\partial_{\eta} \tilde a}
 & \lesssim &   \comii{v}^{\gamma+1} \left(1+\abs v^2+
  |\eta|^2 + |\eta\wedge v|^2 \right)^{s-1/2}\\
   & \lesssim &   \comii{v}^{\gamma/2+s} \tilde a^{1/2}.
\end{eqnarray*}
 Then arguing as above we can use induction on $\abs\alpha+\abs\beta$ to obtain, for $\abs\alpha+\abs\beta\geq 0, $
\begin{eqnarray*}
\abs{\partial_v^\alpha\partial_\eta^\beta\partial_{\eta} \tilde a}
  \lesssim     \comii{v}^{\gamma/2+s} \tilde a^{1/2}.
\end{eqnarray*}
This gives the conclusion for $\partial_\eta\tilde a$.

 Point iii) in Proposition \ref{propsym} is a direct consequence of the computation on $\tilde{a}$,
 since
 \begin{eqnarray*}
 \xi\cdot\partial_\eta \tilde a= s \comii{v}^{\gamma} \left(1+\abs v^2+
  |\eta|^2 + |\eta\wedge v|^2 \right)^{s-1} \sep{
   2\xi\cdot \eta  + 2(v\wedge \xi)\cdot (v\wedge \eta) }.
 \end{eqnarray*}
The proof is complete.
\fin

\subsection{Study of the subprincipal term $\lll_{1,1, \delta}$}

\begin{prop} \label{as2} We can write
 \begin{equation*}
 \lll_{1,1, \delta} f =  - a_s(v, D_v) f,
 \end{equation*}
 where $a_s$, defined by \eqref{asdef} below,  is a (complex valued)  classical symbol in  $(v,\eta)$ satisfying that
 for all $0 < s < 1$ and any $0 <\eps<1$,  we have, with
    semi-norms independent of $\eps$,
 \begin{equation} \label{asii}
 a_s (v, \eta) \in {\red   S\sep{ \eps {a} + \eps^{-1} \seq{v}^{\gamma+2s} , \Gamma}}.
 \end{equation}
 \end{prop}

\preuve
We recall that
$$\lll_{1,1,\delta} f = \iint dv_* d\sigma B \phi_\delta (|v'-v|)  (\mu'_*)^{1/2}  [f' - f]\sep{ (\mu_*)^{1/2} - (\mu'_*)^{1/2}}. $$
We shift to Carleman's representation and get
\begin{equation*}
\begin{split}
\lll_{1,1,\delta} f & =  \int_{\R^3_h} dh \int_{E_{0,h}} d\alpha \tb 1_{|\alpha | \geq |h|} |\alpha+h|^{1+\gamma+2s} \phi_\delta (|h|)  \mu^{1\over 2} (\alpha +v)   [f(v-h) - f (v)] \\
& \ \ \ \ \ \ \ \ \ \ \ \ \
 \sep{ \mu^{1\over 2} (\alpha + v-h) - \mu^{1\over 2} (\alpha +v)} {1\over{| h |^{3+2s}}} \\
& =- \int_{\R^3_\eta} \hat f(\eta )e^{ iv.\eta} a_s(v,\eta )d\eta
\end{split}
\end{equation*}
with
 \begin{multline}\label{asdef}
a_s(v,\eta )   = -\int_{\R^3_h} dh \int_{E_{0,h}} d\alpha \tb 1_{|\alpha | \geq |h|} |\alpha+h|^{1+\gamma+2s} \phi_\delta (|h|)  \mu^{1\over 2} (\alpha +v)   [e^{-ih\cdot\eta} -1]\\
  \sep{ \mu^{1\over 2} (\alpha + v-h) - \mu^{1\over 2} (\alpha +v)} {1\over{| h |^{3+2s}}} .
\end{multline}

For the study of this symbol, we shall essentially follow the same computations as in the $\lll_{1,2,\delta}$ case. We first note that
we have the following  bound for all $h \neq 0$
$$
\abs{ \sep{ \mu^{1\over 2} (\alpha + v-h) - \mu^{1\over 2} (\alpha +v)} {1\over{| h |}}} \leq C .
$$
So that using also that $\abs\alpha \leq |\alpha + h| \leq 2 |\alpha|$  due to the fact that $\alpha \perp h, $ we get
\begin{equation*}
\begin{split}
|a_s(v,\eta )| & \lesssim  \int_{\R^3_h} dh \int_{E_{0,h}} d\alpha\abs{\alpha}^{1+\gamma+2s}   \mu^{1\over 2} (\alpha +v)    \phi_\delta (|h|) {|e^{-ih\cdot\eta} -1| \over{| h |^{2+2s}}}.
\end{split}
\end{equation*}
Now we shift to spherical coordinates taking $h=r \omega$ and we get
\begin{equation}\label{esasveta}
 |a_s(v,\eta )|   \lesssim  \int_0^{+\infty}\int_{ S_\omega^2} d\omega dr \int_{E_{0,\omega}} d\alpha \abs{\alpha}^{1+\gamma+2s}   \mu^{1\over 2} (\alpha +v)    \phi_\delta (r) {|e^{-ir\omega\cdot\eta} -1| \over{r^{2s}}} .
\end{equation}
{\red
We can directly integrate w.r.t. $r$  and this gives
$$
\int_0^\infty \phi_\delta (r) {|e^{-ir\omega\cdot\eta} -1| \over{r^{2s}}}dr \lesssim   \int_0^\delta  {|\cos\inner{r\omega\cdot\eta} -1| \over{r^{2s}}}dr+\int_0^\delta  {|\sin\inner{r\omega\cdot\eta} | \over{r^{2s}}}dr.
  $$
We have proven in the proof of Proposition \ref{moredera} (see the treatment of the case $\abs{\nu_2}=1$ threein)  that
\begin{eqnarray*}
	\int_0^\delta  {|\sin\inner{r\omega\cdot\eta} | \over{r^{2s}}}dr\lesssim \eps \abs{ \omega\cdot\eta}^{2s}+\eps^{-1}
\end{eqnarray*}
for any $0<\eps<1.$ Furthermore if $0<s<1/2$ then
\begin{eqnarray*}
\int_0^\delta  {|\cos\inner{r\omega\cdot\eta} -1| \over{r^{2s}}}dr \lesssim \int_0^\delta  {1 \over{r^{2s}}}dr\leq C_{\delta,s},
\end{eqnarray*}
and if  $1/2<s<1$ then
\begin{eqnarray*}
&&\int_0^\delta  {|\cos\inner{r\omega\cdot\eta} -1| \over{r^{2s}}}dr \lesssim
\abs{\omega\cdot\eta}^{2s-1}\int_0^{\delta\abs{\omega\cdot\eta}}  {|\cos\theta -1| \over{\theta^{2s}}}d\theta \\
&\lesssim&
\abs{\omega\cdot\eta}^{2s-1}\int_0^{\min\set{1,\delta\abs{\omega\cdot\eta}}}  {|\cos\theta -1| \over{\theta^{2s}}}d\theta+\abs{\omega\cdot\eta}^{2s-1}\int_{\min\set{1,\delta\abs{\omega\cdot\eta}}}^{\delta\abs{\omega\cdot\eta}}  {|\cos\theta -1| \over{\theta^{2s}}}d\theta \\
&\lesssim& \abs{\omega\cdot\eta}^{2s-1}\int_0^{1}  {1 \over{\theta^{2s-1}}}d\theta+ \abs{\omega\cdot\eta}^{2s-1}\int_1^{+\infty}  {1 \over{\theta^{2s}}}d\theta\\
  &\lesssim  &\abs{\omega\cdot\eta}^{2s-1} \lesssim \eps\abs{ \omega\cdot\eta}^{2s}+\eps^{-(2s-1)} \lesssim \eps\abs{ \omega\cdot\eta}^{2s}+\eps^{-1},
\end{eqnarray*}
and finally  if  $s=1/2$ then
\begin{eqnarray*}
\int_0^\delta  {|\cos\inner{r\omega\cdot\eta} -1| \over{r^{2s}}}dr
&\leq &
 \int_0^{\min\set{\eps,\delta}} {|\cos\inner{r\omega\cdot\eta} -1| \over{r}}dr +\int_{\min\set{\eps,\delta}} ^\delta {|\cos\inner{r\omega\cdot\eta} -1| \over{r}}dr \\
 &\lesssim  &\abs{\omega\cdot\eta}\int_0^{\min\set{\eps,\delta}}   dr+\eps^{-1}\lesssim \eps\abs{ \omega\cdot\eta} +\eps^{-1}=\eps\abs{ \omega\cdot\eta}^{2s} +\eps^{-1}.
\end{eqnarray*}
Combining the above estimate we have
\begin{eqnarray*}
	\int_0^\infty \phi_\delta (r) {|e^{-ir\omega\cdot\eta} -1| \over{r^{2s}}}dr\lesssim \eps\abs{ \omega\cdot\eta}^{2s} +\eps^{-1},
\end{eqnarray*}
and thus, in view of \eqref{esasveta},
\begin{multline*} \label{doubble}
|a_s(v,\eta )|   \lesssim\eps
\int_{S^2_\omega} d\omega  \int_{E_{0,\omega}} d\alpha\abs{\alpha}^{1+\gamma+2s}   \mu^{1\over 2} (\alpha +v)   |\omega\cdot\eta|^{2s} \\+   \eps^{-1}  \int_{S^2_\omega} d\omega  \int_{E_{0,\omega}} d\alpha\abs{\alpha}^{1+\gamma+2s}   \mu^{1\over 2} (\alpha +v).
\end{multline*}
This enables us to do exactly the same computations as in the $\lll_{1,2,\delta}$ case,
with  the factors $\mu(\alpha+ v)$  in formula
(\ref{typiquea}) replaced respectively by $\mu^{1/2}(\alpha+ v)$  here and the factor $\abs{\omega\cdot\eta}^{2s}$ by $1$.
 We directly get, following the computations after (\ref{typiquea}) , that
\begin{eqnarray*}
|a_s(v,\eta )|   \lesssim\eps
	 \seq{v}^{\gamma}(1+ |\eta|^{2s} + |v \wedge \eta|^{2s})+\eps^{-1} \seq{v}^{\gamma+2s}\lesssim \eps a+\eps^{-1} \seq{v}^{\gamma+2s},
\end{eqnarray*}
the last inequality using \eqref{pegal}.

Again the proof of the estimates for higher order derivatives of $a_s$
is similar to the one of order $0$,   and we skip this part of the
proof for brevity.  This completes the
proof of Proposition \ref{as2}.}
\fin

\section{Proof of the main results}

This section is devoted to the proof of the main results mentionned in the introduction, including in particular Theorems \ref{thmain} and \ref{th1}. We shall use extensively properties of the classical Weyl  and Wick quantizations,   and a brief review of these properties is given in the Appendix.
In Subsection \ref{subsec42} we make the reduction to  the hypoelliptic
problems for  a simplified operator,  by virtue of
Proposition \ref{estaa} whose proof is also presented in this
subsection.   In Subsection \ref{subsec43},  we  give some coercivity estimates, and recover a result of coercivity of \cite{AMUXY1} implying the so-called triple norm.  The proof of the main results is then achieved in the last subsection \ref{subsec44}.

\subsection{Proof of Proposition \ref{estaa} ii) and iii) and related results } \label{subsec42}
In the previous sections, we splitted operator $\lll$ into several pieces  in the
following way,  with $a=a_p+a_m$  defined in Proposition \ref{ap2} and
Proposition \ref{am2},  and $a_s$ defined in Proposition \ref{as2},
\begin{eqnarray*}
    \lll &=&\mathcal L_1+\mathcal L_2= -a(v,D_v)+\mathcal L_2 +\bar {\mathcal L}_{1,\delta,a}+\mathcal
    L_{1,3,\delta}+\mathcal
    L_{1,4,\delta}-a_s(v,D_v)\\
&=&-a^w-\underbrace{\inner{-\mathcal L_2-\bar {\mathcal
      L}_{1,\delta,a}-\mathcal
    L_{1,3,\delta}-{\mathcal L}_{1,4,\delta}+a_s(v,D_v)+(a(v,D_v)-a^w)}}_{\mathcal K},
\end{eqnarray*}
recalling that $\mathcal L_1,\mathcal L_2$ are defined by \eqref{L1L2def},  $a(v,D_v)=-{\mathcal L}_{1,2,\delta}-\bar {\mathcal L}_{1,\delta,b}$ and $a_s(v,D_v)=- {\mathcal L}_{1,1,\delta},$ and
 $\bar {\mathcal L}_{1,\delta,a}, \bar {\mathcal L}_{1,\delta,b} $ and ${\mathcal L}_{1,j,\delta}, 1\leq j\leq 4,$ are given by \eqref{L1detaa}-\eqref{splitL1+}.
Thus we can write
$$
P = v\cdot \partial_x + a^w + \kkk.
$$
Notice that the diffusion term $a^w+\mathcal K$ above is only  an operator
with respect to the velocity variable $v$.  So we only work on the
resulting
operator after performing
partial
Fourier transform in the $x$ variables,
considering the dual variables $\xi$ of $x$ as
parameter.
More precisely we will study the operator
\begin{eqnarray*} 
  \hat{ P}_K=i\inner{v\cdot\xi}+a_K^w,
\end{eqnarray*}
where $a_K$ is given by \eqref{deak}, i.e.,
$$
a_K = a + K  \seq{v}^{2s + \gamma}.
$$
with $K$ a fixed number, constructed in  Lemma \ref{inverse} and Lemma
\ref{coK} below,
depending only on the integer $N$ in \reff{bdness}.   Accordingly we also introduce  the
weight function
$$
\tilde{a}_K = \tilde{a} + K  \seq{v}^{2s + \gamma},
$$
where $\tilde{a}$ is the weight function given in Proposition \ref{estaa}.   We claim that $\tilde a_K$ is temperate  uniformly with respect to $K$.  Indeed,  by  Proposition \ref{estaa} i),  whose proof is given in Subsection \ref{sec33},  we see $\tilde a$ is temperate weight with respect to $\Gamma$, i.e., there exist two constants $C$ and $N$,  both depending only on $\gamma, s$,  such that
    \begin{eqnarray*}
       \forall~(v,\eta), (w, \zeta)\in\mathbb R^6,\quad \frac{\tilde a(v,\eta)}{\tilde a(w, \zeta)}\leq C\inner{\comii{v-w}+\comii{\eta-\zeta}}^N.
    \end{eqnarray*}
Thus for any $ (v,\eta), (w, \zeta)\in\mathbb R^6$,
\begin{eqnarray*}
 \frac{\tilde a_K(v,\eta)}{\tilde a_K(w, \zeta)}&=&\frac{\tilde a(v,\eta)}{\tilde a(w, \zeta)+K\comii{w}^{2s+\gamma}}+\frac{K\comii{v}^{2s+\gamma}}{\tilde a(w, \zeta)+K\comii{w}^{2s+\gamma}}\\
 &\leq &\frac{\tilde a(v,\eta)}{\tilde a(w, \zeta)}+\frac{K\comii{v}^{2s+\gamma}}{K\comii{w}^{2s+\gamma}}\\
 &\leq& C\inner{\comii{v-w}+\comii{\eta-\zeta}}^N+2^{\abs{2s+\gamma}}\comii{v-w}^{\abs{2s+\gamma}}\\
 &\leq& \inner{C+2^{\abs{2s+\gamma}}}\inner{\comii{v-w}+\comii{\eta-\zeta}}^{N+\abs{2s+\gamma}},
\end{eqnarray*}
the second inequality using peetre's inequality \eqref{peetre}. This gives  $\tilde a_K$ is temperate  uniformly with respect to $K$, since the constant $C$ above is independent of $K.$

We note that $ a_K \in S(\tilde{a}_K, \Gamma)$ uniformly in $K$, since for any multi-index $\alpha, \beta\in\mathbb Z_+^3$,  we have
\begin{eqnarray*}
\abs{\partial_v^\alpha\partial_\eta^\beta  a_K(v,\eta)}&\leq& \abs{\partial_v^\alpha\partial_\eta^\beta  a(v,\eta)}+\abs{\partial_v^\alpha\partial_\eta^\beta \inner{ K\comii v^{2s+\gamma}}}\\
&\leq&  C_{\alpha,\beta}  \,a(v,\eta)+K C_{\alpha,\beta} \comii v^{2s+\gamma}\\
&\leq&  2 C_{\alpha,\beta}  \,a_K(v,\eta)\leq   C_{\alpha,\beta, \gamma,s}  \,\tilde a_K(v,\eta) ,
\end{eqnarray*}
with $C_{\alpha,\beta}$ a constant depending only on $\alpha,\beta,$  and $C_{\alpha,\beta,\gamma, s}$ a constant depending only on $\alpha,\beta, \gamma$ and $s.$  Thus $ a_K \in S(\tilde{a}_K, \Gamma)$ uniformly in $K$.   More generally  we can show, for $r\in~[-1,1]$,
\begin{eqnarray*}
\forall~\alpha\in \mathbb Z_+^6, \quad\abs{\partial^\alpha a_K^r} \leq C_\alpha a_K^r\leq \tilde C_\alpha \tilde a_K^r
\end{eqnarray*}
by induction on $\abs\alpha,$
which gives
$
a_K^r \in
 S(\tilde{a}_K^r, \Gamma)$ uniformly w.r.t. $K$  for all $r \in [-1,1]$.  Working with $a_K^w$  instead of $a^w$ will enable us to construct the inverse of
the former, see Lemma \ref{inverse} below.  This  is of big importance in the following analysis of hypo-elliptic estimates.

{\bf Notations.} In the following, let $K$ be fixed, satisfying the
assumptions in  Lemma \ref{inverse} and Lemma
\ref{coK} below, and let $\ell\in \mathbb R$ be an arbitrary number,
fixed and as small as we want.    To simplify the notation, by $A\lesssim B$ we mean
there exists a positive   constant  $C$,  which may depend on $K$  and
$\ell$  but  is  {\it independent of the
  parameters } $\xi$,    such that $A\leq C B$,  and similarly for $A\gtrsim B$. While the notation $A\approx B$ means both $A\lesssim B$ and $B\lesssim A$ hold.
 Given a symbol $q$  and a weight function $M$,   by $q\in S(M,\Gamma)$
we alway mean, in the following discussion, $q$ lies in $S(M,\Gamma)$
\textit{uniformly w.r.t. $K$ and $\xi$. }

Now we state the main result of this subsection, which shows that it is
sufficient to study the operator  $\hat P_K$ instead of  the original one.
\begin{prop}\label{prpmain}
  The conclusion in Theorem \ref{thmain} holds true if the
   estimate
    \begin{equation} \label{resulting}
   \norm{  \tilde a (v,\xi)^{\frac{1}{1+2s}} f }+ \norm{ a_K^w f } \lesssim \norm{
  \hat P_K  f}_{L^2} +  \norm{f}_{L_\ell^2}
    \end{equation}
holds uniformly with respect to $\xi$. Recall $a_K$ is given by \eqref{deak}.
\end{prop}

We proceed to prove the above proposition through several lemmas.
Firstly  we begin with  the construction of the inverses of  operators.

 \begin{lem}\label{inverse}
There exists a $K_0$ sufficiently large, depending only on a fixed finite number
of semi-norms of $a$, 
 such that for all $K\geq
K_0$ we have
 \begin{enumerate}[(i)]
 \item
  $a_K^w$, with $a_K$ defined by \eqref{deak},  is invertible and  its inverse $\inner{a_K^w}^{-1}$  has the
  form
\[
    \inner{a_K^w}^{-1} =H_1 \inner{a_K^{-1}}^w=\inner{a_K^{-1}}^wH_2,
\]
with  $H_1, H_2$ belonging  to  $ \mathcal B (L^2)$,   the space of bounded operators
on $L^2$, and   $\norm{H_j}_{\mathcal B(L^2)}$ bounded from above by some constant independent of $K$ for $j=1,2$;
     \item   $\inner{a_K^{1/2}}^w$ is invertible and   its
       inverse $\com{\inner{a_K^{1/2}}^w}^{-1} $ has the form
\[
   \com{\inner{a_K^{1/2}}^w}^{-1}= G_1
   \inner{a_K^{-1/2}}^w=\inner{a_K^{-1/2}}^wG_2
\]
with  $G_1 ,G_2 \in\mathcal B(L^2)$ and $\norm{G_j}_{\mathcal B(L^2)}$  bounded from above by some constant independent of $K$ for $j=1,2$;
     \item   $\inner{\tilde a_K^{1/2} a_K^{1/2}}^w$ is invertible and
       its
       inverse $\com{\inner{\tilde a_K^{1/2} a_K^{1/2}}^w}^{-1} $ has the form
\[
   \com{\inner{\tilde a_K^{1/2}  a_K^{1/2}}^w}^{-1}= Q_1
   \inner{\tilde a_K^{-1/2} a_K^{-1/2}}^w=\inner{\tilde a_K^{-1/2} a_K^{-1/2}}^wQ_2
\]
with  $Q_1 ,Q_2 \in\mathcal B(L^2)$ and $\norm{Q_j}_{\mathcal B(L^2)}$ bounded from above by some constant independent of $K$ for $j=1,2$.
\end{enumerate}
\end{lem}

\preuve
 Note first that in all what follows, we shall crucially use the fact that only a
finite number $N$ (depending only on the dimension $n=3$ here) of seminorms of a symbol is
needed to control the norm of the corresponding pseudo-differential operator (see  \reff{bdness} here and e.g. \cite[Lemma 2.5.2]{MR2599384}). \color{black}

Let us now prove the conclusion in (i).  Using \reff{11053010} and \reff{expan}, we may write
\begin{equation} \label{+inverse}
   a_K^w (a_K^{-1})^w=  {\rm Id}- R_K^w,
\end{equation}
where
\begin{eqnarray*}
   R_K=-\int_0^1\inner{ \partial_\eta a_K}\sharp_\theta
   \inner{\partial_v( a_K^{-1})}  d\theta+\int_0^1
   \inner{\partial_v a_K}  \sharp_\theta \inner{\partial_\eta(
     a_K^{-1})}  d\theta
\end{eqnarray*}
with $g \sharp_\theta h$ defined by
\begin{equation} \label{sharptheta}
    g \sharp_\theta h  (Y) = \iint e^{-2i \sigma(Y-Y_1, Y-Y_2)/\theta}
   \frac{1}{2i}g(Y_1)
   h(Y_2) dY_1 dY_2/(\pi\theta)^{6}.
\end{equation}
Let now $N$ be the integer which is given in \reff{bdness}
(and therefore depending  only on the dimension $n=3$ here).  By
\cite[Proposition 1.1]{MR1721321}  we can find a constant $C_N$ and a
positive
integer $\ell_N$,
both depending only on $N$ but independent of $K$ and $\theta$,  such that
\begin{eqnarray*}
   \norm{\inner{ \partial_\eta a_K}\sharp_\theta
   \inner{\partial_v( a_K^{-1})} }_{N; S(1,\Gamma)} \leq C_N  \norm{ \partial_\eta a_K}_{\ell_N; S(\tilde{a}_K,\Gamma)} \norm{\inner{\partial_v( a_K^{-1})}}_{\ell_N; S(\tilde{a}_K^{-1},\Gamma)},
\end{eqnarray*}
where the semi-norm $\norm{\cdot}_{k;S(M,\Gamma)}$ is defined by \reff{seminorm}.
Moreover, using  \reff{etader} for $\eps=K^{-1/2}$  yields
\begin{eqnarray*}
   \norm{ \partial_\eta a_K}_{\ell_N; S(\tilde{a}_K,\Gamma)} \leq \tilde C_N K^{-\frac{1}{2}}
\end{eqnarray*}
and from the fact  $a_K\in S(\tilde a_K, \Gamma) $ it follows that $a_K^{-1}\in S(\tilde a_K^{-1}, \Gamma) $,  and thus
\begin{eqnarray*}
   \norm{ \partial_v (a_K^{-1})}_{\ell_N; S(\tilde{a}_K^{-1},\Gamma)} \leq \tilde C_N
\end{eqnarray*}
with $\tilde C_N$ a constant depending only on $N$ but independent of
$K$.
As a result,
\begin{eqnarray*}
      \norm{\inner{ \partial_\eta a_K}\sharp_\theta
   \inner{\partial_v( a_K^{-1})} }_{N; S(1,\Gamma)} \leq C_N \tilde
 C_N^2 K^{-\frac{1}{2}}.
\end{eqnarray*}
Similarly,
\begin{eqnarray*}
      \norm{\inner{ \partial_v a_K}\sharp_\theta
   \inner{\partial_\eta( a_K^{-1})} }_{N; S(1,\Gamma)} \leq C_N \tilde
 C_N^2 K^{-\frac{1}{2}} .
\end{eqnarray*}
Then
\begin{eqnarray*}
   \norm{R_K}_{N; S(1,\Gamma)} \leq 2 C_N  \tilde C_N^2 K^{-\frac{1}{2}},
\end{eqnarray*}
and thus by \reff{bdness}
\[
  \norm{R_K^w}_{\mathcal B(L^2)} \leq 2C C_N  \tilde C_N^2 K^{-\frac{1}{2}}
\]
with $C$ a constant depending only on the dimension.
This  implies  ${\rm Id}- R_K^w$ is invertible in the space  $\mathcal B (L^2)$
of bounded operators on $L^2$ if we choose  $K$ in such a way that
$K\geq \inner{ 4CC_N\tilde C_N^2}^2$.   Moreover
\[
      \inner{{\rm Id}- R_K^w }^{-1}=\sum_{j=0}^\infty \inner{R_K^w}^j\in \mathcal B (L^2) .
\]
Taking into account \reff{+inverse}, we conclude
\[
    a_K^w \inner{(a_K^{-1})^w   \inner{{\rm Id}- R_K^w }^{-1} } ={\rm Id}.
\]
Similarly we can find a $\tilde R_K
\in S(1, \Gamma)$ such that
\[
  \inner{ \inner{{\rm Id}- \tilde R_K^w }^{-1}  (a_K^{-1})^w  }  a_K^w ={\rm Id}.
\]
These facts imply  $a_K^w$ is invertible and its inverse
$\inner{a_K^w}^{-1}$ has the form
\[
   \inner{a_K^w}^{-1} =(a_K^{-1})^w   \inner{{\rm Id}- R_K^w }^{-1} = \inner{{\rm Id}- \tilde R_K^w }^{-1}  (a_K^{-1})^w.
\]
We have proved the conclusion in (i)  in Lemma
\ref{inverse}.
 The remaining  proofs  in  (ii) and (iii)  can be deduced quite
similarly and are therefore omitted. The proof of Lemma \ref{inverse} is thus complete.
\fin

In the following, we always let $K$ be fixed satisfying the
condition in  the above  lemma \ref{inverse}.

\begin{cor}\label{teccor}
     Let $\eps$ be an arbitrarily small number and let   $g \in S\inner {\eps a_K +\eps^{-1} \comii
       v^{2s+\gamma},~\Gamma}$ uniformly with respect to $\eps$.   Then
\[
   \norm{g (v,D_v) f}_{L^2} +\norm{g^w f}_{L^2} \lesssim \eps \norm{a_K^w f}+\eps^{-1} \norm{\comii v^{2s+\gamma}f}_{L^2}.
\]
Recall $a_K$ is given in \eqref{deak}.
\end{cor}

\preuve
We first show that $  \eps a_K +\eps^{-1} \comii
       v^{2s+\gamma} $ is a temperate weight uniformly with respect to $\eps. $ Recall $a_K(v,\eta)=a(v,\eta)+K\comii v^\gamma.$    By  Proposition \ref{estaa} i),  whose proof is given in Subsection \ref{sec33},  we see $a$ is temperate weight with respect to $\Gamma$, i.e., there exist two  constants $N$ and $C$,  both depending only on $\gamma, s$,  such that
    \begin{eqnarray*}
       \forall~(v,\eta), (\tilde v,\tilde\eta)\in\mathbb R^6,\quad \frac{a(v,\eta)}{a(\tilde v,\tilde\eta)}\leq C\inner{\comii{v-\tilde v}+\comii{\eta-\tilde\eta}}^N.
    \end{eqnarray*}
As a result,
    \begin{eqnarray*}
       \forall~(v,\eta), (\tilde v,\tilde\eta)\in\mathbb R^6,\quad \frac{\eps a(v,\eta)}{ \eps a_K(\tilde v,\tilde\eta) +\eps^{-1} \comii
       {\tilde v}^{2s+\gamma}}\leq \frac{\eps a(v,\eta)}{ \eps a(\tilde v,\tilde\eta) } \leq C\inner{\comii{v-\tilde v}+\comii{\eta-\tilde\eta}}^N.
    \end{eqnarray*}
    Moreover,  for any $(v,\eta), (\tilde v,\tilde\eta)\in\mathbb R^6,$
    \begin{eqnarray*}
    \frac{\eps K \comii v^{2s+\gamma}}{ \eps a_K(\tilde v,\tilde\eta) +\eps^{-1} \comii
       {\tilde v}^{2s+\gamma}}\leq \frac{\eps K \comii v^{2s+\gamma}}{ \eps K\comii
       {\tilde v}^{2s+\gamma} }\leq 2^{\abs{2s+\gamma}} \comii{v-\tilde v}^{\abs{2s+\gamma}},
    \end{eqnarray*}
    the last inequality following from  Peetre's inequality \eqref{peetre}. Similarly,
      \begin{eqnarray*}
    \frac{ \eps^{-1} \comii
       v^{2s+\gamma}}{ \eps a_K(\tilde v,\tilde\eta) +\eps^{-1} \comii
       {\tilde v}^{2s+\gamma}}\leq \frac{\eps^{-1} \comii v^{2s+\gamma}}{ \eps^{-1} \comii
       {\tilde v}^{2s+\gamma} }\leq 2^{\abs{2s+\gamma}} \comii{v-\tilde v}^{\abs{2s+\gamma}}.
    \end{eqnarray*}
    The above inequalities yield, for any $(v,\eta), (\tilde v,\tilde\eta)\in\mathbb R^6,$
    \begin{eqnarray*}
   \frac{\eps a_K(v,\eta)+\eps^{-1} \comii v^{2s+\gamma}}{ \eps a_K(\tilde v,\tilde\eta) +\eps^{-1} \comii
       {\tilde v}^{2s+\gamma}}\leq  \inner{C+2^{1+\abs{2s+\gamma}} } \inner{\comii{v-\tilde v}+\comii{\eta-\tilde \eta}}^{N+\abs{2s+\gamma}}.
    \end{eqnarray*}
  Observe the constant $C$  above is independent of $\eps$, and   thus  $  \eps a_K +\eps^{-1} \comii
       v^{2s+\gamma} $ is a temperate weight uniformly with respect to $\eps. $

Now we will prove the conclusion in the corollary.  This is just a consequence of the conclusion
(i) in Lemma \ref{inverse}.    In fact   note that $K\geq K_0$ with  $K_0$ the constant given in Lemma \ref{inverse}, and thus   $K+\eps \geq K_0$. Then the assumption in Lemma  \ref{inverse}  is fulfilled and we may apply the  conclusion
(i) in    Lemma \ref{inverse} to conclude that $a_{K+\eps^{-2}}^w$   is invertible  and  its inverse   has the form
\[
 \inner{a_{K+\eps^{-2}}^w}^{-1}= \inner{a_K ^w+ \eps^{-2}
\comii v^{2s+\gamma}}^{-1} =  \inner{a_{K+\eps^{-2}}^{-1}}^w H
\]
with $H$ a bounded operator in $L^2.$
 The assumption on $g$ shows
\[
\eps^{-1} g \in S\inner { a_K +\eps^{-2} \comii
       v^{2s+\gamma},~\Gamma}, \]
and thus we can write
\[g^w  = \underbrace{\inner{\eps^{-1}g}^w   \inner{a_{K+\eps^{-2}}^{-1}}^w H
 }_{\in~\mathcal B(L^2) }\eps \inner{a_K^w + \eps^{-2}
\comii v^{2s+\gamma}},   \]
which yields the desired estimate for $g^w$.   The estimate for
$g(v,D_v)$ is similar, since
$g(v,D_v)=\inner{J^{-1/2}g}^w$ with $J^{-1/2}g$ belonging to the same
symbol class as $g$.  We have obtained  the desired estimate in Corollary
\ref{teccor}.
The proof is complete.
\fin

We will apply the preceding lemma to specific pseudo-differential operators:

\begin{lem} \label{pseudsmallsymb}
The symbols of   $
a_s(v,D_v)$ and $   a^w - a(v,D_v)$ lie in  $ S\inner {\eps a
  +\eps^{-1} \comii v^{2s+\gamma},~\Gamma} $ for all $\eps>0$ with seminorms independent of $\eps$.
\end{lem}

\preuve For the first operator $a_s(v,D_v)$, this is point ii) of Proposition \ref{as2}.
For the second one $a^w - a(v,D_v)$, we need more facts from the theory of Weyl and classical quantizations. In order to get the result, we use the expansion of $J^{1/2} a$,   which
reads (c.f. \cite[Lemma 4.1.5]{MR2599384} and the appendix)
\[
    a^w-a(v,D_v)=\inner{J^{1/2} a}(v, D_v)-a(v,D_v)=R(v,D_v)
\]
with
\[
R(v,\eta)= \frac{1}{2}\int
   \inner{ J^{\theta/2} \inner{D_\eta\cdot\partial_v  a} }(v, \eta)  d\theta.
\]
Proposition \ref{propsym} implies that $D_\eta\cdot\partial_v  a \in S\inner
{M_\eps,~\Gamma }$ uniformly with respect to
$\eps$,  where
\[ M_\eps=\eps \tilde a
  +\eps^{-1} \comii v^{2s+\gamma}.\]
Then  proceeding as in the proof of \cite[Lemma
4.1.2]{MR2599384},  we conclude that $J^{\theta/2}
\inner{D_\eta\cdot\partial_v  a}$  belongs to the same symbol
class $S\inner
{M_\eps,~\Gamma } $ as $D_\eta\cdot\partial_v  a$,  due to the fact
that
\[M_\eps (v+z,\eta+\zeta) \leq C M_\eps (v, \eta) H(\comii
z,\comii\zeta)\]
with $H(\comii
z,\comii\zeta)$ being some polynomial of $\comii
z,\comii\zeta$ and $C$ a constant independent of $\eps$.   Observe
$\tilde a \lesssim  a_K$ due to Proposition \ref{estaa} i).  Then we have
proven that the classical symbol of the difference $a(v,D_v)-a^w$ lies in $ S\inner {\eps a
  +\eps^{-1} \comii v^{2s+\gamma},~\Gamma}$. The Weyl symbol therefore also belongs to this class by
  direct transformation. The proof is complete.

\begin{prop}\label{120412}
  Let $\xi$ be the dual variable of $x$ and let $\ell$ be an arbitrarily real number.  Then for any
  $\eps$, there exists a constant $C_\eps$   such that
\begin{equation} \label{11052931}
       \forall~ f\in\mathcal S(\mathbb R_v^3),\quad\norm{\comii v^{2s+\gamma}f}
_{L^2}    \leq    \eps \norm{ a_K^{w} f}_{L^2}+C_\eps \inner{
  \norm{(iv\cdot\xi-\mathcal L) f}_{L^2 }+\norm{ f}_{L_\ell^2 }}.
\end{equation}
Recall $a_K$ is given in \eqref{deak}.
\end{prop}

\preuve
   Let us first recall the coercivity estimate (see for instance
   Theorem 1.1 and
   Proposition 2.2 in \cite{AMUXY1} and \cite{mou, mou1}) : for $0<s<1$ and $\gamma>-3$,
\begin{equation*}
\forall~f\in \mathcal S(\mathbb R^6),\quad \norm{\comii
       v^{s+\frac{\gamma}{2}} ({\bf Id}-{\bf P}) f}_{L^2}^2 \lesssim
     \inner{-\mathcal Lf,~f}_{L^2},
\end{equation*}
where ${\bf Id}$ stands for the identity operator and ${\bf P}$ is the
$L^2$-orthogonal projection onto the null space $${\rm Span}
\set{\mu^{1/2}, v_1\mu^{1/2}, v_2\mu^{1/2}, v_3\mu^{1/2}, \abs{v}^2\mu^{1/2}}.$$
Consequently we have, for any $\ell\in\R,$
\begin{equation}\label{innest}
     \forall~f\in \mathcal S(\mathbb R^6),\quad \norm{\comii
       v^{s+\frac{\gamma}{2}} f}_{L^2}^2 \lesssim {\rm Re}
     \inner{(iv\cdot\xi-\mathcal L)f,~f}_{L^2}+\norm{\comii v^{\ell-s-\gamma/2}f}_{L^2}^2.
\end{equation}
Now  applying estimate \reff{innest} to the function $\comii v^{s+\frac{\gamma}{2}}f$ yields
\begin{eqnarray*}
     &&\norm{ \comii{v}^{2s+\gamma }f}_{L^2}^2  \lesssim {\rm Re}\,\biginner{ (iv\cdot\xi-\mathcal L)  \comii{v}^{s+\frac{\gamma}{2}}f,~
\comii{v}^{s+\frac{\gamma}{2}}f}_{L^2} +\norm{f}_{L_\ell^2}^2\\
       &\lesssim&  \abs{\biginner{ (iv\cdot\xi-\mathcal L)
         f,~\comii{v}^{2s+\gamma}f}_{L^2}} +\abs{\biginner{\comii{v}^{s+\frac{\gamma}{2}}
         \com{ \mathcal L ,~ \comii{v}^{s+\frac{\gamma}{2}}} f,
         ~f}_{L^2}}+\norm{f}_{L_\ell^2}^2,
\end{eqnarray*}
and therefore
\begin{equation}\label{firwei}
     \norm{ \comii{v}^{2s+\gamma }f}_{L^2}^2 \lesssim \norm{ (iv\cdot\xi-\mathcal L) f}_{L^2}^2 +\norm{f}_{L_\ell^2}^2 +\abs{\biginner{\comii{v}^{s+\frac{\gamma}{2}}
         \com{ \mathcal L ,~ \comii{v}^{s+\frac{\gamma}{2}}} f,~f}_{L^2}}.
\end{equation}
We have  to treat  the last term in the above estimate,  which is
bounded from above by
\begin{equation} \label{splitL}
   \abs{\biginner{\comii{v}^{s+\frac{\gamma}{2}}
         \com{ a^w ,~ \comii{v}^{s+\frac{\gamma}{2}}} f,~f}_{L^2}}+\abs{\biginner{\comii{v}^{s+\frac{\gamma}{2}}
         \com{ \mathcal K ,~ \comii{v}^{s+\frac{\gamma}{2}}} f,~f}_{L^2}}.
\end{equation}
We apply
\reff{etader} and  \cite[Theorem 2.3.8]{MR2599384} to conclude that
for any $\eps\in ]0,1[$  the
symbol of the operator
\begin{eqnarray*}
 \comii v^{-(2s+\gamma-1)}\comii{v}^{s+\frac{\gamma}{2}}
   \com{ a^w,~ \comii{v}^{s+\frac{\gamma}{2}}}
\end{eqnarray*}
belongs to
  \[
S\inner{  \eps a_K+\eps^{-1}  \comii v^{2s+\gamma},~\Gamma}
  \]
uniformly with respect to $\eps$.   Then Corollary \ref{teccor} gives,
with $\tilde \eps$ arbitrarily small,
\begin{eqnarray*}
\abs{\biginner{\comii{v}^{s+\frac{\gamma}{2}}
         \com{ a^w,~ \comii{v}^{s+\frac{\gamma}{2}}} f,~f}_{L^2}} &\lesssim& \inner{ \eps \norm{a_K^w
         f}_{L^2}+\eps^{-1}\norm{\comii v^{2s+\gamma}f}_{L^2}} \norm{\comii
       v^{2s+\gamma-1}f}_{L^2}\\
&\lesssim&  \eps \norm{a_K^w
         f}_{L^2}^2 + \tilde \eps\norm{\comii v^{2s+\gamma}f}_{L^2}^2+ C_{\eps,\tilde\eps}
       \norm{ f}_{L_\ell^2}^2,
\end{eqnarray*}
 where in the last inequality we used the interpolation  inequality:
\[
  \norm{\comii{v}^{ 2s+ \gamma-1 } f}_{L^2}\leq \tilde \eps
  \norm{\comii{v}^{ 2s+ \gamma} f}_{L^2}+C_{\tilde \eps}
  \norm{f}_{L_\ell^2}.
\]
Now we have to deal with the operator
\begin{equation*} 
 \comii{v}^{s+\frac{\gamma}{2}}
   \com{ \mathcal K,~ \comii{v}^{s+\frac{\gamma}{2}}}
\end{equation*}
in \reff{splitL}.
For this purpose, we split $\kkk$ into three parts :
\begin{eqnarray} \label{splitK}
    {\mathcal K}&=& \underbrace{-\mathcal L_2-\bar {\mathcal
      L}_{1,\delta,a}}_{\mathcal K_{small}}\underbrace{- \mathcal
    L_{1,3,\delta}-{\mathcal L}_{1,4,\delta}}_{\mathcal K_{mult}}+\underbrace{a_s(v,D_v)+\inner{a(v,D_v)-a^w}}_{\mathcal K_{pseudo}}.
\end{eqnarray}
For the second part $\kkk_{mult}$, the estimate is easy since, as recalled in lemma \ref{L1easy} and \ref{l14d}, operators $\mathcal
    L_{1,3,\delta}$ and ${\mathcal L}_{1,4,\delta}$ commute with the multiplication with a function of $v$. We therefore have
\begin{eqnarray*} 
\abs{\biginner{\comii{v}^{s+\frac{\gamma}{2}}
         \com{ \kkk_{mult},~ \comii{v}^{s+\frac{\gamma}{2}}} f,~f}_{L^2}} &=& 0.
\end{eqnarray*}
For the first part $\kkk_{small}$ of $\kkk$ in (\ref{splitK}), we expand the commutators and use
Cauchy-Schwarz inequality to get
\begin{equation*} 
\begin{split}
& \abs{\biginner{\comii{v}^{s+\frac{\gamma}{2}}
         \com{ \kkk_{small},~ \comii{v}^{s+\frac{\gamma}{2}}}
         f,~f}_{L^2}} \\
&\lesssim C_\eps \norm{\comii{v}^{-s-\frac{\gamma}{2}}
         \com{ \kkk_{small},~ \comii{v}^{s+\frac{\gamma}{2}}}
         \comii{v}^{-(s+\gamma)}\comii{v}^{s+\gamma}f}_{L^2}^2+\eps
       \norm{ \comii{v}^{2s+\gamma}f}^2 \\
&\lesssim C_\eps \norm{\comii{v}^{-s-\frac{\gamma}{2}}
         \com{ \kkk_{small},~
           \comii{v}^{-\frac{\gamma}{2}}}\comii{v}^{s+\gamma}f}_{L^2}^2+C_\eps\norm{
         \com{ \kkk_{small},~
           \comii{v}^{-s-\gamma}}\comii{v}^{s+\gamma}f}_{L^2}^2\\
&\quad +\eps \norm{ \comii{v}^{2s+\gamma}f}^2 \\
& \lesssim C_\eps\inner{ \norm{ \comii{v}^{-s-\frac{\gamma}{2}} \mathcal
    L_{2} \comii{v}^{-\frac{\gamma}{2}} \comii{v}^{s+\gamma}f}^2 +  \norm{ \comii{v}^{-s-\gamma}   \mathcal
    L_{2}\comii{v}^{s+\gamma} f}^2+  \norm{  \mathcal
    L_{2}\comii{v}^{-s-\gamma}  \comii{v}^{s+\gamma} f}^2} \\
    &\quad+ C_\eps \norm{\comii{v}^{-s-\frac{\gamma}{2}}
         \com{ \bar {\mathcal
      L}_{1,\delta,a},~
           \comii{v}^{-\frac{\gamma}{2}}}\comii{v}^{s+\gamma}f}_{L^2}^2+C_\eps\norm{
         \com{ \bar {\mathcal
      L}_{1,\delta,a},~
           \comii{v}^{-s-\gamma}}\comii{v}^{s+\gamma}f}_{L^2}^2\\
&\quad +\eps \norm{ \comii{v}^{2s+\gamma}f}^2 .
 \end{split}
\end{equation*}
Then, we use Lemma \ref{l2l2} and  conclusion (ii) in Lemma \ref{lb1da}
with either $\tilde\alpha=-s-\gamma/2$, $\tilde \beta=-\gamma/2$ or
$\tilde \alpha=0$, $\tilde \beta = -s-\gamma$ (for which we have in both cases $\tilde\alpha+\tilde\beta +\gamma+s \leq 0$) and we get
\begin{equation*} 
\begin{split}
 \abs{\biginner{\comii{v}^{s+\frac{\gamma}{2}}
         \com{ \kkk_{small},~ \comii{v}^{s+\frac{\gamma}{2}}} f,~f}_{L^2}} & \lesssim \tilde{C}_\eps \norm{  \comii{v}^{s+\gamma}f}^2 + \eps \norm{ \comii{v}^{2s+\gamma}f}^2 \\
 &\lesssim \widetilde{C}_\eps \norm{  f}_{L^2_\ell}^2 + 2\eps \norm{ \comii{v}^{2s+\gamma}f}^2 \\
 \end{split}
\end{equation*}
since $s>0$.

Next we deal with the last part $\kkk_{pseudo}$ of $\kkk$ in  (\ref{splitK}). From Lemma
\ref{pseudsmallsymb}, we already know that
$\kkk_{pseudo}$
belongs to
  \[
 S\inner {\eps a
  +\eps^{-1} \comii v^{2s+\gamma},~\Gamma}
  \]
with uniform semi-norms with respect to $\eps$. We follow the same strategy as in the lines just after (\ref{splitL}) for commutators involving $a^w$.
Using that $\D_v \seq{v}^\mu = \ooo( \seq{v}^{\mu-1})$ for all $\mu \in \R$, and
applying \cite[Theorem 2.3.8]{MR2599384} (see also appendix), we get that
for any $\eps\in ]0,1[$  the
symbol of the operator
\begin{eqnarray*}
 \comii v^{-(2s+\gamma-1)}\comii{v}^{s+\frac{\gamma}{2}}
   \com{ \kkk_{pseudo},~ \comii{v}^{s+\frac{\gamma}{2}}}
\end{eqnarray*}
belongs to
  \[
S\inner{  \eps a_K+\eps^{-1}  \comii v^{2s+\gamma},~\Gamma}
  \]
uniformly with respect to $\eps$.   Then Corollary \ref{teccor} gives,
with $\tilde \eps$ arbitrarily small,
\begin{eqnarray*}
\abs{\biginner{\comii{v}^{s+\frac{\gamma}{2}}
         \com{ \kkk_{pseudo},~ \comii{v}^{s+\frac{\gamma}{2}}} f,~f}_{L^2}} &\lesssim& \inner{ \eps \norm{a_K^w
         f}_{L^2}+\eps^{-1}\norm{\comii v^{2s+\gamma}f}_{L^2}} \norm{\comii
       v^{2s+\gamma-1}f}_{L^2}\\
&\lesssim&  \eps \norm{a_K^w
         f}_{L^2}^2 + \tilde \eps\norm{\comii v^{2s+\gamma}f}_{L^2}^2+ C_{\eps,\tilde\eps}
       \norm{ f}_{L_\ell^2}^2.
\end{eqnarray*}
Combining these estimates we obtain
\begin{eqnarray*}
\abs{\biginner{\comii{v}^{s+\frac{\gamma}{2}}
         \com{ \mathcal K,~ \comii{v}^{s+\frac{\gamma}{2}}} f,~f}_{L^2}} \lesssim\eps \norm{a_K^w
         f}_{L^2}^2 + \tilde \eps\norm{\comii v^{2s+\gamma}f}_{L^2}^2+ C_{\eps,\tilde\eps}
       \norm{ f}_{L_\ell^2}^2.
\end{eqnarray*}
Now taking into account \reff{firwei},  the desired estimate \reff{11052931}
follows if we choose $\tilde\eps$ small enough. The proof is thus complete.
\fin

In order to  prove  the main result, Proposition \ref{prpmain}, we
will need the conclusion in Proposition \ref{estaa}. So let us firstly
present the proof of this Proposition.

\preuve[ of Proposition \ref{estaa} ii) and iii)]
We have shown   Proposition \ref{estaa} iii)  in Lemma \ref{inverse}.  For
the conclusion ii),  let us rewrite  the linearized Boltzmann operator $\mathcal L$ as
\begin{eqnarray*}
    \lll
&=& -a^w+\underbrace{\mathcal L_2+\bar {\mathcal
      L}_{1,\delta,a}+\mathcal
    L_{1,3,\delta}+{\mathcal
      L}_{1,4,\delta}-a_s(v,D_v)-(a(v,D_v)-a^w)}_{-\mathcal K}.
\end{eqnarray*}
As a direct consequence of Lemma \ref{l2l2}, conclusion (i) in Lemma \ref{lb1da}, Lemma \ref{L1easy},
Lemma  \ref{l14d} we have 
\[
     \norm{\inner{\mathcal L_2+\bar {\mathcal
      L}_{1,\delta,a}+ \mathcal
    L_{1,3,\delta}+ \mathcal
    L_{1,4,\delta}}f}_{L^2}\lesssim \norm{\comii v^{2s+\gamma}f}_{L^2}.
\]
Moreover from
Lemma \ref{pseudsmallsymb}  we know that for any $\eps>0,$
\[
-\kkk_{pseudo} = -  a_s(v,D_v) - 
 (a^w - a(v,D_v)) \in Op_{weyl} \inner {\eps a +\eps^{-1} \comii v^{2s+\gamma},~\Gamma}\]
uniformly with respect to $\eps$, and thus
 \[
     \norm{\kkk_{pseudo} f}_{L^2} \lesssim \eps \norm{a_K^w f}+\eps^{-1} \norm{\comii v^{2s+\gamma}f}_{L^2}
\]
due to Corollary \ref{teccor}.

The proof of point ii) of Proposition \ref{estaa} is
complete.
\fin

The rest  of this subsection is devoted to the

\preuve[ of Proposition \ref{prpmain}]
 Now assuming that \reff{resulting} holds, we have
\[
   \norm{  \tilde a (v,\xi)^{\frac{1}{1+2s}} f }+ \norm{ a_K^w f } \lesssim \norm{
  \inner{iv\cdot\xi -\mathcal L}f}_{L^2}
+\norm{f}_{L_\ell^2}+\norm{\inner{iv\cdot\xi -\mathcal L-
  \hat P_K}  f}_{L^2}.
\]
On the other hand, note that
\[
      iv\cdot\xi -\mathcal L-
  \hat P_K = a^w +\mathcal K-(a+K \comii v^{2s+\gamma})^w=\mathcal K-K\comii v^{2s+\gamma},
\]
and thus  Proposition \ref{estaa} yields, with $ \eps$
arbitrarily small,
\begin{eqnarray*}
    \norm{\inner{iv\cdot\xi -\mathcal L-
  \hat P_K}  f}_{L^2}  &\lesssim&  \eps
\norm{a_K^w  f}_{L^2}+C_{ \eps} \norm{\comii v^{2s+\gamma} f}_{L^2}\\
&\lesssim&  \eps \norm{ a_K^w f
   }+C_\eps\inner{ \norm{(iv\cdot\xi-\mathcal L)f}_{L^2}+\norm{ f}_{L_\ell^2}},
\end{eqnarray*}
the last inequality  following from  \reff{11052931} . Combining these
inequalities and letting the above  $\eps$ be sufficiently small, we get
\begin{eqnarray*}
   \norm{  \tilde a (v,\xi)^{\frac{1}{1+2s}} f }+ \norm{ a_K^w f }
   \lesssim   \norm{(iv\cdot\xi-\mathcal L)f}_{L^2}+\norm{ f}_{L_\ell^2}.
\end{eqnarray*}
Taking into account the facts that
\[
\comii v^{\gamma/(2s+1)}\comii \xi^{2s/(2s+1)}+\comii v^{\gamma/(2s+1)}\comii{v\wedge\xi}^{2s/(2s+1)}\lesssim \tilde a(v,\xi)^{1/(2s+1)}
 \]
 and that
 \[
   \norm{\comii v^{\gamma}\comii {D_v}^{2s} f}_{L^2}+\norm{ \comii
        v^{\gamma}\comii{v\wedge D_v}^{2s}f}_{L^2}+ \norm {\comii
        v^{2s+\gamma} f}_{L^2} \lesssim
    \norm{a_K^{w}f}_{L^2}
 \]
due to the conclusion (i) in Lemma \ref{inverse},  we obtain the desired estimate in Theorem \ref{thmain}.  The proof of  Proposition  \ref{prpmain} is complete.
\fin

\subsection{Proof of Theorem \ref{coermou} and boundedness estimates}\label{subsec43}

In this section we prove first Theorem \ref{coermou} about coercivity. As mentioned in the introduction it can be understood as an exact
estimate for
the so called triple norm introduced in \cite{AMUXY1} and recalled in
Remark \ref{remarktriple} below. It involves the
pseudo-differential part $a^w$, for which we  have elliptic properties
stated  in Proposition \ref{estaa}. Theorem \ref{coermou} is a direct consequence of the following Lemma:
\black

\begin{lem} \label{coerest}
We have for a sufficiently large constant $C$ and for all $l \in \R$
with $l \leq \gamma/2+s$,
\begin{multline*}
 \norm{ \seq{v}^{\gamma/2 } \seq{D_v}^{s}f }^2
 + \norm{\seq{v}^{\gamma/2} \seq{v \wedge D_v}^{s}f}^2
 + \norm{\seq{v}^{\gamma/2 + s}f}^2 \\
  \approx  \sep{ a^w f , f } + C \norm{\seq{v}^{\gamma/2 + s}f}^2
  \approx -\sep{ \lll f, f} + \norm{\seq{v}^{l}f}^2,
\end{multline*}
where in the last equivalence the constant  depends on $l $.
\end{lem}

\preuve  We first show the   second equivalence. To do so
rewrite  the linearized Boltzmann operator $\mathcal L$ as
\begin{eqnarray*}
    \lll
&=& -a^w+\underbrace{\mathcal L_2+\bar {\mathcal
      L}_{1,\delta,a}+\mathcal
    L_{1,3,\delta}+{\mathcal
      L}_{1,4,\delta}-a_s(v,D_v)-(a(v,D_v)-a^w)}_{-\mathcal K}.
\end{eqnarray*}
As a direct consequence of Lemma \ref{l2l2}, conclusion (i) in Lemma \ref{lb1da}, Lemma \ref{L1easy},
Lemma  \ref{l14d} we have 
\[
    \abs{ \inner{\inner{\mathcal L_2+\bar {\mathcal
      L}_{1,\delta,a}+ \mathcal
    L_{1,3,\delta}+ \mathcal
    L_{1,4,\delta}}f, ~f}_{L^2}}\lesssim \norm{\comii v^{\gamma/2+s}f}_{L^2}^2.
\]
Moreover from \reff{asii} and
Lemma \ref{pseudsmallsymb}  we know that
\[
-\kkk_{pseudo} = -  a_s(v,D_v) - 
 (a^w - a(v,D_v)) \in {\red {\rm Op}  \inner {\eps a_K+\eps^{-1} \comii v^{\gamma+2s},~\Gamma}},\]
and thus  for any $\eps>0,$
 \[
    \abs{ \inner{\kkk_{pseudo} f,~f}_{L^2} }\lesssim \eps \norm{\inner{a_K^{1/2}}^w f}^2+\eps^{-1} \norm{\comii v^{\gamma/2+s}f}_{L^2}^2
\]
due to (ii) Lemma \ref{inverse}. Combining these estimates we conclude
\begin{eqnarray*}
-\sep{ \lll f, ~f}_{L^2} =\sep{ a^w f, ~f} _{L^2}+\sep{ \mathcal K f, ~f}_{L^2}
\end{eqnarray*}
with
\begin{eqnarray*}
\abs{\sep{ \mathcal K f, ~f}_{L^2}}\lesssim \eps \norm{\inner{a_K^{1/2}}^w f}^2+\eps^{-1} \norm{\comii v^{\gamma/2+s}f}_{L^2}^2.
\end{eqnarray*}
Moreover by  \reff{coercivity+} we have
\begin{eqnarray*}
\sep{ a^w f, ~f} _{L^2}+\inner{K\comii v^{\gamma+2s}f, f}_{L^2}=\sep{ a_ K^w f, ~f}_{L^2} = \norm{\inner{a_K^{1/2}}^w f}^2
\end{eqnarray*}
and thus  choosing $\eps$ small enough, we get
\begin{eqnarray*}
-\sep{ \lll f, ~f}_{L^2} =\sep{ a^w f, ~f} _{L^2}+\sep{ \mathcal{\tilde K} f, ~f}_{L^2}
\end{eqnarray*}
with
\begin{eqnarray*}
\abs{\sep{ \mathcal {\tilde K} f, ~f}_{L^2}}\lesssim   \norm{\comii v^{\gamma/2+s}f}_{L^2}^2.
\end{eqnarray*}
As a result, combining the estimate (see \reff{innest} for instance)
\begin{eqnarray*}
 \norm{\comii v^{\gamma/2+s}f}_{L^2}^2\lesssim -\sep{ \lll f, ~f}_{L^2} +\norm{\comii v^\ell f}_{L^2}^2
\end{eqnarray*}
we obtain the second  equivalence.

Next we show the first equivalence.   Using the estimate \reff{coercivity+} below,  we see that
\begin{eqnarray*}
\sep{ a^w f, ~f} _{L^2}+C\norm {K\comii v^{\gamma+2s}f, f}_{L^2} = \sep{ a_ K^w f, ~f}_{L^2} \approx \norm{\inner{a_K^{1/2}}^w f}^2.
\end{eqnarray*}
On the other hand  by conclusion (ii) in Lemma \ref{inverse} we
can deduce that
\[
  \norm{\inner{a_K^{1/2}}^w f}_{L^2}^2 \approx  \norm{ \seq{v}^{\gamma/2 } \seq{D_v}^{s}f }^2
 + \norm{\seq{v}^{\gamma/2} \seq{v \wedge D_v}^{s}f}^2
 + \norm{\seq{v}^{\gamma/2 + s}f}^2.
\]
Then the first equivalence follows, and the proof is complete.
\fin

\remark \label{remarktriple}
 In \cite{AMUXY1}, the authors  introduced  the following non-isotropic norm
 \begin{eqnarray*}
   \normm f^2\stackrel{\rm def}{=} \iiint \Phi(\abs{v-v_*}) b(\cos\theta) \mu_*
\inner{f-f'}^2+\iiint \Phi(\abs{v-v_*}) b(\cos\theta) f_*^2 \inner{\sqrt{\mu'}-
\sqrt{\mu}}^2,
 \end{eqnarray*}
 where the integration is over $\mathbb R_v^3\times\mathbb R_{v_*}^3\times\mathbb
S_\sigma^2$. For such a norm,  Theorem 1.1 of \cite{AMUXY1}) says, with
$l\in\mathbb R$ arbitrary (and equivalence norm depending on $l$),
 \begin{equation*} 
    \norm{\seq{v}^{\gamma/2 } \seq{D_v}^{s}f}^2+\norm{\seq{v}^{\gamma/2 +s }f}\lesssim \normm{ f }^2 \lesssim    -\biginner{\mathcal L f,~f} +C_2 \norm{\seq{v}^{l}f}^2,
 \end{equation*}
  provided the Boltzmann cross-section $B$ satisfies \reff{cross-section} with $0
< s < 1$ and $\gamma > -3$.
In Lemma \ref{coerest} above, we are therefore able to exhibit the complete form of this triple norm $ \normm{ f }$.


\bigskip
Now we focus on the more difficult subelliptic estimate stated in \ref{thmain}. We begin with another coercivity estimate for the  Weyl quantization $a_K^w$.

\begin{lem}\label{coK}
    Let $\hat P_K$ be the operator defined at the beginning of
    Subsection \ref{subsec42}. Then there exists a positive number $k_0> 0$   such that for
    all $K\geq k_0$  and any $f\in\mathcal S(\mathbb R^3),$  we have
    \begin{equation}\label{coercivity+}
   \norm{
    \inner{a_K^{1/2}}^w f }^2 \approx    \inner{
  a_K^w  f,~f}_{L^2} =  {\rm Re} \inner{
  \hat P_K  f,~f}_{L^2}
    \end{equation}
and
   \begin{equation}\label{tildecoer}
   \norm{
    \inner{\tilde a_K^{1/2} a_K^{1/2}}^w f }^2 \approx  \inner{
  \inner{\tilde a_K a_K}^w  f,~f}_{L^2}.
    \end{equation}
\end{lem}

\preuve
The arguments are similar to the ones used in the proof of Lemma \ref{inverse}.
Together with  \reff{11053010} and \reff{expan}, we may write
\begin{equation} \label{+inverse+}
   (a_K^{1/2})^w (a_K^{1/2})^w=  a_K^w- R^w,
\end{equation}
where
\begin{eqnarray*}
   R=-\int_0^1\inner{ \partial_\eta (a_K^{1/2})}\sharp_\theta
   \inner{\partial_v( a_K^{1/2})}  d\theta+\int_0^1 \inner{\partial_v(
     a_K^{1/2})} \sharp_\theta
   \inner{\partial_\eta (a_K^{1/2})}  d\theta
\end{eqnarray*}
with $g \sharp_\theta h$ defined in \reff{sharptheta}.
Using  \reff{etader} for $\eps=K^{-1/4},$  we conclude that
$$\partial_\eta (a_K^{1/2}) \in S(
K^{-1/4}a_K^{1/2},~\Gamma )$$  uniformly
with respect to $K$.  On the other hand, it is clear that  $\partial_v (a_K^{1/2}) \in S(
a_K^{1/2},~\Gamma )$.  As a result,  \cite[Proposition 1.1]{MR1721321}
yields
\[
   \inner{ \partial_\eta (a_K^{1/2})}\sharp_\theta
   \inner{\partial_v( a_K^{1/2})} , \inner{\partial_v(
     a_K^{1/2})} \sharp_\theta
   \inner{\partial_\eta (a_K^{1/2})} \in S(K^{-1/4}
a_K,~\Gamma )
\]
uniformly w.r.t. $K$. Thus  $R \in S(K^{-1/4}
a_K,~\Gamma )$ uniformly w.r.t. $K$.   Then    the conclusion (ii) in Lemma  \ref{inverse}
allows us to rewrite $R^w$ as
\[
    R^w= K^{-1/4}(a_K^{1/2})^w \underbrace{K^{1/2}\com{(a_K^{1/2})^w}^{-1} R^w
      \com{(a_K^{1/2})^w}^{-1}}_{\in \mathcal B(L^2)~\textrm{uniformly
      w.r.t.}~K} (a_K^{1/2})^w,
\]
  which gives
\[
    \abs{\inner{R^w f,~f}_{L^2}} \leq C_0 K^{-1/4}   \norm{(a_K^{1/2})^w
      f}^2
\]
with  $C_0$ some  constant independent of $K$.
Taking into account    the relation \reff{+inverse+}  we obtain
\begin{eqnarray*}
    \inner{ a_K^w f,~f}_{L^2} \leq \inner{(a_K^{1/2})^w
      (a_K^{1/2})^w f,~f}_{L^2}   +  C_0 K^{-1/4}   \norm{(a_K^{1/2})^w
      f}^2\leq (C_0+1)    \norm{(a_K^{1/2})^w
      f}^2
\end{eqnarray*}
and
\begin{eqnarray*}
     \inner{(a_K^{1/2})^w (a_K^{1/2})^w f,~f}_{L^2}  \leq \inner{ a_K^w f,~f}_{L^2}+  C_0 K^{-1/4}   \norm{(a_K^{1/2})^w
      f}^2.
\end{eqnarray*}
The desired estimate \reff{coercivity+}  follows  if we take $K$ sufficiently large such
that $K\geq k_0\stackrel{\rm def}{=} 16C_0^4$.   Since the second
estimate \reff{tildecoer} can be deduced similarly by virtue of (iii)
in Lemma \ref{inverse}, we omit it here. The proof is thus complete.
\fin

\begin{cor}
Let $\ell$  be an arbitrary real number.  The following   estimate
  \begin{equation}\label{coest}
       \forall~ f\in\mathcal S(\mathbb R_v^3),\quad\norm{\comii v^{2s+\gamma}f}
_{L^2}+\norm{ (a_K^{1/2})^{w} \comii{v}^{s+
    \gamma/2  }f}_{L^2}
         \lesssim      \norm{\hat P_K  f}_{L^2 }+\norm{ f}_{L_\ell^2 }
\end{equation}
holds uniformly with respect to $\xi$.  Recall $a_K$ is defined by \eqref{deak}.
\end{cor}

\preuve
     We have obtained in the proof of Lemma \ref{coerest} the estimate
   \begin{eqnarray*}
   \norm{\comii v^{2s+\gamma}f}
_{L^2}\lesssim \norm{ (a_K^{1/2})^{w} \comii{v}^{s+
    \gamma/2  }f}_{L^2}.
   \end{eqnarray*}
 Moreover using the coercivity estimate \reff{coercivity+} applied to the function $\comii
v^{s+\gamma/2}f$, we have
\begin{eqnarray*}
  &&\norm{\comii v^{2s+\gamma}f}
_{L^2}^2+\norm{ (a_K^{1/2})^{w} \comii{v}^{s+
    \gamma/2  }f}_{L^2} ^2\\
    &&
         \lesssim       \abs{\inner{\hat P_K \comii{v}^{s+
    \gamma/2  } f,~\comii{v}^{s+
    \gamma/2  }f}_{L^2 }} \\
    && \lesssim       \abs{\inner{ \com{\hat P_K,~ \comii{v}^{s+
    \gamma/2  } }f,~ \comii{v}^{s+
    \gamma/2  }f}_{L^2 }} + \abs{\inner{  \hat P_K f,~ \comii{v}^{2s+
    \gamma  }f}_{L^2 }} \\
    && \lesssim       \abs{\inner{ \com{a^w,~ \comii{v}^{s+
    \gamma/2  } }f,~ \comii{v}^{s+
    \gamma/2  }f}_{L^2 }} + \eps^{-1}\norm{\hat P_K f}_{L^2}^2+\eps\norm{\comii{v}^{2s+
    \gamma  }f}_{L^2 }^2.
\end{eqnarray*}
We apply
\reff{etader} and  \cite[Theorem 2.3.8]{MR2599384} to conclude that the
symbol of the operator
\begin{eqnarray*}
   \com{ a^w,~ \comii{v}^{s+\frac{\gamma}{2}}}
\end{eqnarray*}
belongs to
  \[
S\inner{ a_K^ {1/2}\comii v^{2s+\gamma-1},~\Gamma}.
  \]
This fact, together with Lemma \ref{inverse} (ii),  allows us to write
\begin{eqnarray*}
 &&\com{ a^w,~ \comii{v}^{s+\frac{\gamma}{2}}}\\
 &&=\epsilon^{-1} \comii{v}^{s-1+\frac{\gamma}{2}}\underbrace{\comii{v}^{-\inner{s-1+\frac{\gamma}{2}}}\com{ a^w,~ \comii{v}^{s+\frac{\gamma}{2}}} \comii{v}^{-\inner{s+\frac{\gamma}{2}}}\left[ \inner{a_K^ {1/2}}^{w}\right]^{-1}}_{\in~\mathcal B(L^2)}\epsilon\inner{a_K^ {1/2}}^{w}\comii{v}^{s+\frac{\gamma}{2}}.
\end{eqnarray*}
Then
\begin{eqnarray*}
 \abs{\inner{ \com{a^w,~ \comii{v}^{s+
    \gamma/2  } }f,~ \comii{v}^{s+
    \gamma/2  }f}_{L^2 }} &\lesssim &\eps\norm{\inner{a_K^ {1/2}}^{w}\comii{v}^{s+\frac{\gamma}{2}}f}_{L^2 }^2+\eps^{-1}\norm{\comii{v}^{2s+\gamma-1} f}_{L^2}^2\\
    &\lesssim &\eps\norm{\inner{a_K^ {1/2}}^{w}\comii{v}^{s+\frac{\gamma}{2}}f}_{L^2 }^2+\eps \norm{\comii{v}^{2s+\gamma} f}_{L^2}^2\\
    &&+C_\eps \norm{\comii{v}^{\ell} f}_{L^2}^2.
\end{eqnarray*}
Letting $\eps$ be small sufficiently gives the conclusions.
\fin

\begin{cor}
\begin{equation}\label{+11051518}
\inner{ \comi{ \comii v^{2s+\gamma}}^{\rm Wick} f,~f}_{L^2}\lesssim \inner{ \comi{\tilde
a(v,\eta)}^{\rm Wick} f,~f}_{L^2}
   \lesssim  \abs{\biginner{\hat P_K  f,~ f}_{L^2}}.
\end{equation}
\end{cor}

\preuve
    The first inequality is due to the positivity of Wick
    quantization. The second one is just an immediate consequence of
    \reff{coercivity+} and Lemma \ref{inverse}, since we may write
  \[
   \comi{\tilde
a(v,\eta)}^{\rm Wick}
= \comi{a_K^{1/2}}^w \underbrace{\com{\comi{a_K^{1/2}}^w}^{-1}\comi{\tilde
a(v,\eta)}^{\rm Wick}\com{\comi{a_K^{1/2}}^w}^{-1}}_{\in\mathcal B(L^2)}\comi{a_K^{1/2}}^w,
 \]
 where we use the fact (see the appendix) that $\tilde a^{\rm Wick}=b^w$ with $b$ belonging to the same symbol class as $\tilde a.$
 The proof is complete.

\fin


\subsection{Hypoelliptic estimates and proof of Theorems \ref{thmain} and \ref{th1}}\label{subsec44}

This last subsection is devoted to the proofs of the main results, Theorem
\ref{thmain} and Theorem \ref{th1}.   As explained in Proposition
\ref{prpmain},  we only work on $\hat P_K$ instead of $P$. Therefore, in this
subsection, $\xi$ and $\tau$ are considered as parameters. Recall
that $\tilde a$ is defined in \reff{atilde},  whose explicit form,  as to be seen
below,  will be of convenient use.   The main result to be shown here can be
stated as follows

\begin{prop} \label{pkreg}
 Under the conditions of Theorem 1, we have, for any $\ell\in\mathbb R,$
$$
\norm{ \tilde a(v,\xi)^{\frac{1}{1+2s}} f }+\norm{ a_K^w f } \lesssim \norm{
  \hat P_K   f}_{L^2} + \norm{\comii v^\ell f}_{L^2}.
$$
 Recall $a_K$ is defined by \eqref{deak}.
\end{prop}

The above proposition will be proved in several steps,  following the multiplier
strategy introduced in \cite{HK2011}.  To this end,  throughout
this section,  we let $\chi\in C_0^\infty(\mathbb R;~[0,1])$ such that
$\chi=1$ in
$[-1,1]$ and supp~$\chi \subset[-2,2]$,  and let $g$ be a  symbol given by
\begin{equation} \label{cutcut}
g(v,\eta) =
  g_{\xi}(v,\eta)=\frac{  a_3(v,\eta)}{\tilde a(v,\xi)^{\frac{2s}{1+2s}}}\psi(v,
\eta),
\end{equation}
where
\begin{equation}\label{psi}
\psi(v,\eta)=\chi\sep{ \frac{\tilde a(v,\eta)}{ \tilde a(v,\xi)^{\frac{1}{1+2s}}}}
\end{equation}
and
\begin{eqnarray}
 a_3(v,\eta) = \comii v^{\gamma}\inner{1+\abs v^2+\abs\xi^2+\abs{v\wedge \xi}^2}
^{s-1}\biginner{\xi\cdot\eta+(v\wedge \xi)\cdot(v\wedge \eta)}.
\end{eqnarray}
The main property  linking  $a_3$ and $\tilde a$ is that
\begin{equation} \label{mainprop}
\set{a_3(v,\eta),  v\cdot\xi} = \tilde a( v, \xi)-\comii
v^{\gamma+2}\inner{1+\abs v^2+\abs\xi^2+\abs{v\wedge
    \xi}^2}^{s-1}.
\end{equation}
where   $\set{\cdot,~\cdot}$ is the Poisson bracket defined in
  \reff{11051505}.
Thanks to the explicit  symbolic estimates  for $\tilde a$,   $g$ and $\psi$ also have good
behavior as  symbols, that is,
\begin{eqnarray*} 
      g , ~ \psi\in S(1, \abs{dv}^2+\abs{d\eta}^2)
      \end{eqnarray*}
       uniformly with respect to $\xi$, where  we use the estimate
       \begin{eqnarray*}
       \abs {a_3(v,\eta)}
       \lesssim \tilde a(v,\xi)^{\frac{2s-1}{2s}} \tilde a(v,\eta)^{\frac{1}{2s}}.
       \end{eqnarray*}
Moreover  direct computation
       shows that
  \begin{eqnarray} \label{11040806}
    \abs{ \xi\cdot\partial_\eta\psi } \lesssim \tilde a(v,\eta).
   \end{eqnarray}

 \begin{lem} \label{crucial} Under the conditions in Theorem 1, we have
$$
\forall ~f\in
\mathcal S(\mathbb R^3),\quad \norm{ \tilde a(v,\xi)^{\frac{1}{1+2s}} f }^2 \lesssim \norm{
  \hat P_K   f}_{L^2}^2+\norm{f}_{L_\ell^2}.
$$
\end{lem}

\preuve
The proof is divided into three steps.

 {\it Step 1)} ~~ Let $g ^{\rm Wick}$ be the Wick quantization of the
 symbol $g $ given in \reff{cutcut}.  We claim
\begin{equation} \label{11051518}
      \abs{\biginner{a_K^w f,~ g ^{\rm Wick}
          f}_{L^2}}\lesssim \abs{\biginner{\hat P_K
          f,~f}_{L^2}}.
  \end{equation}
Indeed,  let us  write,  denoting by   $H$ the inverse
of $\inner{a_K^{1/2}}^w$,
\[
     \biginner{a_K^w f,~ g ^{\rm Wick}
          f}_{L^2}=\biginner{ H a_K^w H   \inner{a_K ^{1\
2}}^w   f,~ \inner{a_K ^{1/2}}^w g^{\rm Wick} H \inner{a_K ^{1/2}}^w
          f}_{L^2}.
\]
Note  that   $H a_K^w H $ and $\inner{a_K ^{1/2}}^w g^{\rm Wick} H $   are
bounded  operators on $L^2$  due to Lemma \ref{inverse} and the fact that $g^{\rm Wick}=\tilde g^w$ with $\tilde g\in S(1,\Gamma)$ (see the appendix).   Then  one
has
\[
     \abs{\biginner{a_K^w f,~ g ^{\rm Wick}
          f}_{L^2}}\lesssim \norm{ \inner{a_K^{1/2}}^{w}
          f}_{L^2}^2 \lesssim  \abs{\biginner{\hat P_K
          f,~f}_{L^2}},
\]
the last inequality following from \reff{coercivity+}.

{\it Step 2)~}
 We now  prove
\begin{equation}\label{11reg}
     \norm{\tilde a(v,\xi)^{\frac{1}{2+4s}} f}_{L^2}
   \lesssim  \norm{\tilde a(v,\xi)^{-\frac{1}{2+4s}}\hat P_K   f}_{L^2}.
\end{equation}
Note that   $g\in S(1,\Gamma)$   and $
   \tilde a(v,\xi)^{r}\in S\inner{\tilde
     a(v,\xi)^{r}, \Gamma}$ for any $r\in\mathbb R$.      Then the
   above estimate will follow if we can show that
  \begin{equation}\label{11041105}
   \norm{\tilde a(v,\xi)^{\frac{1}{2+4s}} f}_{L^2}^2
   \lesssim  \abs{\biginner{\hat P_K   f,~ f}_{L^2}}+\abs{\biginner{\hat P_K   f,~ g
^{\rm Wick}  f}_{L^2}}.
  \end{equation}
To prove the above inequality   we make  use of  the relation
  \begin{eqnarray*}
     {\rm Re}~\biginner{i \inner{v\cdot\xi} f,~g ^{\rm Wick}
       f}_{L^2} = {\rm Re}~\biginner{\hat P_K   f,~g ^{\rm
         Wick} f}_{L^2}-{\rm Re}~\biginner{a_K^w f,~ g ^{\rm Wick} f}_{L^2}
  \end{eqnarray*}
  and \reff{11051518},
  to conclude that
  \begin{equation}\label{11040807}
    {\rm Re}~\biginner{i \inner{v\cdot\xi} f,~g ^{\rm Wick} f}_{L^2}\lesssim
\abs{\biginner{\hat P_K   f,~ f}_{L^2}}+\abs{\biginner{\hat P_K   f,~ g ^{\rm Wick}  f}
_{L^2}}.
  \end{equation}
   Next we will give a lower bound of the term on the left hand side.
   Observe that by \reff{ww},
  \[
    v\cdot\xi= \inner{v\cdot\xi}^{\rm Wick}.
  \]
  Then we have, by \reff{11082406},
  \begin{equation}\label{11040808}
    {\rm Re}~\biginner{i \inner{v\cdot\xi} f,~g ^{\rm Wick} f}_{L^2}=\frac{1}
{4\pi}\biginner{\big\{g,~v\cdot\xi
    \big\}^{\rm Wick}f,~f}_{L^2}.
  \end{equation}
  Using \reff{mainprop} we compute
  \begin{eqnarray*}
     &&\big\{g,~v\cdot\xi\big\}\\
    &=& \tilde a(v,\xi)^{\frac{1}{1+2s}}\psi-\frac{\comii
v^{\gamma+2}\inner{1+\abs v^2+\abs\xi^2+\abs{v\wedge
    \xi}^2}^{s-1}}{\tilde a(v,\xi)^{\frac{2s}{1+2s}}}\psi+\frac{  a_3(v,\eta)}
{\tilde a(v,\xi)^{\frac{2s}{1+2s}}}\xi\cdot\partial_\eta\psi\\  &=& \tilde a(v,
\xi)^{\frac{1}{1+2s}}-\tilde a(v,\xi)^{\frac{1}{1+2s}}\inner{1-\psi}-\frac{\comii
v^{\gamma+2}\inner{1+\abs v^2+\abs\xi^2+\abs{v\wedge
    \xi}^2}^{s-1}}{\tilde a(v,\xi)^{\frac{2s}{1+2s}}}\psi\\
&&+\frac{  a_3(v,\eta)}{\tilde a(v,\xi)^{\frac{2s}{1+2s}}}\xi\cdot\partial_\eta
\psi.
    \end{eqnarray*}
   This   along with \reff{11040807} and \reff{11040808} yields
 \begin{eqnarray}\label{11052920}
   \inner{ \comi{\tilde a(v,\xi)^{\frac{1}{1+2s}}}^{\rm Wick} f,~f}_{L^2}
\lesssim \sum_{j=1}^3 T_j +  \abs{\biginner{\hat P_K   f,~ f}_{L^2}}+
\abs{\biginner{\hat P_K   f,~ g ^{\rm Wick}  f}_{L^2}},
\end{eqnarray}
with
\begin{eqnarray*}
  T_1&=&\inner{\inner{\tilde
      a(v,\xi)^{\frac{1}{1+2s}}\inner{1-\psi}}^{\rm Wick}
    f,~~f}_{L^2},\\
T_2&=&\inner{\inner{\comii
v^{\gamma+2}\inner{1+\abs v^2+\abs\xi^2+\abs{v\wedge
    \xi}^2}^{s-1}\tilde a(v,\xi)^{-\frac{2s}{1+2s}}\psi}^{\rm Wick} f,~~f}_{L^2},
\\
  T_3&=&\inner{\biginner{-\frac{  a_3(v,\eta)}{\tilde a(v,\xi)^{\frac{2s}{1+2s}}}
\xi\cdot\partial_\eta\psi}^{\rm Wick} f,~~f}_{L^2}.
\end{eqnarray*}
Note that
$\tilde a(v,\xi)^{\frac{1}{1+2s}} \leq \tilde a(v,\eta)$
on the support of $ 1-\psi$, and thus
\[
   \tilde a(v,\xi)^{\frac{1}{1+2s}}(1-\psi) \leq  \tilde a(v, \eta).
\]
Then the  positivity of Wick quantization gives
\begin{eqnarray}\label{+11052921}
   T_1\lesssim \comi{\inner{\tilde a(v,\eta)}^{\rm Wick}f,~f}_{L^2} \lesssim
\abs{\biginner{\hat P_K   f,~ f}_{L^2}},
\end{eqnarray}
where the last inequality follows from \reff{+11051518}.
Similarly, observing that
\[
    \comii
v^{\gamma+2}\inner{1+\abs v^2+\abs\xi^2+\abs{v\wedge
    \xi}^2}^{s-1}\tilde a(v,\xi)^{-\frac{2s}{1+2s}}\psi\leq \comii{ v}^{2s+
\gamma},
\]
we have
\begin{eqnarray}\label{K2}
   T_2\lesssim  \inner{\inner{\comii v^{2s+\gamma}}^{\rm Wick}f,~f}_{L^2}
    \lesssim \abs{\biginner{\hat P_K   f,~ f}_{L^2}}.
\end{eqnarray}
As for $T_3$,  it follows from  \reff{11040806} that
\[
   -\frac{  a_3(v,\eta)}{\tilde
     a(v,\xi)^{\frac{2s}{1+2s}}}\xi\cdot\partial_\eta\psi\lesssim
   \tilde a(v,\eta).
\]
Thus
\begin{eqnarray*}
   T_3\lesssim  \comi{\inner{\tilde a(v,\eta)}^{\rm Wick}f,~f}_{L^2}
    \lesssim \abs{\biginner{\hat P_K   f,~ f}_{L^2}}.
\end{eqnarray*}
This, together with \reff{11052920}, \reff{+11052921} and \reff{K2},  gives
\begin{eqnarray*}
    \inner{ \comi{\tilde a(v,\xi)^{\frac{1}{1+2s}}}^{\rm Wick} f,~f}_{L^2}
   \lesssim  \abs{\biginner{\hat P_K   f,~ f}_{L^2}}+\abs{\biginner{\hat P_K   f,~ g
^{\rm Wick}  f}_{L^2}}.
  \end{eqnarray*}
Moreover   by  \reff{ww},
\[
     \comi{\tilde a(v,\xi)^{\frac{1}{1+2s}}}^{\rm Wick}= \int \tilde
     a(v-\tilde v,\xi)^{\frac{1}{1+2s}} e^{-2\pi \tilde v^2} 2^3
     d\tilde v,
\]
which is bounded from below by $\tilde a(v,\xi)^{1/\inner{1+2s}}$ by
a direct check. In fact observe
$$\int \tilde
     a(v-\tilde v,\xi)^{\frac{1}{1+2s}} e^{-2\pi \tilde v^2} 2^3
     d\tilde v\gtrsim  \int_{1/4\leq \abs{\tilde v}\leq 1/2} \tilde
     a(v-\tilde v,\xi)^{\frac{1}{1+2s}}
     d\tilde v,$$  and  for any $1/4\leq \abs{\tilde v}\leq 1/2$ we use  Peetre's inequality to compute
\begin{eqnarray*}
	\tilde
     a(v-\tilde v,\xi)&\gtrsim&  \comii{v-\tilde v}^{\gamma+2s}+\comii{v-\tilde v}^\gamma\inner{\abs{\xi}^{2}+ \abs{(v-\tilde v)\wedge \xi}^2}^{s}\\
   &   \gtrsim  &\comii{v}^{\gamma+2s}+ \comii{v}^{\gamma}\inner{  \abs{\xi}^{2}+\inner{\abs{v \wedge \xi}-\abs\xi/2}^{2} }^s\\
     &   \gtrsim  &\comii{v}^{\gamma+2s}+ \comii{v}^{\gamma}    \inner{\abs{v \wedge \xi}+\abs\xi/2}^{2s}\approx \tilde a(v,\xi) .
\end{eqnarray*}
 As  a result, the desired estimate \reff{11041105} follows.

{\it Step 3)} ~~ Now applying inequality \reff{11041105} to the
function $\tilde a(v,\xi)^{\frac{1}{2+4s}}f$,  we get
\begin{eqnarray*}
  \norm{\tilde a(v,\xi)^{\frac{1}{1+2s}}f}_{L^2}
   &\lesssim&  \norm{\tilde a(v,\xi)^{-\frac{1}{2+4s}}\hat P_K  \tilde
     a(v,\xi)^{\frac{1}{2+4s}} f}_{L^2}\\& \lesssim & \norm{\hat P_K  f}_{L^2}+\norm{\tilde a(v,\xi)^{-\frac{1}{2+4s}}
\com{a_K^w ,~ \tilde
     a(v,\xi)^{\frac{1}{2+4s}}} f}_{L^2}.
\end{eqnarray*}
In view of  \reff{expan},  the symbol of $ \tilde a(v,\xi)^{-1/\inner{2+4s}}\com{a_K^w,~\tilde
       a(v,\xi)^{1/\inner{2+4s}}}$ has the form
\[
     \tilde a(v,\xi)^{-\frac{1}{2+4s}} \int_0^1 \inner{\partial_\eta
       a_K}\sharp_\theta\inner{\partial_v (\tilde
       a^{1/\inner{2+4s}})}d\theta,
\]
which,  arguing as in the proof of  Lemma \ref{inverse},    belongs  to
  \[
       S\inner{ a^{1/2} \comii v^{s+\gamma/2},~\Gamma}.
 \]
As a result, we can use (ii) in Lemma \ref{inverse} to  write
\begin{eqnarray*}
      && \tilde a(v,\xi)^{-\frac{1}{2+4s}}\com{a_K^w  ,~
     \tilde a(v,\xi)^{\frac{1}{2+4s}}}\\
    &=&\underbrace{ \tilde a(v,\xi)^{-\frac{1}{2+4s}} \com{a_K^w  ,~
     \tilde a(v,\xi)^{\frac{1}{2+4s}}} \comii{v}^{-(s+ \gamma/2)  }    \inner{\inner{a_K ^{1/2}}^w}^{-1}
   }_{\in~\mathcal B (L^2)}
   (a_K^{1/2})^{w} \comii{v}^{s+ \gamma/2  }.
  \end{eqnarray*}
This gives
 \begin{eqnarray*}
 \norm{\tilde a(v,\xi)^{-\frac{1}{2+4s}}\com{a_K^w ,~ \tilde
     a(v,\xi)^{\frac{1}{2+4s}}} f}_{L^2} &\lesssim& \norm{  (a_K^{1/2})^{w}
\comii{v}^{s+
     \gamma/2  }f}_{L^2}\\
&\lesssim&\norm{  \hat P_K   f}_{L^2}+\norm{  f}_{L_\ell^2},
  \end{eqnarray*}
where the last inequality follows from \reff{coest}.   Combining
these inequalities, we get the desired estimate
\begin{eqnarray*}
  \norm{\tilde a(v,\xi)^{\frac{1}{1+2s}}f}_{L^2} \lesssim  \norm{  \hat P_K
  f}_{L^2}+\norm{  f}_{L_\ell^2}.
\end{eqnarray*}
 The  proof of
Lemma
\ref{crucial} is thus complete.
\fin



\begin{lem}\label{regvel}
Under the conditions in Theorem 1, we have, for any $\ell\in\mathbb R$,
$$
\norm{  a_K^w f }_{L^2}  \lesssim
   \norm{\hat  P_K f}_{L^2} + \norm{\comii v^{\ell}f}_{L^2}.
$$
\end{lem}

\preuve
The proof is divided into four steps. In the following, let $\eps>0$ be
an arbitrarily small number,  to be fixed later on, and  denote by $C_\eps$ different suitable constants
depending only on $\eps$ and appearing in the the estimations below.

{\it Step 1)}~
We define
$\rho_\eps$   by
  \begin{eqnarray*}
  \rho_\eps(v,\eta)=\chi\left(\frac{\tilde a(v,\xi)^{\frac{1}{1+2s}}}{\eps
\tilde a(v,\eta)}\right),
\end{eqnarray*}
where  $\chi\in C_0^\infty(\mathbb R;~[0,1])$ such that $\chi=1$ in
$[-1,1]$ and supp~$\chi \subset[-2,2]$.
Let $\lambda_{1,\eps}$ and $\lambda_{2,\eps}$  be two symbols defined by
\begin{eqnarray}\label{11082310}
  \lambda_{1,\eps} (v,\eta)= \rho_\eps(v,\eta) \tilde a(v,\eta)
\end{eqnarray}
and
\begin{eqnarray}\label{11082311}
  \lambda_{2,\eps}(v,\eta)=\comi{1-\rho_\eps(v,\eta)}\tilde a(v,\eta).
\end{eqnarray}
Then
$\rho_\eps(v,\eta)\in S\inner{1,~\Gamma}$,
\begin{eqnarray}\label{11082312}
   \lambda_{1,\eps},~\lambda_{2,\eps} \in S\inner{\tilde a(v,
     \eta),~\Gamma} ~~{\rm and}~~\lambda_{2,\eps} \in S\inner{\eps^{-1}\tilde
a(v, \xi)^{\frac{1}{1+2s}},~\Gamma},
\end{eqnarray}
uniformly with respect to $\xi$ and $\eps$, due to the conclusion (i) in
Proposition \ref{estaa} and the fact that
$\tilde a(v, \eta)\leq \eps^{-1}\tilde a(v, \xi)^{\frac{1}{1+2s}}$
on the support of $\lambda_{2,\eps}$.

{\it Step 2)}~~Let
  $\lambda_{1,\eps}(v,\eta)$ be given in \reff{11082310}.   In this
  step we show that
 \begin{eqnarray}\label{11081305}
      \abs{\Big(\com{ v\cdot\xi,~~\lambda_{1,\eps}^w} f,~  f\Big)_{L^2}}
     \leq   \eps \norm{ a_K^w f}_{L^2}^2.
 \end{eqnarray}
 In fact, the  symbol of the above commutator $\com{ v\cdot\xi,~~\lambda_{1,\eps}
^w}$  is
 \[
    -\frac{1}{2i\pi}\xi\cdot \partial_\eta \lambda_{1,\eps}(v,\eta),
 \]
 which belongs to
 $S\inner{\eps^{(1+2s)/2s}\tilde a(v,\eta)^2,~\Gamma}$
 uniformly with respect to $\xi$ and $\eps$,  due to  \reff{derx} and
 the fact that
 \[
    \abs{\xi}+\abs{v\wedge\xi}\lesssim \tilde a(v,\xi)^{\frac{1}{2s}}\comii
    v^{-\frac{\gamma}{2s}}\leq \eps^{\frac{1+2s}{2s}}\tilde a(v,
\eta)^{\frac{1+2s}{2s}}\comii
    v^{-\frac{\gamma}{2s}}
 \]
 on the support of $\lambda_{1,\eps}$.
Thus writing
 \[
      \com{ v\cdot\xi,~~\lambda_{1,\eps}^w}=\eps  a_K^w
\underbrace{ (a_K^w)^{-1}\com{ v\cdot\xi,~~\lambda_{1,\eps}^w} (a_K^w)^{-1}}_{\in\,
\mathcal B(L^2)} a_K^w,
 \]
we obtain
 \begin{eqnarray*}
   \abs{\inner{\com{ v\cdot\xi,~~\lambda_{1,\eps}^w} f,~
       f}_{L^2}} \lesssim \eps\norm{ a_K^w f}_{L^2}^2.
 \end{eqnarray*}
 This gives the desired upper bound and therefore the proof of \reff{11081305}.

{\it Step 3)} ~Let
  $\lambda_{2,\eps}(v,\eta)$ be given in \reff{11082311}. We claim that
 \begin{eqnarray}\label{11082302}
     \abs{\Big(\com{ v\cdot\xi,~~\lambda_{2,\eps}^w} f,~  f\Big)_{L^2}}
     \lesssim  \eps \norm{(v\cdot\xi)f}_{L^2}^2
    +C_\eps \inner{\norm{\hat P_K  f}_{L^2}^2+\norm{ f}_{L_\ell^2}^2}.
 \end{eqnarray}
Indeed, we write $\com{
   v\cdot\xi,~~\lambda_{2,\eps}^w}=
   v\cdot\xi\lambda_{2,\eps}^w-\lambda_{2,\eps}^w  v\cdot\xi$ to get
 \begin{eqnarray*}
     \abs{\inner{ \com{ v\cdot\xi,~~\lambda_{2,\eps}^w} f,~  f}_{L^2}}
     \leq 2\norm{\inner{v\cdot\xi} f}_{L^2}  \norm{\lambda_{2,\eps}^w f}_{L^2}.
 \end{eqnarray*}
 Moreover it follows from  \reff{11082312} that
 \[
    \norm{\lambda_{2,\eps}^w f}_{L^2}\lesssim
      \eps^{-1} \norm{\tilde a(v,\xi)^{1/(1+2s)}f}_{L^2}\lesssim
\eps^{-1}\inner{ \norm{ \tilde {\mathcal P}f}_{L^2}+       \norm{f}_{L_\ell^2}},
 \]
 the last inequality using Lemma \ref{crucial}.
 Combining these inequalities, we obtain the desired estimate  \reff{11082302}.

{\it Step 4)}~~Now we are ready to prove that
\begin{eqnarray}\label{11081802}
   \norm{ a_K^w f}_{L^2}^2 \lesssim    \norm{\hat P_K  f}_{L^2}^2+\norm{ f}_{L_\ell^2}^2 .
 \end{eqnarray}
This inequality will be obtained if we can show that
 \begin{equation}\label{11081802bis}
      \abs{{\rm Re}\inner{i(v\cdot\xi)f,~\tilde a_K^wf}_{L^2}}
      \lesssim  \eps\norm{  a_K^wf}_{L^2}^2 +C_\eps\inner{\norm{\hat P_K  f}_{L^2}^2+
\norm{ f}_{L_\ell^2}^2}
 \end{equation}
 and
 \begin{equation}\label{aakw}
    \norm{a_K^w f}^2 \leq  {\rm Re}~ \inner{a_K ^w
         f, ~ \tilde a_K^w f}_{L^2}+\eps\norm{a_K^w
         f}^2+C_\eps \inner{\norm{\hat P_K f}^2+\norm{f}^2},
\end{equation}
due to the relation
\begin{eqnarray*}
 {\rm Re}\inner{\hat P_K  f,~\tilde a_K^w f}_{L^2} = {\rm
   Re}\inner{i(v\cdot\xi)f,~  \tilde a_K ^wf}_{L^2} +{\rm
   Re}\inner{ a_K^w f,~\tilde a_K^w f}_{L^2}.
\end{eqnarray*}
To prove \reff{11081802bis},  we compute
\begin{eqnarray*}
     &&\abs{{\rm Re}\inner{i(v\cdot\xi)f,~\tilde a_K^wf}_{L^2}}
      =\abs{ \frac{i}{2} \Big(  \com{v\cdot\xi,~  \tilde a_K^w} f,~   f\Big)_{L^2}}=\abs{ \frac{i}{2} \Big(  \com{v\cdot\xi,~  \tilde a^w} f,~   f\Big)_{L^2}}
\\
        &\lesssim& \abs{ \Big( \com{v\cdot\xi,~  \lambda_{1,\eps}^w} f,~   f
\Big)_{L^2}}+\abs{ \Big( \com{v\cdot\xi,~  \lambda_{2,\eps}^w} f,~   f
\Big)_{L^2}}
 \end{eqnarray*}
 with $\lambda_{1,\eps}$,   $\lambda_{2,\eps}$ defined in
 \reff{11082310} and \reff{11082311}.  Combining the above
 inequalities and the conclusion in the previous two steps,  we have
 \begin{eqnarray*}
      \abs{{\rm Re}\inner{i(v\cdot\xi)f,~\tilde a_K^wf}_{L^2}}
      \lesssim  \eps\norm{  a_K^wf}_{L^2}^2+ \eps\norm{(v\cdot\xi)f}_{L^2}^2+C_\eps
\inner{\norm{\hat P_K  f}_{L^2}^2+\norm{ f}_{L_\ell^2}^2}.
 \end{eqnarray*}
This inequality along with the relation
\[
   \norm{(v\cdot\xi)f}_{L^2}^2\lesssim \norm{\hat P_K  f}_{L^2}^2+\norm{ a_K^w f}
_{L^2}^2
\]
 implies the desired estimate \reff{11081802bis}.

 We now   prove \reff{aakw}.  In view of  \reff{expan}   we may write
\begin{equation}\label{lassharp}
   \inner{ \tilde a_K \sharp a_K  }^w =  \inner{ \tilde a_K  a_K  }^w +r^w,
\end{equation}
where
\[
   r(Y)=\int_0^1\iint e^{-2i \sigma(Y-Y_1, Y-Y_2)/\theta}
   \frac{1}{2i}\sigma (\partial_{Y_1}, \partial_{Y_2}) \tilde a(Y_1)
   a_K(Y_2) dY_1 dY_2 d\theta/(\pi\theta)^6.
\]
Note that \reff{etader} also holds true, with $a$ replaced by $\tilde a_K$
or $a_K$.  Then in view of \cite[Proposition 1.1]{MR1721321}, we can check that
\[
   r\in S\inner{ a_K^{3/2}\comii v^{s+\gamma/2}, ~\Gamma},
\]
and thus we may use  Lemma \ref{inverse} to rewrite  $r^w$ as
\[
r^w= \eps^{1/2}  a_K^w \underbrace{ \inner{a_K^w}^{-1} r^w  \comii
v^{-(s+\gamma/2)}\com{\inner{a_K^{1/2}}^w}^{-1} }_{\in \mathcal B(L^2)} \eps^{-1/2}\inner{a_K^{1/2}}^w \comii v^{s+\gamma/2}.  \]
This gives
\begin{eqnarray*}
     \abs{\inner{r^w f,~f}_{L^2}} &\lesssim& \eps \norm{
  a_K^wf}_{L^2}^2+\eps^{-1}
\norm{  \inner{a_K^{1/2}}^w \comii v^{s+\gamma/2} f}_{L^2}^2\\
&\lesssim& \eps \norm{
  a_K^wf}_{L^2}^2+\eps^{-1} \inner{
\norm{\hat P_K f}_{L^2}^2+\norm{f}_{L_\ell^2}^2},
\end{eqnarray*}
the last inequality following from \reff{coest}.  Taking into account \reff{lassharp},  one has
\begin{eqnarray*}
{\rm Re} \inner{\inner{\tilde a_K a_K}^w f,~f}_{L^2}
\lesssim {\rm Re}
\inner{  a_K ^w f,~\tilde a_Kf}_{L^2}+\eps \norm{
  a_K^wf}_{L^2}^2+\eps^{-2} \inner{
\norm{\hat P_K f}_{L^2}^2+\norm{f}_{L_\ell^2}^2},
\end{eqnarray*}
which along with \reff{tildecoer} yields
\[
     \norm{
  \inner{\tilde a_K^{1/2}a_K^{1/2}}^wf}_{L^2}^2 \lesssim   {\rm Re}
\inner{  a_K ^w f,~\tilde a_Kf}_{L^2}+\eps \norm{
  a_K^wf}_{L^2}^2+\eps^{-2} \inner{
\norm{\hat P_K f}_{L^2}^2+\norm{f}_{L_\ell^2}^2}.
\]
Moreover note that
\[
           \norm{a_K^wf}_{L^2}^2\lesssim  \norm{
  \inner{\tilde a_K^{1/2}a_K^{1/2}}^wf}_{L^2}^2
\]
due to the conclusion (iii) in Lemma \ref{inverse}. Then the desired
estimate
\reff{aakw} follows from the above inequalities,  completing the proof
of Lemma  \ref{regvel}.
\fin

Combining  the conclusions in Lemma \ref{crucial} and Lemma \ref{regvel},  we
obtain   Proposition
\ref{pkreg}. Thus  Theorem \ref{thmain}  follows due to Proposition \ref{prpmain}.    Now
it remains to do the

\preuve[of Theorem \ref{th1}]
Let $\tau$ be the dual variable of $t$ and let $\hat P_\tau$ be the
operator defines as follows
\[
  \hat P_\tau=i \tau+ i v\cdot\xi-\mathcal L=i\inner{\tau+v\cdot\xi}+a^w+\mathcal K .
\]
Just proceeding as in the proof of  Lemma \ref{crucial} and
Lemma \ref{regvel},  we have the maximal hypoelliptic estimate
\begin{equation}\label{+++11052931}
\norm{  \comii v^{2s+\gamma} f}_{L^2}+\norm{  \inner{\tau+v\cdot\xi}
  f}_{L^2}+\norm{   a^w f}_{L^2}+\norm{\comii v^{\frac{\gamma}{1+2s}}\abs{\xi}^{\frac{2s}{1+2s}} f}_{L^2}\lesssim
\norm{ \hat P_\tau  f}_{L^2}+\norm{ f}_{L_\ell^2}.
\end{equation}
Now it remains to prove
\begin{eqnarray*}
   \norm{\comii v^{\frac{\gamma-2s}{1+2s}}\comii\tau^{\frac{2s}{1+2s}}f}_{L^2}
\lesssim \norm{\hat P_\tau  f}_{L^2}+\norm{ f}_{L_\ell^2}.
\end{eqnarray*}
To do so,   we compute
 \begin{eqnarray*}
   \comii
   v^{\frac{\gamma -2s}{1+2s}}\abs\tau^{\frac{2s}{1+2s}}
   &\lesssim&  \comii v^{\frac{\gamma -2s}{1+2s}}\abs{\tau+v\cdot\xi}^{\frac{2s}
{1+2s}}
   +\comii v^{\frac{\gamma-2s}{1+2s}}\abs{v\cdot\xi}
^{\frac{2s}{1+2s}}\\
   &\lesssim & \comii v^{\frac{\gamma -2s}{1+2s}} \abs{\tau+v\cdot\xi}^{\frac{2s}
{1+2s}}  +\comii v^{\frac{\gamma}{1+2s}}\abs{\xi}^{\frac{2s}{1+2s}}\\
   &\lesssim & \comii v^{\gamma-2s} + \abs{\tau+v\cdot\xi}
   +\comii v^{\frac{\gamma}{1+2s}}\abs{\xi}^{\frac{2s}{1+2s}},
 \end{eqnarray*}
 where the last inequality follows from  the Young's inequality:
 \[
     \comii v^{\frac{\gamma-2s}{1+2s}} \abs{\tau+v\cdot\xi}^{\frac{2s}{1+2s}} \leq
\frac{\inner{\comii v^{\frac{\gamma-2s}{1+2s}}}^{1+2s}}{1+2s}+\frac{2s}{1+2s}
\biginner{\abs{\tau+v\cdot\xi}^{\frac{2s}{1+2s}}}^{(1+2s)/(2s)}.
 \]
 As a result  we have,
 \begin{eqnarray*}
   \norm{\comii v^{\frac{\gamma-2s}{1+2s}}\abs
   \tau^{\frac{2s}{1+2s}}f}_{L^2}
   & \lesssim& \norm{ \inner{\tau+v\cdot\xi} f}_{L^2} +
   \norm{ \comii v^{\gamma-2s} f}_{L^2}
    +\norm{\comii v^{\frac{\gamma}{1+2s}}\abs{\xi}^{\frac{2s}{1+2s}}
      f}_{L^2} \\
   & \lesssim& \norm{\hat P_\tau  f}_{L^2}+\norm{ f}_{L_\ell^2} ,
\end{eqnarray*}
where the last inequality follows from   \reff{+++11052931}.
The proof of Theorem  \ref{th1} is complete.
\fin

\appendix

\section{Appendix}
In this section we briefly review some tools used through the proofs. The first section is devoted to the links between some integrals concerning the Boltzmann kernel. In the second one, we recall some basic facts about the Weyl quantization and Weyl-H\"ormander calculus, and the last section will recall some  ideas and results about the Wick quantization.

\subsection{Schur's Lemma}
Let  $\bm K$ be an operator whose corresponding integral kernel $(y,z)\rightarrow k(y,z)$ satisfies
\begin{eqnarray*}
	M_1=\sup_{y\in\mathbb R^d} \int_{\mathbb R^d} \abs{k(y,z)}dz<+\infty, \\
	M_2=\sup_{z\in\mathbb R^d} \int_{\mathbb R^d} \abs{k(y,z)}dy<+\infty.
\end{eqnarray*}
Then  $\bm K$ can be extended from $C_0^\infty(\mathbb R^d)$ to a linear continuous operator on  $L^2(\mathbb R^d)$ (still denoted by $\bm K$) whose norm satisfies
\begin{eqnarray*}
	\norm{\bm K}\leq \sqrt{M_1M_2}.
\end{eqnarray*}

\subsection{Integral representations}

\subsubsection{Principal values}
Let  $q(\theta)$ be  a given measurable   function such that
\[\int_\mathbb R \abs{q(\theta) }d\theta d\theta= \infty,\quad \int_\mathbb R \theta^2 \abs{q(\theta)} d\theta<\infty.\]
Then for any $\psi(\theta)\in C^2(\mathbb R)$,  the function
\[
\theta\longrightarrow q(\theta)\inner{\psi(\theta)+\psi(-\theta)-2\psi(0)}
\]
belongs to $L^1$ locally.   In particular, when $q(\theta)$ is moreover an even and compactly supported function,  we use the notation
\begin{eqnarray*}
\int_{\mathbb R} q(\theta) \psi(\theta)d\theta\stackrel{{\rm def}}{=}\frac{1}{2}\int_{\mathbb R} q(\theta) \inner{\psi(\theta)+\psi(-\theta)-2\psi(0)} d\theta.
\end{eqnarray*}
In our paper, we use it for the function $q(\theta)=\abs\theta^{-1-2s}{\bf 1}_{\abs\theta\leq \pi/2}$.

\subsubsection{A basic formula}
The first tool we use
is the following Fubini-type formula, derived by rather explicit computation.

Consider a measurable function $0 \leq F(\alpha, h)$ of variables  $h$ and $\alpha \in \R^3$. For any $h \in \R^3$, we denote by $E_{0,h}$ the (hyper-)vector plane orthogonal to $h$. Then
\begin{equation} \label{interversion}
\int_{\R^3_h} dh \int_{E_{0,h}} d\alpha F(\alpha, h) = \int_{\R^3_\alpha} d\alpha \int_{E_{0,\alpha}} d h \frac{|h|}{|\alpha|} F(\alpha, h).
\end{equation}

\subsubsection{Carleman representation}
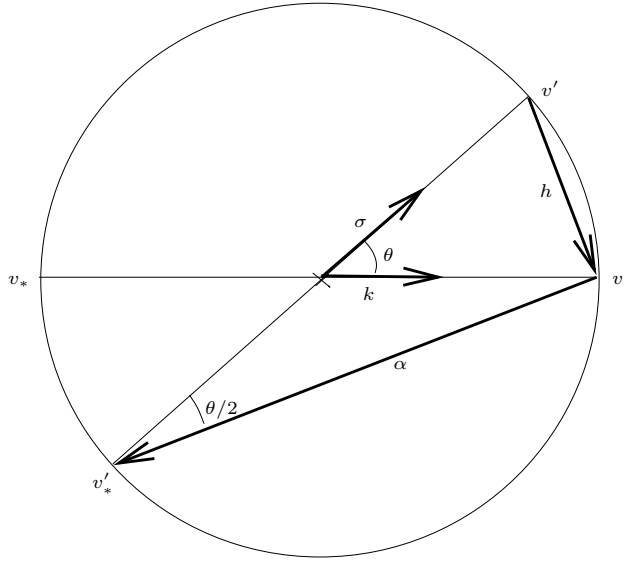
\begin{figure}[!t]
\centering
\begin{picture}(0,0)%
\includegraphics{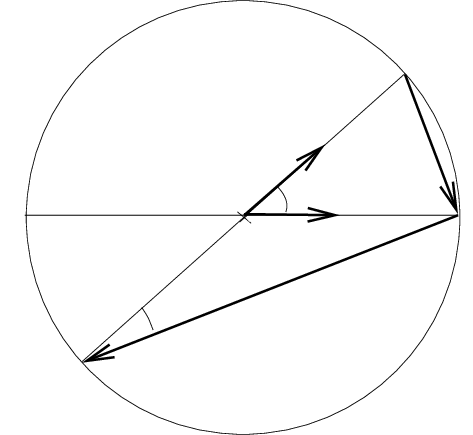}%
\end{picture}%
%
%
\setlength{\unitlength}{2368sp}%
\begingroup\makeatletter\ifx\SetFigFont\undefined%
\gdef\SetFigFont#1#2#3#4#5{%
  \reset@font\fontsize{#1}{#2pt}%
  \fontfamily{#3}\fontseries{#4}\fontshape{#5}%
  \selectfont}%
\fi\endgroup%
\begin{picture}(6285,5798)(1291,-6410)
\put(5191,-3316){\makebox(0,0)[lb]{\smash{{\SetFigFont{7}{8.4}{\rmdefault}{\mddefault}{\updefault}{\color[rgb]{0,0,0}$\theta$}%
}}}}
\put(6826,-1591){\makebox(0,0)[lb]{\smash{{\SetFigFont{7}{8.4}{\rmdefault}{\mddefault}{\updefault}{\color[rgb]{0,0,0}$v'$}%
}}}}
\put(7561,-3556){\makebox(0,0)[lb]{\smash{{\SetFigFont{7}{8.4}{\rmdefault}{\mddefault}{\updefault}{\color[rgb]{0,0,0}$v$}%
}}}}
\put(4981,-3721){\makebox(0,0)[lb]{\smash{{\SetFigFont{7}{8.4}{\rmdefault}{\mddefault}{\updefault}{\color[rgb]{0,0,0}$k$}%
}}}}
\put(4891,-2956){\makebox(0,0)[lb]{\smash{{\SetFigFont{7}{8.4}{\rmdefault}{\mddefault}{\updefault}{\color[rgb]{0,0,0}$\sigma$}%
}}}}
\put(1306,-3556){\makebox(0,0)[lb]{\smash{{\SetFigFont{7}{8.4}{\rmdefault}{\mddefault}{\updefault}{\color[rgb]{0,0,0}$v_*$}%
}}}}
\put(2176,-5686){\makebox(0,0)[lb]{\smash{{\SetFigFont{7}{8.4}{\rmdefault}{\mddefault}{\updefault}{\color[rgb]{0,0,0}$v_*'$}%
}}}}
\put(6811,-2641){\makebox(0,0)[lb]{\smash{{\SetFigFont{7}{8.4}{\rmdefault}{\mddefault}{\updefault}{\color[rgb]{0,0,0}$h$}%
}}}}
\put(5296,-4441){\makebox(0,0)[lb]{\smash{{\SetFigFont{7}{8.4}{\rmdefault}{\mddefault}{\updefault}{\color[rgb]{0,0,0}$\alpha$}%
}}}}
\put(3346,-4906){\makebox(0,0)[lb]{\smash{{\SetFigFont{7}{8.4}{\rmdefault}{\mddefault}{\updefault}{\color[rgb]{0,0,0}$\theta/2$}%
}}}}
\end{picture}%
\caption{$\sigma$ and Carleman representations}
\label{cercle}
\end{figure}

The second formula is the so-called $\omega$-representation. It says that
we have the following (almost everywhere) equalities when all sides are well-defined :

\begin{multline*}
\iint dv_* d\sigma b(\cos \theta) |v-v_*|^\gamma F(v,v_*, v',v'_*)
\\ = 4   \int_{\R^3_h} dh \int_{E_{0,h}} d\alpha \frac{1}{|\alpha +h|\ |h|} b(\cos \theta) |\alpha-h| ^\gamma F(v, v+\alpha-h, v-h, v+\alpha) \\
\\ \approx   \int_{\R^3_h} dh \int_{E_{0,h}} d\alpha \un_{|\alpha| \geq |h|} \frac{1}{|\alpha +h|\ |h|} b \sep{\frac{|\alpha|^2 -|h|^2}{|\alpha+h|^2} } |\alpha+h|^\gamma \\
F(v, v+\alpha-h, v-h, v+\alpha).
\end{multline*}
These formulae are consequences of the following properties (see picture \ref{cercle}):
\begin{enumerate}
\item We make the change of variables $(v_*, \sigma) \longmapsto (\alpha, h)$ with $v'=v-h$, $v_*=v+\alpha-h$, $v_*'= v+\alpha$.
\item Since we restricted by symmetrization to the case  $\sigma\cdot (v-v_*) \geq 0$ (which is equivalent to $\cos \theta \geq 0$), this implies $|\alpha| \geq |h|$. Note also that $h \perp \alpha$ and therefore $|\alpha+h|^2= |\alpha-h|^2 = |\alpha|^2 + |h|^2$.
\item By immediate trigonometric properties we have $\cos\theta = \frac{|\alpha|^2 -|h|^2}{|\alpha+h|^2}$ and $\sin\theta = \frac{2|\alpha| \ |h|}{|\alpha+h|^2}$.
\end{enumerate}
From the singular behavior of the singular kernel we deduce
$$
0 \leq b(\cos \theta) \approx K \theta^{-2-2s} \approx \tilde{K} (\sin\theta)^{-2-2s} \approx \tilde{K} \frac{|\alpha+h|^{4+4s}}{|\alpha|^{2+2s}|h|^{2+2s}} \approx \frac{|\alpha+h|^{2+2s}}{|h|^{2+2s}},
$$
since $|\alpha|^2 \leq |\alpha+h|^2 \leq 2|\alpha|^2$.
At the end we get
\begin{multline} \label{carl}
\iint dv_* d\sigma b(\cos \theta) |v-v_*|^\gamma F(v,v_*, v',v'_*) \\
 =   \int_{\R^3_h} dh \int_{E_{0,h}} d\alpha \tb(\alpha,h) \un_{|\alpha| \geq |h|} \frac{|\alpha+h|^{\gamma+1+2s}}{ |h|^{3+2s}}   F(v, v+\alpha-h, v-h, v+\alpha).
\end{multline}
where $\tb(\alpha, h)$ is bounded from below and above by positive constants, and
$ \tb(\alpha, h) = \tb(\pm \alpha, \pm h)$.
Figure \ref{cercle}  shows the preceding relations between all vectors and angles.

\subsubsection{The cancellation lemma}

We give here an other formula,  in a slightly different version than the original one presented in \cite{AV02}. We consider a function
$G(|v-v_*|, |v-v'|)$. Then for smooth $f$, we have
$$
 \sep{ \iint dv_* d\sigma  G(|v-v_*|, |v-v'|) b(\cos \theta)   \sep{ f'_* -  f_* }}  = S\ast_{v_*} f (v),
 $$
where for all $z \in \R^3$, $S$ has the following expression
\begin{equation*}
\begin{split}
S (z) = &  2\pi \int^{\pi /2}_0 d\theta \sin \theta  b(\cos \theta ) \sep{ G \bigg( {{|z|}\over{\cos {\theta \over 2}}}, {{|z|}\over{\cos {\theta \over 2}}} \sin {\theta \over 2} \bigg) \cos^{-3 } {\theta \over 2} -G (|z|, | z| \sin {\theta \over 2})}
\end{split}
\end{equation*}
This applies in particular to functions of type
 $$
 G(|v-v_*|, |v-v'|, \cos \theta) = b(\cos \theta) |v-v_*|^\gamma \phi(v-v').
 $$

\subsection{Weyl-H\"ormander calculus}\label{subsec41}

We recall here  some  notations and basic facts  of  symbolic
calculus, and refer to
\cite[Chapter 18]{Hormander85} or \cite{MR2599384}  for detailed discussions on the pseudo-differential
calculus.

From now on,  we set $\Gamma = \abs{dv}^2+\abs{d\eta}^2$,
and let $M$ be an admissible weight function w.r.t. $\Gamma$, that is the weight function $M$ satisfies the following conditions:
\begin{enumerate}[(a)]
\item (slowly varying condition) there exists a constant $\delta$ such that
\begin{eqnarray*}
\forall X, Y ~\abs{X-Y}\leq \delta, \quad M(X) \approx M(Y);
\end{eqnarray*}
\item (temperance)  there exist two constants $C$ and $N$ such that
\begin{eqnarray*}
\forall~ X, Y\in \mathbb R^6, \quad M(X)/M(Y) \leq C
\comii{X-Y}^N.
\end{eqnarray*}
\end{enumerate}
Considering  symbols $q(\xi, v,\eta)$ as a function of $(v,\eta)$ with
parameters $\xi$,    we say that  $q\in S\inner{M,\Gamma}$ uniformly
with respect to $\xi$, if
\[
    \forall~ \alpha, \beta\in\mathbb Z_+^3,~~\forall~v,\eta\in\mathbb R^3,\quad
\abs{\partial_v^\alpha\partial_\eta^\beta q(\xi,v,\eta)}\leq C_{\alpha,
\beta} M,
\]
with $ C_{\alpha,\beta}$ a constant depending only on $\alpha$ and $\beta$, but
independent of $\xi$.
For simplicity of notations, in the following discussion, we
omit the parameters dependence in the  symbols, and by $q\in S(M,\Gamma)$ we always mean
that $q$  satisfies the above inequality, uniformly with respect to
$\xi$. The space $S(M,\Gamma)$ endowed with the semi-norms
\begin{eqnarray}\label{seminorm}
 \norm{q}_{k; S(M, \Gamma)}= \max_{0 \leq \abs\alpha+\abs\beta \leq k} \sup_{(v,\eta)\in \mathbb
 R^6} \abs{M(v,\eta)^{-1}\partial_v^\alpha\partial_\eta^\beta q (v,\eta)},
\end{eqnarray}
becomes a Fr\'echet space.
 Let $q\in \mathcal S'(\mathbb R_v^3\times\mathbb R_\eta^3)$ be a tempered distribution and let $t\in\mathbb R$.   the operator ${\rm op}_t q$ is an operator from  $ \mathcal S(\mathbb R_v^3)$ to $ \mathcal S'(\mathbb R_v^3),$  whose Schwartz kernel $K_t$ is defined by the oscillatory integral:
 \[
    K_t  (z, z')= (2\pi)^{-3} \int_{\mathbb R^3} e^{i (z-z') \cdot\zeta}q((1-t)z+tz', \zeta) d\zeta.
  \]
In particular we denote $q(v, D_v)={\rm op}_0q$ and $q^w={\rm op}_{1/2}q$.
 Here $q^w$ is called  the Weyl quantization of symbol $q$.

 An elementary property to be used
frequently is  the $L^2$ continuity  theorem in the class
$S\inner{1,~g}$, see \cite[Theorem 2.5.1]{MR2599384} for instance,  which says that  there exists a constant $C$
and  a positive integer $N$ depending only the dimension,
such that
\begin{equation}\label{bdness}
    \forall ~u\in L^2,\quad \norm{ q^w u}_{L^2}\leq C \norm{q}_{N;
      S(1,\Gamma)}\norm{u}_{L^2}.
\end{equation}
Let us also recall here the composition formula of Weyl
quantization.  Given $p_i\in S(M_i,\Gamma)$ we have
\begin{eqnarray}\label{11053010}
  p_1^wp_2^w=(p_1\sharp p_2)^w
\end{eqnarray}
with $p_1\sharp p_2 \in S\inner{M_1M_2,~\Gamma}$  admitting the expansion
\begin{equation}\label{expan}
   p_1\sharp p_2=p_1 p_2+ \int_0^1\iint e^{- i \sigma(Y-Y_1, Y-Y_2)/(2\theta)}
   \frac{i}{2}\sigma (\partial_{Y_1}, \partial_{Y_2}) p_1(Y_1)
   p_2 (Y_2) dY_1 dY_2 d\theta/\theta^6,
\end{equation}
where $\sigma$ is the symplectic form in $\mathbb R^6$ given by
\[
   \sigma \inner{(z,\zeta), (\tilde z,\tilde \zeta)}=\zeta\cdot
   \tilde z-\tilde \zeta \cdot z.
\]
 For the relation between  the classical pseudo-differential operator $q(v,D_v)$ and Weyl quantization $q^w$,  we  have the formula:
  \begin{equation}\label{11081101}
    q^w =\comi{{J^{1/2} q}}(v,D_v),
  \end{equation}
  where $J^{1/2}: \mathcal S'\rightarrow \mathcal S'$ is defined by
 \begin{equation}\label{11081106}
    ( J^{1/2}q) (v.\eta)=(2\pi)^{-3} \iint  e^{-i z\cdot\zeta}
    q(v+z,\eta+\zeta)dz d\zeta.
  \end{equation}

\subsection{Wick quantization}
Finally let us  recall some  basic properties of the Wick quantization, which is  also called anti-Wick in \cite{shubin}.   The importance in studying the Wick quantization lies in the facts that  positive symbols give rise to positive operators.   There are several equivalent ways of defining Wick quantization and one is defined in terms of coherent states.   The coherent states method  essentially reduces the partial differential operators to  ODEs,  by virtue of  the Wick calculus.  We  refer
the readers to the
works   \cite{Berezin, Berezin1, degosson, MR1957713, MR2477145, MR2599384, shubin} and references therein  for extensive presentations of this
quantization and  its applications in  mathematics and mathematical physics.

 Let $Y = (v, \eta)$  be a point in $
\mathbb R^{6}$.  The Wick quantization of a symbol $q$  is  given by
\begin{eqnarray*}
	q^{{\rm Wick}} =(2\pi)^{-3}\int_{\mathbb R^6} q(Y) \Pi_{Y}\ dY,
\end{eqnarray*}
where  $\Pi_Y$ is the projector associated to the Gaussian $\varphi_Y$ which is defined by
\[
   \varphi_Y(z)=\pi^{-3/4}e^{-\frac{1}{2}\abs{z-v}^2}e^{i z \cdot\eta /2},\quad z\in
\mathbb R^3.
\]
 The main property of the Wick quantization is its positivity, i.e.,
\begin{eqnarray*}
  q(v,\eta)\geq 0 ~~\,\textrm{for all}~ (v,\eta)\in\mathbb R^{6} ~{\rm implies}
~~\,q^{\rm Wick}\geq 0.
\end{eqnarray*}
According to Theorem 24.1 in \cite{shubin}, the Wick and Weyl
quantizations of a symbol $q$ are linked by
the following identities
\begin{eqnarray}\label{ww}
  q^{\rm Wick}=\inner{q* \pi^{-3} e^{-
\abs{\cdot}^2 }}^w=q^w+r^w
\end{eqnarray}
with
\[
  r(Y)= \pi^{-3} \int_0^1 \int_{\mathbb R^{6}} (1-\theta)q''(Y+\theta Z )Z^2e^{-
\abs{Z}^2}\di Zd\theta.
\]
As a result,  $q^{\rm Wick} $ is a bounded operator in $L^2$ if $q\in S(1, g)$ due to \reff{bdness}.

We also recall the following composition
formula obtained in the proof of Proposition 3.4 in \cite{MR1957713}
\begin{eqnarray} \label{11082406}
  q_1^{\rm Wick}q_2^{\rm Wick}=\com{q_1q_2- q_1'\cdot q_2'+
  \frac{1}{i}\big \{q_1,~q_2\big\}}^{\rm Wick}+T,
\end{eqnarray}
with $T$ a bounded operator in $L^2(\mathbb R^{2n})$,  when $q_1\in L^{\infty}
(\mathbb R^{2n})$ and $q_2$ is a smooth symbol whose derivatives of order $\geq2$
are bounded on $\mathbb R^{6}$. The notation $\set{q_1,q_2}$ denotes the Poisson
bracket defined by
\begin{eqnarray} \label{11051505}
  \bigset{q_1,~q_2}=\frac{\partial q_1}{\partial\eta}\cdot\frac{\partial q_2}{\partial
v}-\frac{\partial q_1}{\partial v}\cdot\frac{\partial q_2}{\partial \eta}.
\end{eqnarray}


\end{document}